\documentclass[11pt,dvips]{article}

\usepackage{latexsym,url}
\usepackage{mdwlist}
\usepackage{export} %sert a exporter des compteurs, utile pour mettre un enonce 2 fois.
\usepackage{amsmath,amsthm,amsfonts,amssymb}
\usepackage{rotating} %utilise pour les fleches \actson \actedon
\usepackage{graphicx}
\usepackage{color}
\usepackage{epsfig}
\usepackage[all]{xy}
\usepackage{pstricks,pst-node,pst-ellipse}

\newtheorem{thm}{Theorem}[section]
\newtheorem*{pb*}{Problem}

\newtheorem*{thm*}{Theorem}
\newtheorem*{conj*}{Conjecture}
\newtheorem{dfn}[thm]{Definition} 
\newtheorem*{dfn*}{Definition}

\newtheorem{cor}[thm]{Corollary}
\newtheorem*{cor*}{Corollary}

\newtheorem{prop}[thm]{Proposition} 
\newtheorem*{prop*}{Proposition} 
 
\newtheorem{lem}[thm]{Lemma} 
\newtheorem*{lem*}{Lemma} 
 
\newtheorem{claim}[thm]{Claim} 
\newtheorem*{claim*}{Claim} 
 
\newtheorem*{fact*}{Fact}

\newtheorem*{conv*}{Convention} 
\newtheorem*{propdouble}{Propositions \ref{universal=closed} and \ref{prop_univ2}}

\theoremstyle{remark}
\newtheorem*{rem*}{Remark}
\newtheorem*{rems*}{Remarks}
\newtheorem*{example*}{Example}

\newenvironment{SauveCompteurs}[1]{%
\newcommand{\monparametre}{#1}
\openexport{\monparametre_sauve}
  \Export{thm}\Export{section}\Export{subsection}\Export{subsubsection}
\closeexport}{}

\newenvironment{UtiliseCompteurs}[1]{%
\newcommand{\monparametre}{#1}% utile car en fait #1 n'est pas accessible dans la deuxieme partie de l'environnement.
\openexport{\monparametre_aux}
  \Export{thm}\Export{section}\Export{subsection}\Export{subsubsection}
\closeexport
\PackageInfo{export}{\MessageBreak
Importations from \monparametre_sauve.xpt\MessageBreak}%
\InputIfFileExists{\monparametre_sauve.xpt}{\relax}{\relax}%
\renewcommand{\label}[1]{}%But=eviter les multiply defined. 
}{%
\PackageInfo{export}{\MessageBreak
Importations from \monparametre_aux.xpt\MessageBreak}%
\InputIfFileExists{\monparametre_aux.xpt}{\relax}{\relax}}

 \newcommand{\es}{\emptyset}
\renewcommand{\phi}{\varphi} 
\newcommand{\m} {^{-1}} 
\newcommand{\pmu} {^{\pm 1}} 
\newcommand{\eps} {\varepsilon} 
\let\ts\textstyle

\newcommand{\eneq} {\,\makebox[0pt][l]{\raisebox{0.7ex}{$\scriptstyle
 \neq$}}\raisebox{-0.7ex}{$\scriptstyle =$}\,}
\newcommand {\ra} {\rightarrow}

\newcommand {\onto} {\twoheadrightarrow}
\newcommand {\into} {\hookrightarrow}
\newcommand {\xra} {\xrightarrow}

\newcommand{\eqv} {{\sim}}% les accolades reduisent l'espace autour (il comprend qu'il n'est pas une relation binaire...
\newcommand{\ultra}[1]{{{}^*#1}}% les accolades reduisent l'espace autour (il comprend qu'il n'est pas une relation binaire...

\newcommand{\ol}[1]{\overline{#1}}

\newcommand{\dunion}{\sqcup}

\newcommand{\ie} {i.~e.\ }

\newcommand {\calc} {{\mathcal {C}}}   
\newcommand {\cald} {{\mathcal {D}}}

\newcommand {\calg} {{\mathcal {G}}}

\newcommand {\call} {{\mathcal {L}}}

\newcommand {\calp} {{\mathcal {P}}}

\newcommand {\calt} {{\mathcal {T}}}

\newcommand{\bbZ} {{\mathbb{Z}}}   
\newcommand{\bbN} {{\mathbb{N}}}   
\newcommand{\bbR} {{\mathbb{R}}}   
   
\newcommand{\bbC} {{\mathbb{C}}}

\newcommand{\Fix}{\mathop{\mathrm{Fix}}}

\newcommand{\Aut} {\mathop{\mathrm{Aut}}}

\newcommand{\Hom} {\mathop{\mathrm{Hom}}}
\newcommand{\Epi} {\mathop{\mathrm{Epi}}}
\newcommand{\Mod} {\mathop{\mathrm{Mod}}}
\newcommand{\can} {\mathrm{can}}

\newcommand{\Id} {\mathrm{Id}}
\newcommand{\id} {\mathrm{id}}
\newcommand{\IGD}{\mathcal{IGD}}

\newcommand{\et}{\wedge}
\newcommand{\ou}{\vee}
\newcommand{\non}{\neg}
\newcommand{\implq}{\rightarrow}
\newcommand{\satis}{\models}
\newcommand{\Elem}{\mathrm{Elem}}
\newcommand{\Univ}{\mathrm{Univ}}
\newcommand{\Exist}{\mathrm{Exist}}

\newcounter{monparact}[subsection]
\renewcommand{\themonparact}{(\alph{monparact})}
\newcommand{\monpara}[1]{\refstepcounter{monparact}\paragraph{\themonparact. #1}}

\begin{document}

\title{Limit groups as limits of free groups:\\ compactifying the set of free groups}
\author{Christophe Champetier, Vincent Guirardel}
%\date{\small \today. Fichier \texttt{\jobname.tex}}
\date{\today}
\maketitle

\begin{abstract}
  We give a topological framework for the study of Sela's \emph{limit groups}:
limit groups are limits of free groups in a compact space of marked
groups.
Many results get a natural interpretation in this setting.
The class of limit groups is known to coincide with the class of
finitely generated fully residually free groups.
The topological approach gives some new insight on the relation between
fully residually free groups, the universal theory of free groups,
ultraproducts and non-standard free groups.
\end{abstract}

\tableofcontents

\section{Introduction}
Limit groups have been introduced by Z. Sela in the first paper of his
solution of Tarski's problem \cite{Sela_diophantine1}. These groups appeared
to coincide with the long-studied class of finitely generated fully
residually free groups
(definition~\ref{residual_freeness}, see \cite{Baumslag_residually}, \cite{Baumslag_generalised},
\cite{KhMy_irreducible1,KhMy_irreducible2}, see
\cite{BMR_algebraic_in_algorithmic} and
\cite{Chi_introduction} and references).
In this paper, we propose a new approach of limit groups,
in a topological framework, that sheds further light on these groups.
We survey the equivalent definitions of limit groups and
their elementary properties, and we detail the Makanin-Razborov diagram,
and general ways of constructing limit groups.
This article is aimed to be
self-contained, and some short classical proofs are rewritten for completeness.
\\

\paragraph{Tarski's problem.}
Tarski's problem asks whether all the free groups of rank $\geq 2$
have the same elementary theory.  The elementary theory of a group $G$
is the set of all sentences satisfied in $G$ (see section \ref{sec_logic}
for a short introduction, and
\cite{Chang_Keisler_model,Hodges_shorter} for further references).
Roughly speaking, a \emph{sentence} in the language of groups is a
``usual'' logical sentence where one quantifies only on individual elements of a group;
to be a little bit more precise, it is a string of symbols made of
quantifiers, variables (to be interpreted as elements of a group), 
the identity element ``$1$'', the group
multiplication and inverse symbols (``.'' and ``$^{-1}$''), equality
``$=$'', and logical connectives ``$\non$'' (not), ``$\wedge$'' (and),
``$\vee$'' (or) and without free variables.
Any sentence is equivalent (assuming the axioms of groups) to a sentence where all quantifiers are
placed at the beginning, followed by a disjunction of systems of
equations or inequations:

$$\forall x_1,\dots, x_p\, \exists y_1,\dots, y_q \,\forall z_1,\dots,
z_r \dots \left\{
\begin{matrix}
  u_1(x_i,y_j,z_k,\dots ) \eneq 1 \\ \vdots \\ u_n (x_i,y_j,z_k,\dots
  )\eneq 1
\end{matrix} \right.
\vee \left\{
\begin{matrix}
  v_1(x_i,y_j,z_k,\dots ) \eneq 1 \\ \vdots \\ v_m (x_i,y_j,z_k,\dots
  )\eneq 1
\end{matrix} \right.
\vee \dots
$$

Such a sentence is a \emph{universal sentence} if it can be written
$\forall x_1,\dots, x_p\ \phi(x_1,\dots ,x_p)$
for some quantifier free formula $\phi(x_1,\dots,x_p)$.
For example, $\forall x, y ~ xy=yx$
is a universal sentence in the language of groups. This
sentence is satisfied in a group $G$ if and only if $G$ is abelian.
The universal theory $\Univ(G)$ of a group $G$ is the set of 
universal sentences satisfied by $G$.

\paragraph{From equations to marked groups.}
The first step in the study of the elementary theory of a group is the
study of systems of equations (without constant) in that group.
Such a system is written
$$\left\{\begin{matrix} w_1(x_1,\dots
    ,x_n) =1 \\ \vdots \\ w_p (x_1,\dots ,x_n)= 1
\end{matrix} \right.$$ 
where each $w_i(x_1,\dots ,x_n)$ is a reduced word on
the variables $x_1,\dots ,x_n$ and their inverses, \ie an element of
the free group $F_n = \langle x_1,\dots ,x_n\rangle$.

Solving this system of equations in a group $G$ consists in finding all tuples
$(a_1,\dots ,a_n)\in G^n$ such that for all index $i$, $w_i(a_1,\dots ,a_n) = 1$
in $G$.
There is a natural correspondence between solutions of this system of equations 
and morphisms $h:E\ra G$, where 
$E$ is the group presented by $E=\langle x_1,\dots x_n ~|~
w_1(x_1,\dots ,x_n), \dots , w_p(x_1,\dots ,x_n)\rangle$:
to a solution $(a_1,\dots,a_n)\in G^n$ corresponds the morphism $E\ra G$
sending $x_i$ on $a_i$; and conversely, given a morphism $h:E\ra G$, 
 the corresponding solution is the tuple $(h(x_1),\dots , h(x_n))\in G^n$.

Thus, to study the set of solutions of a system of equations in free groups,
one has to understand the set of morphisms from a finitely presented group $E$
to free groups.
Any morphism $h$ from $E$ to a free group $F$ is obtained by composing an
epimorphism from $E$ onto the free group $h(E)$ with a morphism from
$h(E)$ to $F$, and morphisms from a free group to an arbitrary group
are well-known. The study of all morphisms from $E$ to $F$ thus reduces
to the study of epimorphisms from $E$ onto free
groups (of rank at most $n$).

The very beginning of Sela's study of equations in free groups 
may be viewed as a compactification of the set of all
epimorphisms from $E$ onto free groups.
This compactification consists in epimorphisms from $E$ onto possibly
non-free groups, which are called \emph{limit groups} by Sela.

An epimorphism $h$ from the group $E = \langle x_1,\dots x_n ~|~
w_1(x_1,\dots ,x_n), \dots , w_p(x_1,\dots ,x_n)\rangle$ onto a group
$G$ gives a preferred generating family $(h(x_1),\dots , h(x_n))$ of
$G$. In other words, $h$ defines a marking on $G$ in the following sense:
a \emph{marked group} $(G,S)$ is a group $G$ together with an ordered generating family $S$.
There is a natural topology on the set of marked groups (topology of Gromov-Hausdorff and Chabauty)
which makes it compact.
This topology can be roughly described as follows:
two marked groups $(G,S)$ and $(G',S')$ are closed to each other
if large balls of their Cayley graphs are isomorphic.
Spaces of marked groups have been used in \cite{Gromov_polynomial}, 
\cite{Grigorchuk_degrees}, and studied in \cite{Champetier_espace}. 
Some elementary properties of this topology will be presented in section \ref{sec_topology}.

\paragraph{Limit groups.}
In this paper, we propose a new definition of limit groups as limit of marked free groups.
We give five equivalent characterizations of limit groups. Three
of them are well known in model theory (see \cite{Remeslennikov_exist_siberian}).

\begin{thm}\label{thm_equiv}
Let $G$ be a finitely generated group.
The following assertions are equivalent:
\begin{enumerate*}
\item $G$ is a limit group in the sense of Sela (\cite{Sela_diophantine1}) \label{limit_sela}

\item\label{marked}
Some marking (or equivalently any marking) of $G$ is a limit
of markings of free groups in a
compact space of marked groups.
  
\item\label{universal}
$G$ has the same universal theory as a free group.
  
\item\label{non_standard} $G$ is a subgroup of a non standard free group.
  
\item\label{fully} $G$ is fully residually free.
\end{enumerate*}
\end{thm}

\begin{rem*}
Fully residually free groups have another interpretation in the
language of algebraic geometry over free groups : they are precisely
the coordinate groups of irreducible algebraic sets in free groups (see
\cite{BMR_algebraic_in_algorithmic} Theorem D2 and \cite{KhMy_description}
 Lemma 4) 
\end{rem*}

The proof of this theorem will follow from our propositions
\ref{limit_residually}, \ref{limit_universal}, \ref{limit_ultra} and
\ref{ultra_residually}. 
The equivalence between \ref{limit_sela} and \ref{fully}
is due to Sela (\cite{Sela_diophantine1}).
The equivalence between \ref{universal} and
\ref{non_standard} is a particular case of more general results in
model theory (see for example \cite{Bell_Slomson_models} Lemma 3.8 Chap.9).
The equivalence between \ref{universal},
\ref{non_standard} and \ref{fully} is shown in
Remeslennikov \cite{Remeslennikov_exist_siberian}
(groups that have the same universal theory as free groups are called
\emph{universally free groups} in \cite{FGMRS_n-free}
and \emph{$\exists$-free groups} by Remeslennikov in
\cite{Remeslennikov_exist_siberian}).
Assertion \ref{marked} is a reformulation of lemma 1.3(iv) in \cite{Sela_diophantine1}.
The topological point of view allows us to give direct proofs for the
equivalences between \ref{marked}, \ref{universal}, and
\ref{non_standard} (see sections \ref{sec_universal} and 
\ref{ultraproduct_topology}). 
More generally, we relate the topology on the set of marked groups
to the universal theory of groups and to ultraproducts by the two following propositions
(see sections \ref{sec_universal}, \ref{ultraproduct_topology}):

\begin{UtiliseCompteurs}{cpteursUniv}
\begin{propdouble}
If $\Univ(G)\supset\Univ(H)$, then for all generating family $S$ of $G$,
$(G,S)$ is a limit of marked subgroups of $H$.
Moreover, if a sequence of marked groups $(G_i,S_i)$ converge to a marked group $(G,S)$, then
$\Univ (G)\supset \limsup\Univ (G_i).$
\end{propdouble}
\end{UtiliseCompteurs}

\begin{UtiliseCompteurs}{cpteursUltra}
\begin{prop}
The limit of a converging sequence of marked groups $(G_i,S_i)$ embeds in any (non-principal)
ultraproduct of the $G_i$'s. Moreover, any finitely generated subgroup of an ultraproduct of $G_i$
is a limit of a sequence of markings of finitely generated subgroups of the $G_i$'s.
\end{prop}
\end{UtiliseCompteurs}

\paragraph{Simple properties of limit groups.}
The topological point of view on limit groups gives natural proofs
of the following
(well known) simple properties of limit groups:
\newcommand{\elementaryProperties}{%
  \begin{prop}
    \label{elementary_properties}
    Limit groups satisfy the following properties:
    \begin{enumerate*}
    \item\label{item_CSA}
      A limit group is torsion-free, commutative transitive, and CSA (see Definitions
  \ref{dfn_commutative} and \ref{dfn_CSA}).
    \item \label{subgroups} Any finitely generated subgroup of a limit group is a limit group.
    \item \label{betti} If a limit group is non-trivial
      (resp. non-abelian),
      then its first Betti number is at least $1$ (resp.\ at least 2).
    \item \label{2_generated}
      Two elements of a limit group generate a free abelian group
      ($\{1\}$, $\bbZ$ or $\bbZ ^2$) or a non-abelian free group of rank $2$.
    \item \label{item_ordered}
      A limit group $G$ is bi-orderable: there is a total order on $G$ which is left and right invariant.
    \end{enumerate*}
  \end{prop}
}  
\begin{UtiliseCompteurs}{elementary_properties}
\elementaryProperties
\end{UtiliseCompteurs}

\paragraph{First examples of limit groups.}
Limit groups have been extensively studied as finitely generated fully
residually free groups (\cite{Baumslag_generalised},
\cite{Baumslag_residually}, and \cite{Chi_introduction} and its references).
The class of fully residually free groups
is clearly closed under taking subgroup and free products.
The first non free finitely generated examples of fully
residually free groups, including all the non-exceptional surface groups,
have been given by Gilbert and Benjamin Baumslag in
\cite{Baumslag_generalised} and \cite{Baumslag_residually}.
They obtained fully residually free groups 
by \emph{free extension of centralizers} in free groups (see section \ref{sec_construction}).
A free extension of centralizers of a limit group $G$ is a group of the form $G*_C (C\times \bbZ^p)$
where $C$ is a maximal abelian subgroup of $G$.
A free (rank $p$) extension of centralizers of a limit group $G$ is a group of
the form $G*_C (C\times \bbZ^p)$ where $C$ is a maximal abelian
subgroup of $G$.
As a corollary, the fundamental group of a closed surface with Euler
characteristic at most $-2$ is a limit group:
indeed, such a group embeds in a double of a free group over a maximal cyclic subgroup, 
and such a double occurs as a subgroup of an extension of centralizers
of a free group (see section \ref{sec_construction}). 

As seen above, fully residually free groups of rank at most 2 are known to be either one of the free
abelian groups $1$, $\bbZ$, $\bbZ ^2$, or the non-abelian free group $F_2$
(see \cite{Baumslag_generalised}).
In \cite{FGMRS_classification}, a classification of
$3$-generated limit group is given:
a $3$-generated limit group is either a free group of rank $3$, a free
abelian group of rank $3$ or a free rank one extension of centralizer in a free
group of rank $2$.

\paragraph{Finiteness properties.}
The works of Kharlampovich and Myasnikov (\cite{KhMy_irreducible1,KhMy_irreducible2}),
and of Sela \cite{Sela_diophantine1} show that limit groups 
can be obtained recursively from 
free groups, surface groups and free abelian groups by a finite
sequence of free products or amalgamations over $\bbZ$ (see also \cite{Gui_limit}
where another proof is presented using actions of groups on
$\bbR^n$-trees).
Such a decomposition of a limit group implies its finite presentation.
In the topological context, the finite presentation of limit groups
allows to give short proofs of the two following finiteness results.

\begin{UtiliseCompteurs}{cpteur_Finitude}
\begin{prop}[\cite{BMR_algebraic_in_algorithmic} Corollary 19, \cite{Sela_diophantine1}]
Given a finitely generated group $E$, there
exists a finite set of epimorphisms $E\onto G_1$,\dots ,$E\onto G_p$ from $G$
to limit groups $G_1,\dots ,G_p$ such that any morphism from $E$
to a free group factorizes through one of these epimorphisms.
\end{prop}
\end{UtiliseCompteurs}

\begin{UtiliseCompteurs}{cpteur_sequence_term}
\begin{prop}[\cite{Razborov_systems}, \cite{KhMy_irreducible1}, \cite{Sela_diophantine1}]
  Consider a sequence of quotients of limit groups 
$$G_1\onto G_2\onto\dots\onto G_k\onto\dots$$
Then all but finitely many epimorphisms are isomorphisms.
\end{prop}  
\end{UtiliseCompteurs}

In  sections \ref{sec_logic} and \ref{sec_finitude}, we give alternative proofs
avoiding the use of the finite presentation of limit groups. 
The first one is an argument by Remeslennikov (\cite{Remeslennikov_exist_siberian}).
The second result is due to Razborov (\cite{Razborov_systems}, see also \cite{KhMy_irreducible2,BMR_algebraicI})
and the simple argument we give is inspired by a point of view given in \cite{Chatzidakis_limit}.

These two finiteness results are important steps in the construction of Makanin-Razborov diagrams.

\paragraph{Makanin-Razborov diagrams.}
To understand the set of solutions a given system of equations $(S)$ in free groups, 
or equivalently, the set of morphisms from a fixed group $E$ into free groups,
Sela introduces a Makanin-Razborov diagram associated to $E$ (or equivalently to $(S)$). 
This diagram is a finite rooted tree
whose root is labelled by $E$. Its other vertices are labelled by quotients of $E$ which are limit groups.
Its essential feature is that
any morphism $h$ from $E$ to a free group can be \emph{read} from this diagram (see section \ref{subsec_MR}).

\paragraph{Construction and characterization of limit groups.}

The first characterization of limit groups is due to Kharlampovich-Myasnikov (\cite{KhMy_irreducible2}).
A finitely generated group is an \emph{iterated extension of centralizers of a free group} if
it can be obtained from a free group by a sequence of free extension of centralizers.
\newcommand{\caracUn}{%
\begin{thm}[First characterization of limit groups {\cite[Th.4]{KhMy_irreducible2}}]\label{thm_sbg_ext}
  A finitely generated group is a limit group if and only if it is a
  subgroup of an iterated extension of centralizers of a free group.
\end{thm}
}
\begin{UtiliseCompteurs}{caracUn}
 \caracUn 
\end{UtiliseCompteurs}

The second characterization we give does not require to pass to a subgroup.
It is defined in terms of what we call \emph{generalized double}, which
is derived from Sela's strict MR-resolutions. 

\newcommand{\dfnGeneDouble}{%
\begin{dfn}[Generalized double]
  A \emph{generalized double over a limit group $L$} is a group 
 $G=A*_C B$ (or $G=A*_C$) such that both vertex groups $A$ and $B$ are finitely generated and
 \begin{enumerate*}
 \item $C$ is a  non-trivial abelian group
whose images under both embeddings are maximal abelian in the vertex groups
\item there is an epimorphism $\phi:G\onto L$ which is one-to-one in restriction to each vertex
group (in particular, each vertex group is a limit group).
 \end{enumerate*}
\end{dfn}
}
\begin{UtiliseCompteurs}{dfnGeneDouble}
 \dfnGeneDouble 
\end{UtiliseCompteurs}

The terminology comes from the fact that a genuine double $G=A*_{C=\ol C}\ol A$ over a maximal abelian group
is an example of a generalized double where $L=A$ and $\phi:G\ra A$ is the morphism restricting to the identity on $A$
and to the natural map $\ol x\mapsto x$ on $\ol A$.

\begin{UtiliseCompteurs}{caracDeux}
\begin{thm}[Second characterization of limit groups] (compare Sela's strict MR-resolution).
  The class of limit groups coincides with the class $\IGD$ defined as the smallest class 
containing finitely generated free groups,
and stable under free products and under generalized double over a group in $\IGD$.
\end{thm}
\end{UtiliseCompteurs}

\paragraph{Other constructions of limit groups.}
There are more general ways to construct limit groups in the flavor of
generalized double, and closer to Sela's strict MR-resolution. Our  most
general statement is given in Proposition
\ref{prop_general_to_simple}.
A slightly simpler statement is the following proposition.

\begin{UtiliseCompteurs}{prop_sgolg_subice}
\begin{prop}[{\cite[Th.5.12]{Sela_diophantine1}}]
  Assume that $G$ is the fundamental group of a graph of groups
  $\Gamma$ with finitely generated vertex groups such that:
\begin{itemize*}
\item each edge group is a non-trivial abelian group whose images
  under both edge morphisms are maximal abelian subgroups of the
  corresponding vertex groups;
\item $G$ is commutative transitive;
\item there is an epimorphism $\phi$ from $G$ onto a limit group $L$
  such that $\phi$ is one-to-one in restriction to each vertex group.
\end{itemize*}
Then $G$ is a limit group.
\end{prop}
\end{UtiliseCompteurs}

Another version of this result has a corollary which is worth noticing
(Proposition \ref{prop_sgolg_twist} and \ref{prop_cut_surf}).
Consider a marked surface group $(G,S)$ of Euler characteristic at most $-1$.
The modular group of $G$ acts on the set of marked quotients of $(G,S)$ as follows:
given $\phi:G\onto H$ and $\tau$ a modular automorphism of $G$,
$(H,\phi(S)).\tau=(H,\phi\circ\tau(S))$.
For any marked quotient $(H,\phi(S))$ of $(G,S)$ such that $H$ is a non-abelian limit group,
the orbit of $(H,\phi(S))$ under the modular group of $G$
accumulates on $(G,S)$. In other words, one has the following:

%\begin{UtiliseCompteurs}{orbits_in_surfaces}
\begin{cor*}
  Let $G$ be the fundamental group of a closed surface $\Sigma$ with
  Euler characteristic at most $-1$.  Let $\phi$ be any morphism from
  $G$ onto a non-abelian limit group $L$.
  
  Then there exists a sequence of elements $\alpha_i$ in the modular
  group of $\Sigma$ such that $\phi\circ\alpha_i$ converges to $\id_G$
  in $\calg(G)$.
  \end{cor*}
%\end{UtiliseCompteurs}

\paragraph{Fully residually free towers.}
Particular examples of limit groups are the \emph{fully residually
  free towers} (in the terminology of Z. Sela, see also
  \cite{KhMy_irreducible1} where towers appear as coordinate groups of
  particular systems of equations), \ie the class of groups $\calt$, containing all the
finitely generated free groups and surface groups and stable under the
following operations :
\begin{itemize*}
\item free products of finitely many elements of $\calt$.
\item free extension of centralizers.
\item glue on a base group $L$ in $\calt$ a surface group that
  retracts onto $L$ (see section \ref{sec_frft} for the precise
  meaning of this operation).
\end{itemize*}

\begin{UtiliseCompteurs}{cpteur_frftowers}
\begin{thm}[\cite{Sela_diophantine1}]
A  fully residually free tower is a limit group.
\end{thm}
\end{UtiliseCompteurs}

It follows from Bestvina-Feighn combination theorem that a fully residually free tower is
Gromov-hyperbolic if and only if it is constructed without using
extension of centralizers. 

A positive answer to Tarski's problem has been announced in
\cite{Sela_diophantine6} and \cite{KhMy_Tarski}.
We state the result as a conjecture since the referring
processes are not yet completed.

\begin{conj*}[{\cite[Th.7]{Sela_diophantine6}}]
  A finitely generated group is elementary equivalent to a non abelian
  free group if and only if it is a non elementary hyperbolic fully
  residually free tower.
\end{conj*}

We would like to warmly thank Frederic Paulin for his great encouragements and his
careful reading of the paper. We also thank O.~Kharlampovich and A.~Myasnikov for their
very instructive remarks and suggestions.

%%%%%%%%%%%%%%%%%%%%%%%%%%%%%%%%%%%%%%%%%%%%%%%%%%%%%%%%%%%%%%%%%%%%%
%
%%%%%%%%%%%%%%%%%%%%%%%%%%%%%%%%%%%%%%%%%%%%%%%%%%%%%%%%%%%%%%%%%%%%%
%
%%%%%%%%%%%%%%%%%%%%%%%%%%%%%%%%%%%%%%%%%%%%%%%%%%%%%%%%%%%%%%%%%%%%%
%
%%%%%%%%%%%%%%%%%%%%%%%%%%%%%%%%%%%%%%%%%%%%%%%%%%%%%%%%%%%%%%%%%%%%%
%
%%%%%%%%%%%%%%%%%%%%%%%%%%%%%%%%%%%%%%%%%%%%%%%%%%%%%%%%%%%%%%%%%%%%%
%
%%%%%%%%%%%%%%%%%%%%%%%%%%%%%%%%%%%%%%%%%%%%%%%%%%%%%%%%%%%%%%%%%%%%%
%

\section{A topology on spaces of marked groups}
\label{sec_topology}

The presentation given here extends \cite{Champetier_espace}.
We give four equivalent definitions for the space of marked groups.

\subsection{Definitions}
\label{definition}

\monpara{Marked groups.}
A {\it marked group} $(G,S)$ consists in a group $G$ with a prescribed
 family $S = (s_1,\dots ,s_n)$ of generators. Note that the family
is ordered, and that repetitions are allowed.

Two marked groups will be identified if they are isomorphic in the natural sense
for marked groups: two marked groups
$(G,(s_1,\dots ,s_n))$ and $(G',(s'_1,\dots ,s'_n))$ are isomorphic as
marked groups if and only if the bijection
that sends $s_i$ on $s'_i$ for all $i$
extends to a isomorphism from $G$ to $G'$ 
(in particular the two generating families must have the same cardinality $n$).

\begin{dfn}
For any fixed $n$, the set of marked groups $\calg _n$ is the set of groups marked by
$n$ elements up to isomorphism of marked groups.
\end{dfn}

\monpara{Cayley graphs.}
A marked group $(G,S)$ has a natural Cayley graph,
whose edges are labeled by integers in $\{ 1,\dots , n\}$.
Two marked groups are isomorphic as marked groups if and only if
their Cayley graphs are isomorphic as labeled graphs.
Thus $\calg _n$ may be viewed as the set of labeled Cayley graphs on
$n$ generators up to isomorphism.

\monpara{Epimorphisms.}
Fix an alphabet $\{s_1,\dots ,s_n\}$ and
consider the free group $F_n = \langle s_1,\dots ,s_n\rangle$ marked
by the free basis $(s_1,\dots ,s_n)$.
Generating families of cardinality $n$ for a group $G$ are in
one-to-one correspondence with epimorphisms from $F_n$ onto $G$.
In this context, two epimorphisms $h_1: F_n\ra G_1$, $h_2: F_n\ra G_2$
correspond to isomorphic marked groups if and only if $h_1$ and $h_2$
are \emph{equivalent} in the following sense: there is an
isomorphism $f:G_1 \ra G_2$ making the diagram below commutative.
$$\xymatrix{ 
 F_n \ar[r]^{h_1} \ar[rd]_{h_2} & G_1\ar[d]^{f}\\
&G_2 }
$$

\monpara{Normal subgroups and quotients of $F_n$.}
Two epimorphisms $h_1: F_n\ra G_1$, $h_2: F_n\ra G_2$ represent the
same point in $\calg _n$ if and only if they have the same kernel.
Thus $\calg _n$ can be viewed as the set of normal subgroups of $F_n$.
Equivalently $\calg _n$ can be viewed as the set of quotient groups of
$F_n$. A quotient $F_n/N$ corresponds to the group $F_n /N$ marked by
the image of $(s_1,\dots ,s_n)$.
As a convention, we shall sometimes use a presentation 
$\langle  s_1,\dots ,s_n ~|~ r_1,r_2,\dots\rangle$ to
represent the marked group $(\langle  s_1,\dots ,s_n ~|~
r_1,r_2,\dots\rangle, (s_1,\dots ,s_n))$ in $\calg _n$.

\begin{rem*}
More generally a quotient of a marked group $(G,S)$ is naturally
marked by the image of $S$.
\end{rem*}

Note that in the space $\calg _n$, many marked groups have
isomorphic underlying groups.
An easy example is given by
$\langle  e_1,e_2 ~|~ e_1 = 1\rangle$ and
$\langle  e_1,e_2 ~|~ e_2 = 1\rangle$.
These groups are isomorphic to $\bbZ$, but their presentations give non
isomorphic marked groups
(once marked by the generating family obtained from the presentation).
They are indeed different points in $\calg _2$.
However, there is only one marking of free groups of rank $n$ and of rank $0$
in $\calg _n$, since any generating family of cardinality $n$ of $F_n$
is a basis, and two bases are mapped onto one another by an automorphism of $F_n$.
But there are infinitely many markings of free groups
of rank $k$, for any $0< k < n$: for instance $(\langle a,b\rangle,(a,b,w(a,b)))$
give non isomorphic markings of the free group $\langle a,b\rangle$ for different
words $w(a,b)$.

\subsection{Topology on $\calg _n$}
\label{definition_topology}

\monpara{The topology in terms of normal subgroups.}
The generating set $S$ of a marked group $(G,S)$
induces a word metric on $G$.
We denote by $B_{(G,S)} (R)$ its ball of radius
$R$ centered at the identity element of $G$.

Let $2^{F_n}$ be the set of all subsets of the free group $F_n$.
For any subsets $A,A' \in 2^{F_n}$, consider the maximal radius of the
balls on which $A$ and $A'$ coincide:
$$ v(A,A' ) = \max \big\{R\in \bbN \cup \{+\infty  \}  ~|~
A\cap B_{(F_n,(s_1,\dots ,s_n))}(R) =
A' \cap B_{(F_n,(s_1,\dots ,s_n))} (R)  \big\}.$$
It induces a metric $d$ on $2^{F_n}$ defined by $ d(A,A' ) = e^{-v(A,A' )}$.
This metric is ultrametric and makes $2^{F_n}$
a totally discontinuous metric space, which is compact
by Tychonoff's theorem.

The set $\calg _n$ viewed as the set of normal subgroups of $F_n$
inherits of this metric.
The space of normal subgroups of $F_n$ is easily seen to be a closed
subset in $2^{F_n}$. Thus $\calg _n$ is compact.

\monpara{The topology in terms of epimorphisms}\label{topo_epimorphisms}
Two epimorphisms $F_n \ra G_1$ and $F_n \ra G_2$ of $\calg _n$ are close to
each other if their kernels are close in the previous topology.

\monpara{The topology in terms of relations and Cayley graphs.}
A \emph{relation} in a marked group $(G,S)$ is an $S$-word
representing the identity in $G$. Thus two marked groups $(G,S)$,
$(G',S')$ are at
distance at most $e^{-R}$ if they have exactly the same relations of
length at most $R$.
This has to be understood under the following
abuse of language:

\begin{conv*}
Marked groups are always considered up to isomorphism of marked
groups.
Thus for any marked groups $(G,S)$ and $(G',S')$ in $\calg _n$,
we identify an $S$-word with the corresponding
$S'$-word under the canonical bijection induced by $s_i \mapsto s'_i$,
$i=1\dots n$.
\end{conv*}

In a marked group $(G,S)$, the set of relations of length at most $2L+1$
contains the same information as the ball of radius $L$ of its Cayley graph.
Thus the metric on $\calg _n$ can be expressed in term of the Cayley
graphs of the groups:
two marked groups $(G_1,S_1)$ and $(G_2,S_2)$ in $\calg _n$ are at
distance less than $e^{-2L+1}$ if their labeled Cayley graphs have the
same labeled balls of radius $L$.

\subsection{Changing the marker}
\label{marker}

Let $E$ be a finitely generated group. 
We introduce three equivalent definitions of the set of groups marked by $E$.

\monpara{Normal subgroups.}
The set $\calg(E)$ of \emph{groups marked by $E$} is the set of normal subgroups of $E$.

\monpara{Epimorphisms.}
$\calg(E)$ is the set of equivalence classes of epimorphisms from $E$ to variable groups.

\monpara{Markings.}
Given a marking $S_0$ of $E$ with $n$ generators, $\calg (E)$
corresponds to the closed subset of $\calg _n$ consisting of marked groups
$(G,S)\in \calg _n$ such that any relation of $(E,S_0)$ holds in
$(G,S)$.\\

In the last definition, the marking of $E$ corresponds to a morphism
$h_0: F_n \ra E$, and the embedding $\calg (E) \hookrightarrow \calg
_n$ corresponds, in terms of epimorphisms, to the map $f\mapsto f\circ h_0$. 

The set $\calg (E)$ (viewed as the set of normal subgroups of $E$)
 is naturally endowed with the topology induced by
Tychonoff's topology on $2^E$. This topology is the same as the
topology induced by the embeddings into $\calg _n$ described
above. Therefore $\calg (E)$ is compact.

The following lemma is left as an exercise:

\begin{lem}\label{lem_open}
  Let $h_0:E'\onto E$ be an epimorphism and let $h_0^*:\calg(E)\ra\calg(E')$
be the induced map defined in terms of epimorphisms by
$h_0^*:h\mapsto h\circ h_0$.

Then $h_0^*$ is an homeomorphism onto its image. Moreover, $h_0^*$ is open if and only if $\ker h_0$
is the normal closure of a finite set.
In particular, when $E'$ is a free group, $h_0^*$ is open if and only if $E$ is finitely presented.
\end{lem}

Because of this property, we will sometimes restrict to the case where the
marker $E$ is finitely presented.
A typical use of lemma \ref{lem_open} is to embed $\calg(E)$ into
$\calg _n$, or $\calg _n$ into $\calg _m$ for $n\leq m$, as an
open-closed subset.
For example consider the epimorphism
$h_0:F_{n+1} = \langle e_1,\dots,e_{n+1}\rangle \onto F_n = \langle f_1,\dots,f_n\rangle$ that
sends $e_i$ to $f_i$ for $i=1,\dots ,n$, and $e_{n+1}$ to $1$.
Then $h_0^*$ embeds $\calg _n$ into $\calg _{n+1}$ in the following way:
a marked group $(G, (g_1,\dots ,g_n))$ of $\calg _n$ will correspond to the marked
group $(G, (g_1,\dots ,g_n,1))$ of $\calg _{n+1}$.

\subsection{Examples of convergent sequences}
\label{examples_sequences}

\monpara{Direct limits.}
An infinitely presented group
$\langle s_1,\dots ,s_n ~|~ r_1,r_2,\dots ,r_i,\dots \rangle$,
marked by the generating family $(s_1,\dots ,s_n)$,
is the limit of the finitely presented groups
$$\langle s_1,\dots ,s_n ~|~ r_1,r_2,\dots ,r_i\rangle$$
when $i\ra +\infty$, since for any radius $R$, balls of radius $R$ in
the Cayley graphs
eventually stabilize when adding relators.

\begin{figure}[htbp]
\begin{center}
\includegraphics[width=10cm]{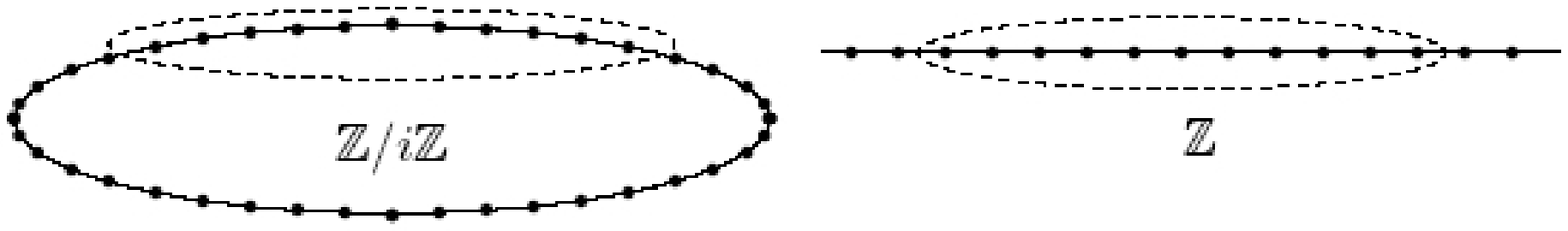}
\caption{$\bbZ/i\bbZ$ converging to $\bbZ$}
\end{center}
\end{figure}

\monpara{$\bbZ$ as a limit of finite cyclic groups.}
For any $i$, the ball of radius $i/3$ in the marked group
$(\bbZ / i\bbZ , (\overline{1}))$ is the same as the ball of radius 
$i/3$ in $(\bbZ , ( 1 ))$.

In other words, the sequence of marked groups
$(\bbZ / i\bbZ , (\overline{1} ))$ converges
to the marked group $(\bbZ , (1))$ when $i\ra +\infty$

\begin{figure}[htbp]
\begin{center}
\includegraphics[width=10cm]{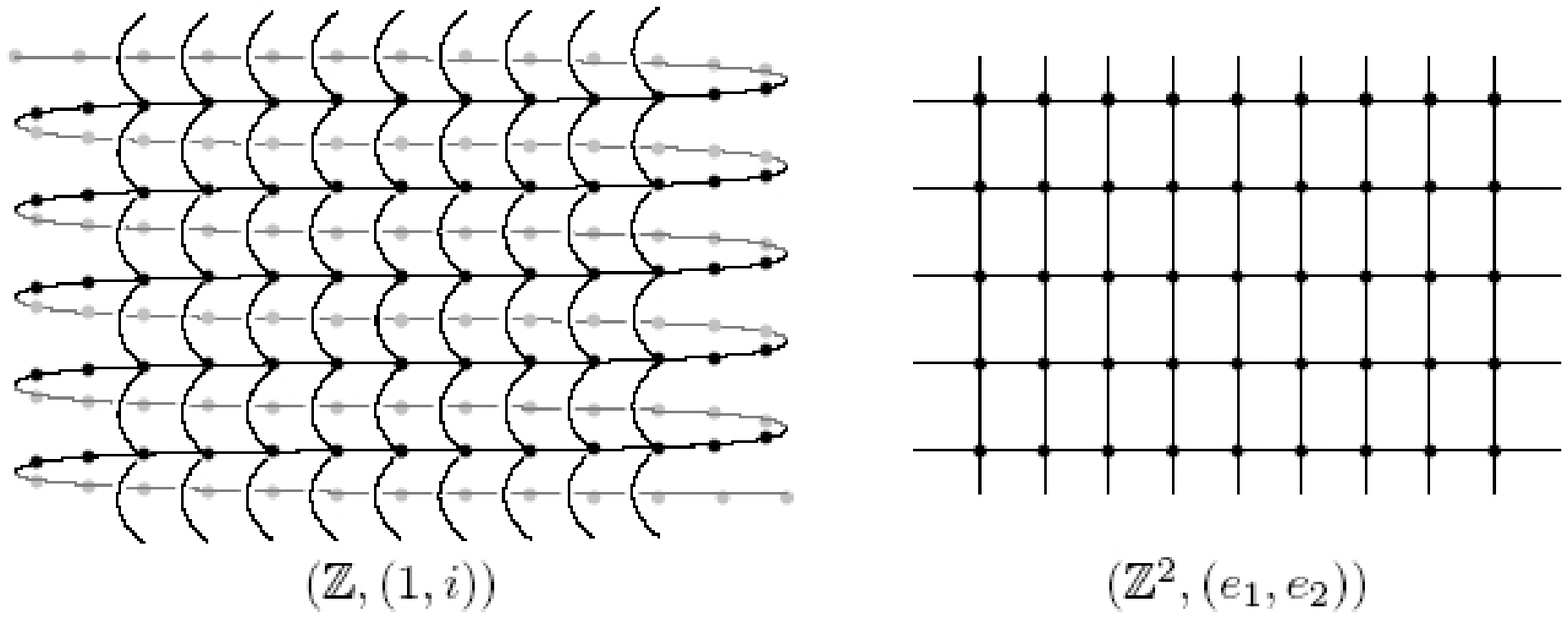}
\caption{$\bbZ^2$ as a limit of markings of $\bbZ$}
\end{center}
\end{figure}

\monpara{$\bbZ^2$ as a limit of markings of $\bbZ$.}
\label{Z_converge_to_Z2}
A less trivial example shows that a fixed group
with a sequence of different markings may converge to a non
isomorphic group.
For any $i\in\bbN$, consider the marked group $(G_i,S_i) = (\bbZ , (1, i ))\in\calg_2$.
For any $R$, if $i\geq 100R$, the only relations between $1$ and $i$
in the ball of radius $R$ of $(G_i,S_i)$ are relations of commutation.
Thus the ball of radius $R$ in $(G_i,S_i)$ is the same as the 
ball of radius $R$ in $(\bbZ ^2 , ((1,0),(0,1)))$.

In other words, the marked groups $(\bbZ , (1, i ))$ converge when $i\ra +\infty$
to the marked group $(\bbZ ^2 , ((1,0),(0,1)))$.
Of course, the same argument shows that $\bbZ^n$ marked by its
canonical basis is a limit of markings of $\bbZ$ in $\calg_n$.

Now consider the marking  $(\bbZ^k,S)$ in $\calg_n$ ($k\leq n$) defined by
taking $(s_1,\dots,s_k)$ a basis of $\bbZ^k$ and $s_i=1$ for $i>k$.
One can deduce from the argument above that there is a sequence of markings $(\bbZ,S_i)$ converging
to  $(\bbZ^k,S)$ in $\calg_n$ where the last $n-k$ generators of $S_i$ are trivial in $\bbZ$.
More conceptually, we use the continuity of the embedding $\calg_h\into\calg_n$
defined by $(G,(s_1,\dots,s_k))\mapsto(G,(s_1,\dots,s_k,1,\dots,1))$
to prove in corollary \ref{cor_abelian} that \emph{any} marking of $\bbZ^k$ is a limit
of markings of $\bbZ$.

\monpara{$F_k$ as a limit of markings of $F_2$.}
\label{F_2convergetoF_k}
A similar argument could be used for non abelian free groups to prove
that a sequence of markings of $F_2$ can converge to a free group of
rank $n$ in $\calg _n$.
Consider a free group $F_2 = \langle a,b\rangle$.
For any large $L$, choose random words $w_1 (a,b),\dots ,w_{n-2}(a,b)$ of
length $L$ in $\langle a,b\rangle$, so that $w_1 (a,b),\dots
,w_{n-2}(a,b)$ satisfy small cancellation $C'(1/100)$ property.
Consider the marking $S_L = (a,b, w_1 (a,b),\dots ,w_{n-2}(a,b))$ of $F_2$. 
Classical arguments from small cancellation
theory show that there are no relations of
length less than $L/2$ between elements of $S_L$ (since there are no
relations of length less than $L/2$ in the small cancellation group
$<a,b~|~w_1 (a,b),\dots ,w_{n-2}(a,b)>$).
Therefore, as $L$ tends to infinity, the sequence $(F_2, S_L)$ 
converges to a free group of rank $n$ marked by a free basis.

The same argument as above shows that 
there is a marking of $F_k$ is $\calg_n$ for $k\leq n$ which is a limit 
of markings of $F_2$.
It will be proved in corollary \ref{cor_indep}
that \emph{any} marking of $F_k$ is a limit
of markings of $F_2$ in $\calg_n$.

\monpara{A residually finite group is a limit of finite groups.}
\label{finite_converge_to_fp_residually_finite}
Let $G$ be a residually finite group, and $(G,S)$ a marking of $G$.
For any $i$, there is a finite quotient $G_i$ of $G$
in which the ball of radius $i$
of $(G,S)$ embeds. Denote by $S_i$ the image of $S$ in $G_i$.
Thus $(G_i,S_i)$ has the same ball of radius $i$ as $(G,S)$. Therefore,
$(G,S)$ is the limit of the marked finite groups $(G_i,S_i)$.

There is a partial converse of this result.
First note the following easy lemma:

\begin{lem}[Neighbourhood of a finitely presented group]
\label{fp_voisinage}
Let $(G,S)$ be a marking of a finitely presented group.
There exists a neighbourhood of $(G,S)$ containing only marked quotients
of $(G,S)$.
\end{lem}

\begin{proof}
Let $\langle s_1,\dots s_n ~|~ r_1(s_1,\dots ,s_n),\dots ,r_k(s_1,\dots ,s_n)\rangle$
be a finite presentation of $(G,S)$.
If $(G',S')$ is close enough to $(G,S)$, these two groups have
sufficiently large isomorphic balls to show that
$r_1(s'_1,\dots ,s'_n),\dots ,r_k(s'_1,\dots ,s'_n)$
are trivial in $G'$.
\end{proof}

As a corollary, a finitely presented group which is a limit of finite groups
is residually finite.

\begin{pb*}
Describe the closure of the set of finite groups in $\calg _n$.
\end{pb*}

\subsection{Limit groups: first approach}
\label{first_approach}

We propose the following definition of limit groups.
It follows from lemma 1.3 (iv) of \cite{Sela_diophantine1} 
that this definition is equivalent to the original definition of Sela.

\begin{dfn}
A marked group in $\calg _n$ is a \emph{limit group}
if it is a limit of marked free groups.
\end{dfn}

We shall see in corollary \ref{cor_indep} that
being a limit group does not depend on the
marking, nor of the space $\calg_n$ where this marking is chosen.
As first examples, we have seen that finitely generated free abelian groups 
 are limit groups as limit of markings of $\bbZ$.

Example \ref{finite_converge_to_fp_residually_finite} concerning
residual finiteness can be generalized to other residual properties.
In particular, residual freeness will play a central role, 
giving first a criterion for being a limit group
(proposition \ref{frfg}).

\begin{dfn}[Residual freenesses]
\label{residual_freeness}
A group $G$ is \emph{residually free} (or \emph{$1$-residually free})
if for any element $x\in G\setminus\{1\}$,
there exist a morphism $h$ from $G$ to a free group such that
$h(x)\neq 1$.

A group $G$ is  \emph{fully residually free}
if for any finite set of distinct elements $x_1$,\dots ,$x_i$,
there exist a morphism  $h$ from $G$ to a free group such that
$h(x_1),\dots ,h(x_i)$ are distinct.
\end{dfn}

Since subgroups of free groups are free,
we could assume in this definition that the morphisms are onto.

A residually finite group is fully residually finite, since finite
direct products of finite groups are finite. For freeness, the two notions are different.
The group $F_2 \times \bbZ$ is residually free, since any non trivial element $(g_1,g_2)$
of $F_2 \times \bbZ$ has at least one coordinate $g_1$ or $g_2$ which is non trivial.
But  $F_2 \times \bbZ$ is not fully residually free
because it is not commutative transitive
(see definition \ref{dfn_commutative} and corollary
\ref{cor_F2xF2}).

As in example \ref{finite_converge_to_fp_residually_finite} above,
the following property is immediate:

\begin{prop}
\label{frfg}
A marked fully residually free group $(G,S)$ is a limit group.
\end{prop}

It results from lemma \ref{fp_voisinage} that finitely presented limit groups are
fully residually free.
A theorem of Kharlampovich-Myasnikov and Sela (\cite{KhMy_irreducible1},
\cite{Sela_diophantine1}) shows that
limit groups are finitely presented. Thus the converse of the proposition \ref{frfg} is true:
the limit groups are precisely the finitely generated fully residually free groups.
This result will be proved without the use of the finite presentation of limit groups in section \ref{sec_residual}.

\subsection{Open and closed algebraic properties}
\label{open_closed}

In this paragraph, we investigate whether a given property of groups
defines a closed or open subset of $\calg _n$.

\setcounter{monparact}{0}
\monpara{Finite groups are isolated in $\calg_n$.}
Indeed, if $(G,S)$ is a finite group of cardinal $R$,
the ball of radius $R$ of $(G,S)$ determines the group law of $G$.
Therefore any marked group with a
ball of radius $R$ isomorphic to this ball is isomorphic to $G$.

Thus finiteness is an open property. But it is not closed:
the group $\bbZ$ is a limit of markings of the finite groups $\bbZ / i\bbZ$.

\monpara{Being abelian is open and closed.}
A group generated by $S$ is abelian if and only if the elements of $S$ commute.
In other words,
a marked group $(G,S)$ is abelian if and only if a certain finite collection
of words of length $4$ (the commutators in the generators)
define the identity element in $G$.
As soon as two marked groups are close enough
(at distance less than $e^{-4}$)
in the space $\calg _n$, they are either both abelian or both non-abelian.
Thus being abelian is an open and closed property.

\monpara{Nilpotence.}
By the same argument, the property of being nilpotent of class less than a given $k$
is also both open and closed in $\calg_n$,
since this property is satisfied if and only if a finite number of
words in the generators are trivial.
Thus, being nilpotent (of any class) is an open property.
On the other hand, being solvable of length at most $k$ is a closed property but not open.

\monpara{Torsion.}
The property of having torsion is open.
Indeed, suppose $g^i=1$ in $G$, for some $g\neq 1$.
Then in any marked group $(G',S')$ close enough to $(G,S)$,
the element $g'$ in $G'$ corresponding to $g$
(in isometric balls of their Cayley graphs)
is non trivial and verifies the relation $g'^i=1$
(as soon as the balls are large enough to ``contain'' this relation).
Having torsion is not a closed property:
the group $\bbZ$ is a limit of finite groups $\bbZ / i\bbZ$.

\monpara{Rank.}\label{kurosh}
``\emph{Being generated by at most $k$ elements}'' is an open property in
$\calg _n$.
Indeed, consider generators $a_1,\dots,a_k$ of a marked group $(G,S)$
and write each $s_i\in S$ as a word $w_i(a_1,\dots,a_k)$.
In a ball of radius large enough, one reads the relation
$s_i=w_i(a_1,\dots,a_k)$. In any marked group $(G',S')$
close enough to $(G,S)$, consider $a'_1,\dots,a'_k$ the elements
corresponding to $a_1,\dots,a_k$ under the bijection between their
balls in the Cayley graphs. One can read the relation
$s'_i=w_i(a'_1,\dots,a'_k)$ in $(G',S')$, so that
$a'_1,\dots,a'_k$ generate $G'$.

In other words, the property ``\emph{being generated by less than $k$
elements}'' can be read in a finite ball of the Cayley graph of $(G,S)$.

\monpara{Commutative transitivity and CSA.}

Commutative transitivity has been introduced by B. Baumslag in
\cite{Baumslag_residually} as a criterium for a residually free group
to be fully residually free (see proposition \ref{criterium}).
CSA-groups (or \emph{Conjugately Separated Abelian} groups)
has been defined by A. Myasnikov and
V. Remeslennikov in their study of exponential groups
(\cite{MyRe_exponential2,GiKhMy_CSA}).
These two properties are satisfied by free groups, and they will be
shown to be closed, thus satisfied by limit groups.

\begin{dfn}[Commutative transitivity.]\label{dfn_commutative}
A group $G$ is said to be
{\it commutative transitive} if commutativity is a transitive relation
on $G\setminus\{ 1 \}$. In other words:
$$\forall a,b,c \in G\setminus \{ 1 \}, ~~~ [a,b] = [b,c] =1
~\Rightarrow ~ [a,c] =1.$$
\end{dfn}

A commutative transitive group has the following properties
(each one being trivially equivalent to the definition):
\begin{itemize*}
\item the centralizer of any non trivial element is abelian. 
\item if two abelian subgroups intersect non trivially, their union generates an abelian subgroup.
In other words, different maximal abelian subgroups intersect trivially.
\end{itemize*}

In particular, for any maximal abelian subgroup $H$ and
any element $g$ of a commutative transitive group $G$,
the subgroup $H$ and its conjugate $gHg\m $ are equal or intersect trivially.
A stronger property is given by the following definition:

\begin{dfn}[CSA]\label{dfn_CSA}
A group $G$ is said to be {\it CSA} 
if any maximal abelian
subgroup $H<G$ is malnormal, \ie for all $g\in G\setminus H$,
 $H\cap gHg\m = \{1\}$.
\end{dfn}

It is elementary to check that property CSA implies commutative
transitivity.
The property CSA can be expressed by universal sentences
(see section \ref{sec_universal},
and proposition 10 of \cite{MyRe_exponential2}):

\begin{prop}
A group $G$ is CSA if and only if it satisfies both following properties:
\begin{itemize*}
\item $\forall a,b,c \in G\setminus \{ 1 \}, \quad [a,b] = [b,c] =1
\ \Rightarrow \ [a,c] =1$
(commutative transitivity),
\item $\forall g,h\in G\setminus \{ 1 \}, ~~~ [h,ghg\m ] =1 ~\Rightarrow ~ [g,h] =1$.
\end{itemize*}
\end{prop}

The proof is straightforward.

\begin{cor}\label{cor_F2xF2}
Commutative transitivity and CSA are closed properties.
\end{cor}

\begin{rem*}
  This result extends with the same proof to the fact that any universal formula
  defines a closed property in $\calg _n$ (see proposition \ref{universal=closed}).
\end{rem*}

\begin{proof}
Suppose that elements $a,b,c$ in $G\setminus \{ 1 \}$ are such that
$[a,c] \neq 1$, and $[a,b] = [b,c] =1$.
Write $a,b,c$ as $S$-words and let $L$ be the maximum of their lengths.
Consider a ball of radius $4L$ in the Cayley graph of
$(G,S)$. Since the relations $[a,b] = [b,c] =1$ and the non-relations
$a\neq 1$, $b\neq 1$, $c \neq 1$ and $[a,c] \neq 1$ can be read in this ball, a marked
group close enough to $(G,S)$ is not commutative transitive.
The proof is similar for the property CSA.
\end{proof}

\begin{cor}\label{limit_csa}
Limit groups are commutative transitive and CSA.
\end{cor}

In particular the group $F_2 \times \bbZ $, marked by any generating
set, is not a limit group,
and hence is not fully residually free.

We conclude by quoting the following theorem of B. Baumslag
(theorems 1 and 3 of \cite{Baumslag_residually}).

\begin{thm}[B. Baumslag, \cite{Baumslag_residually}]
\label{criterium}
Let $G$ be a finitely generated group.
The following properties are equivalent:

\begin{enumerate*}
\item
$G$ is fully residually free.
\item
$G$ is residually free and commutative transitive.
\item
$G$ is residually free and does not contain a subgroup isomorphic to
$F_2 \times \bbZ $.
\end{enumerate*}
\end{thm}

\monpara{Orderable groups.}
A group $G$ is said to be \emph{left-orderable} (resp. \emph{bi-orderable})
if there is a total order on $G$ which is left-invariant (resp. left and right-invariant).

\begin{prop} \label{prop_ordered}
    The property of being a left-orderable (resp. bi-orderable) is closed in $\calg_n$.
\end{prop}

  \begin{proof}
    Take a sequence $(G_i,S_i)$ of ordered marked groups converging to $(G,S)$.
Let $R>0$, consider the restriction of the total order on the ball of radius $R$ of $(G_i,S_i)$,
and let $\leq_{i,R}$ be the corresponding total order on the ball of radius $R$ of $(G,S)$ (for $i$ large enough).
Since there are only finitely many total orders on this ball, a diagonal argument shows that one can
take a subsequence so that on each ball of radius $R$, the orders $\leq_{i,R}$ are eventually constant. 
Thus this defines a total order on $(G,S)$ by $g\leq h$ if for $R\geq \max(|g|,|h|)$, 
one has $g\leq_{i,R} h$ for all but finitely many indexes $i$. This order is clearly left-invariant
(resp. bi-invariant) if the orders on $G_i$ are.
  \end{proof}

Free groups are bi-orderable (this non-trivial fact uses the Magnus
embedding in a ring of formal series, see for instance
\cite{BottoMuraRhemtulla_orderable}).
Thus we get the following corollary (we thank T. Delzant who pointed
out this fact) :

\begin{cor}
Limit groups are bi-orderable.
\end{cor}

This result is well-known in the context of model theory, using that
limit groups are subgroups of a non-standard free group (see for
instance in \cite{Chi_introduction}).

It will be shown in section \ref{sec_construction} that non-exceptional surface groups
are limit groups. Therefore, the corollary implies that these groups are bi-orderable 
 which is not immediate (see \cite{Baumslag_generalised,Rolfsen_Wiest_orderings}). 
%We thank T.~Delzant who pointed out the fact that limit groups are bi-orderable.

\monpara{Properties stable under quotient.}
Consider a property $(P)$ stable under taking quotients.
Let $(G,S)$ be a finitely presented marked group satisfying $(P)$.
Then according to lemma \ref{fp_voisinage}, $(P)$ is satisfied in a
neighbourhood of $(G,S)$.

For example if a sequence of marked groups converges to a finitely
presented solvable (resp. amenable, Kazhdan-(T)) group, then all
groups in this sequence are eventually solvable (resp. amenable,
Kazhdan-(T)).

Shalom has proved in \cite{Shalom_rigidity} that any finitely
generated Kazhdan-(T) group is a quotient of a finitely presented
Kazhdan-(T) group.
This implies the following result:

\begin{prop}
Kazhdan's property (T) is open in $\calg_n$.
\end{prop}

\begin{proof}
Let $(G,S)$ be a marked group having property $(T)$.
Let $H$ be a finitely presented Kazhdan-(T) group such that $G$ is a
quotient of $H$.
Write $(G,S)$ as a direct limit of $n$-generated finitely presented
groups
$(G_1,S_1)\onto (G_2,S_2) \onto \dots$.
Then for $i$ large enough, $G_i$ is a quotient of $H$ and hence has property (T).
Now the set of marked quotients of $(G_i,S_i)$ is an open set (see
lemma \ref{lem_open}) containing $(G,S)$,
all elements of which have property (T).
\end{proof}

\subsection{The isomorphism equivalence relation}\label{sec_isomorphism}

The space $\calg_n$ is naturally endowed with the \emph{isomorphism}
equivalence relation: two marked groups are equivalent if their
underlying group (forgetting about the marking) are isomorphic.
We will denote by $[G]_{\calg_n}$ (resp. $[G]_{\calg(E)}$)
the equivalence class of $G$ in $\calg_n$ (resp. $\calg(E)$).
From the definition of $\calg(E)$ by isomorphism classes of
epimorphisms,
$[G]_{\calg(E)}$ is naturally in bijection with $\Epi(E\onto G)/\Aut(G)$ 
for the natural action of $\Aut(G)$ on $\Epi(E\onto G)$.

The dynamical properties of this equivalence relation  have been
studied in \cite{Champetier_espace}, in particular on the closure of
the set of marked hyperbolic groups.
It is shown that this equivalence relation is generated
by a pseudo-group of homeomorphisms on $\calg_n$.
This implies the following lemma, but we present here a direct proof.

\begin{dfn}[saturation]
  A subset $F\subset\calg(E)$ is \emph{saturated} if
it is a union of equivalence classes for the
isomorphism relation.
The \emph{saturation} of a subset $F\subset\calg(E)$ is the union of
equivalence classes meeting $F$.
\end{dfn}

\begin{lem}  Consider a finitely presented group $E$. 
The saturation of an open set $U\subset\calg(E)$ is open in $\calg(E)$. 

The closure of a saturated set $F\subset \calg(E)$ is saturated.
\end{lem}

\begin{proof}
The second statement follows from the first one: consider a saturated
set $F\subset \calg(E)$. 
Then the interior $U$ of its complement is the largest open set which
does not intersect $F$. 
Since $F$ is saturated, the saturation $V$ of $U$ does not meet $F$,
and since $V$ is open, $V=U$ and $\overline{F}=\calg(E)\setminus V$ is saturated.

We may restrict to the case $E=F_n$
since $\calg(E)$ is an open and closed subset of $\calg_n$.
Consider an open set $U\subset\calg_n$, $V$ its saturation, and
consider $(G,S)\in U$ and $(G,S')\in V$.
Consider a radius $R$ such that any marked group having the same ball
of radius $R$ as $(G,S)$ lies in $U$.
We need to prove that there exists a radius $R'$ such that
any marked group $(H,T')$ having the same ball of radius $R'$ as $(G,S')$
has a marking $(H,T)$ having the same ball of radius $R$ as $(G,S)$.

Express the elements of $S$ as $S'$-words $s_i=w_i(s'_1,\dots,s'_n)$.
Let $L$ and $L'$ be the maximum length of the words $w_i$ and let $R'=RL$.
Let $t_i$ be the element of $H$ corresponding to the word
$w_i(t'_1,\dots,t'_n)$.
Given a word $r(e_1,\dots,e_n)$ of length at most $2R$ in the alphabet
$\{e_1,\dots , e_n\}$,
the following properties are equivalent:

\begin{enumerate*}
\item
the word $r(t_1,\dots,t_n)$ defines a relation in $(H,T)$,
\item
the word $r(w_1(t'_1,\dots,t'_n),\dots w_n(t'_1,\dots,t'_n))$ of length
at most $2R'$ defines a relation in $(H,T')$,
\item
the word $r(w_1(s'_1,\dots,s'_n),\dots w_n(s'_1,\dots,s'_n))$ of length
at most $2R'$ defines a relation in $(G,S')$,
\item
the word $r(s_1,\dots,s_n)$ defines a relation in $(G,S)$.
\end{enumerate*}
Thus $(H,T)$ has the same ball of radius $R$ as $(G,S)$.
\end{proof}

% \begin{lem}[
%   Consider a set $\calu$ of isomorphism classes of finitely generated groups.
% Let $E$ be a quotient of $F_n$.
% Denote by $[\calu]_{\calg_n}]$ and $[\calu]_{\calg(E)}$ the set of markings of elements of $\calu$ in $\calg_n$ and $\calg(E)$.
%
% A group $G$ has a marking in the closure of $[\calu]_{\calg_n}$
% if and only if $G$ has a marking in the closure of $[\calu]_{\calg(E)}$.
%
% If $E$ is finitely presented, then the same statement holds for the interiors instead of the closures.
% \end{lem}
%
% \begin{proof}
%   Choose an epimorphism $h_0:F_n\onto E$, and let $\phi:E\ra\calg_n$ be the induced map, defined in terms of epimorphisms by
% $\phi:h\mapsto h\circ h_0$.
% Remember that $\phi$ is an homeomorphism onto a compact subset of $\calg_n$ (see lemma \ref{lem_open}),
% and $\phi$ is open if $E$ is finitely presented.
% Thus $[\calu]_{\calg(E)}=\phi\m([\calu]_{\calg_n}]$ is open (resp. closed) if and only if $[\calu]_{\calg_n}\cap \Im\phi$ is.
%
% Thus, if $[\calu]_{\calg_n}$ is closed (resp. open) then so is $[\calu]_{\calg(E)}$.
% Reciprocally, if $[\calu]_{\calg(E)}$ is open, then so is  $[\calu]_{\calg_n}\cap \Im\phi$,
% hence so is its saturated
% \end{proof}

\begin{cor}\label{cor_indep}
  Being a limit group does not depend on the marking,
nor of the set $\calg_n$ (or $\calg(E)$ for $E$ finitely presented) 
in which this marking is chosen.
\end{cor}

\begin{proof}
  The lemma shows that 
the set of limit groups is saturated, \ie
that if a marking of a group $G$ in $\calg(E)$ is a
  limit of markings of free groups,
then any other marking of $G$ is also a limit of markings of free groups in $\calg(E)$.

Consider an embedding $\calg(E)\subset\calg_n$ given by a marking of $E$.
It is clear that a limit of a marked free groups in $\calg(E)$ 
is a limit of marked free groups in $\calg_n$.
For the converse, if a group $(G,S)$ in $\calg(E)$ is a limit of marked free groups $(G_i,S_i)\in\calg_n$,
then for $i$ large enough, $(G_i,S_i)$ lies in the open set $\calg(E)$, so $(G,S)$ is a limit of marked
free groups in $\calg(E)$.
\end{proof}

\begin{rem*}
  The characterization of limit groups as finitely generated fully residually free groups
allows to drop the restriction on the finite presentation of $E$.
\end{rem*}

As a consequence of the previous results, we also get the following remark:

\begin{cor}[closure of markings of a free abelian group]\label{cor_abelian}
For $k = 1,\dots ,n$, the closure of all the markings of the free abelian group $\bbZ^k$
in $\calg _n$ is
the set of all markings of the groups $\bbZ^k$, $\bbZ^{k+1}$, \dots , $\bbZ^n$:
$$\ol{[\bbZ^k]}_{\calg_n}=[\bbZ^k]_{\calg_n}\cup [\bbZ^{k+1}]_{\calg_n}\cup\dots\cup [\bbZ^{n}]_{\calg_n}$$
\end{cor}

\begin{proof}
Consider $p\in\{k,\dots,n\}$.
As in example \ref{Z_converge_to_Z2}, it is easy to construct
sequences of markings of $\bbZ^k$ converging to some particular marking
of $\bbZ^p$.
Now since the closure of all markings of $\bbZ^k$ is saturated,
all the markings of $\bbZ^p$ are limits of markings of $\bbZ^k$.
Conversely, if a marked group $(G,S)$ is a limit of markings of $\bbZ^k$, then it is abelian and torsion-free
since these are closed properties.
Moreover, its rank is at least $k$ according to property \ref{kurosh} in section \ref{open_closed}.
\end{proof}

\subsection{Subgroups}\label{sec_subgroups}
In this section, we study how subgroups behave when
going to the limit.

\begin{prop}[Marked subgroups]
\label{prop_subgroups}
Let $(G_i,S_i)$ be a sequence of marked groups  converging to a marked group $(G,S)$.
Let $H$ be a subgroup of $G$, marked by a generating family $T=(t_1,\dots,t_p)$.

Then for $i$ large enough, there is a natural family $T_i=(t_1^{(i)},\dots,t_p^{(i)})$ of elements of $G_i$,
such that the sequence of subgroups $H_i=\langle T_i\rangle\subset G_i$ marked by $T_i$
converges to $(H,T)$ in $\calg_n$.
\end{prop}

$$\xymatrix{ 
  (G_i,S_i) \ar@{~>}[r]^{i\ra\infty}  & (G,S)\\
  \exists (H_i,T_i)\ar@{~>}[r]^{i\ra\infty} \ar@{^{(}.>}[u]     & (H,T)\ar@{^{(}->}[u] 
}
$$

Since subgroups of free groups are free, we get the following corollary:

\begin{cor}\label{cor_limit_subgroups}
  A finitely generated subgroup of a limit group is a limit group.
\end{cor}

\begin{proof}[Proof of the proposition.]
Consider $R>0$ such that the ball $B$ of radius $R$ in $(G,S)$ contains $T$.
Let $i$ be large enough so that the ball $B_i$ of radius $R$ in $(G_i,S_i)$
is isomorphic to $B$.
Let $T_i$ be the family of elements corresponding to $T$
under the canonical bijection between $B$ and $B_i$.
Then any $T$-word is a relation in $H$ if and only if for $i$ large enough,
the corresponding $T_i$-word is a relation in $H_i$.
\end{proof}

\begin{rem*}
The proposition claims that for any $R$, the
ball $B_{(H_i, T_i)} (R)$ converges to the ball 
$B_{(H, T)} (R)$ for $i$ large enough.
But one should be aware that the trace of the group $H_i$
in a ball of $(G_i,S_i)$ might not converge to the trace of $H$ in a ball of $(G,S)$.
For example, take $(G,S)=(\bbZ ^2, ((1,0),(0,1))$,
$T=((1,0))$ and $(G_i,S_i)= (\bbZ , (1,i))$. Then the trace of $H=\bbZ\times\{0\}$
in a ball of $(G,S)$ is a small subset of this ball. But on the other hand, since $H_i=G_i$,
the trace of $H_i$ in any ball of $(G_i,S_i)$ is the entire ball.
This is due to the fact that elements of $H_i$ may be short
in the word metric associated to $S_i$, but long in the word metric
associated to $T_i$.

This phenomenon, occuring here with the subgroup generated by a finite set,
does not occur if one consider the centralizer of a finite set.
This is the meaning of the next lemma.
It will be used in the proof of Proposition \ref{prop_extension_centralizer}
to build examples of limit groups.
\end{rem*}

\begin{dfn*}[Hausdorff convergence of subgroups]
  Let $(G_1,S_1)$ and $(G_2,S_2)$ be two marked groups, and
$H_1$, $H_2$ two subgroups of $G_1$, $G_2$.
We say that $H_1$, $H_2$, or more precisely that the pairs $((G_1,S_1),H_1)$ and $((G_2,S_2),H_2)$,
are $e^{-R}$-Hausdorff close if 
\begin{enumerate*}
\item the balls of radius $R$ of $(G_1,S_1)$ and $(G_2,S_2)$ coincide
\item the traces of $H_1$ and $H_2$
on these $R$-balls coincide.
\end{enumerate*}
\end{dfn*}

Denote by $Z_G(x)$ the centralizer of an element $x$ in a group $G$.

\begin{lem}[Hausdorff convergence of centralizers]
\label{centralizers}
Consider a sequence of mar\-ked groups $(G_i,S_i)$ converging to $(G,S)$
and fix any $x\in G$. For $i$ large enough, consider the element $x_i\in G_i$
corresponding to $x$ under the natural bijection between balls of Cayley graphs.

Then $Z_{G_i}(x_i)$ Hausdorff-converge to  $Z_G (x)$.
\end{lem}

\begin{proof}
The commutation of an element of length $R$ with $x$ is read in the
ball of length $2(R+|x|)$.
\end{proof}

\subsection{Free and amalgamated products}

Given two families $S=(s_1,\dots ,s_n)$ and $S'=(s'_1,\dots ,s'_{n'})$,
we denote by $S\vee S'$ the family $(s_1,\dots ,s_n,s'_1,\dots , s'_{n'})$.

\begin{lem}
\label{free_products}
Let $(G_i,S_i)\in\calg_n$ and $(G'_i,S'_i)\in\calg_{n'}$ be two sequences of marked groups
converging respectively to $(G,S)$ and $(G',S')$.
Then the sequence $(G_i*G'_i,S_i\vee S'_i)$ converges to
 $(G*G',S\vee S')$ in $\calg_{n+n'}$.
\end{lem}

This is clear using normal forms in a free product.
To generalize this statement for amalgamated products,
we first need a definition.

\begin{dfn*}[convergence of gluings]\label{dfn_gluings}
   Let $(A,S)$ and $(A',S')$ be two marked groups, and
$C$, $C'$ two subgroups of $A$, $A'$.

A {\emph gluing} between the pairs  $((A,S), C)$ and $((A',S'), C')$
 is an isomorphism
$\phi : C \ra C'$.

We say that two gluings $\phi_1$ and $\phi_2$ between the pairs
$((A_1,S_1), C_1)$ and $((A'_1,S'_1), C'_1)$,
$((A_2,S_2), C_2)$ and  $((A'_2,S'_2), C'_2)$,
are $e^{-R}$-close if :
\begin{enumerate*}
\item  $C_1$ (resp. $C'_1$) is   $e^{-R}$-Hausdorff close 
to $C_2$ (resp. to $C'_2$)
\item the restrictions $\phi_{1| C_1\cap B_{R}(A_1,S_1)}$ and
$\phi_{2| C_2\cap B_{R}(A_2,S_2)}$ coincide using the natural identification
between $ B_{R}(A_1,S_1)$ and  $B_{R}(A_2,S_2)$
\end{enumerate*}
In other words, this means that the following ``diagram'' commutes
on balls of radius $R$
\newcommand{\fl}{\wr{\scriptscriptstyle R}}
$$ \begin{matrix}
A_1  &\supset& C_1& \xra{\ \phi_1\ } & C'_1 & \subset & A'_1\\
\fl  &       & \fl&                  & \fl  &         & \fl \\
A_2&\supset& C_2& \xra{\ \phi_2\ } & C'_2 & \subset & A'_2
\end{matrix}
$$
\end{dfn*}

\begin{prop}[Convergence of amalgamated products]\label{prop_amalgamated}
  Consider a sequence of groups $G_i=A_i\underset{C_i=\phi_i(C_i)}{*} A'_i$ 
and $G=A\underset{C=\phi(C)}{*} A'$, some markings $S_i$, $S'_i$, $S$, $S'$ of $A_i$, $A'_i$, $A$, $A'$  such that
\begin{enumerate*}
  \item $(A_i,S_i)$ converges to $(A,S)$
  \item $(A'_i,S'_i)$ converges to $(A',S')$
  \item $C_i$ converges to $C$ in the Hausdorff topology
  \item $C'_i$ converges to $C'$ in the Hausdorff topology
  \item the gluing $\phi_i$ converges to $\phi$.
\end{enumerate*}

Then $(G_i,S_i\vee S'_i)$ converges to $(G,S\vee S')$.
\end{prop}

The proof is straightforward using normal forms in amalgamated products and is left to the reader.
A similar statement holds for the convergence of HNN-extensions.

\subsection{Quotients}

The following proposition shows that limits of quotients are
quotients of the limits.

\begin{prop}[Limits of quotients]
\label{prop_quotients}
Consider a sequence of marked groups $(G_i,S_i)\in\calg_n$ converging
to a marked group $(G,S)$.
For any $i$, let $H_i$ be a quotient of $G_i$, marked by the image $T_i$ of
$S_i$.
Assume that $(H_i,T_i)$ converge to $(H,T)$
(which may always be assumed up to taking a subsequence).

Then $(H,T)$ is a marked quotient of $(G,S)$.
\end{prop}

$$\xymatrix{ 
  (G_i,S_i) \ar@{~>}[r]^{i\ra\infty}\ar@{->>}[d]             & (G,S)\ar@{.>>}[d]_{\exists}\\
  (H_i,T_i)\ar@{~>}[r]^{i\ra\infty}      & (H,T) 
}
$$

\begin{proof}
Up to the canonical bijection between the families $S_i$, $T_i$, $S$, $T$,
the only thing to check is that any relation between the elements of $S$ in $G$ is
verified by the corresponding elements of $T$ in $H$.
But for $i$ large enough, the relation is verified by the elements of $S_i$ in $G_i$,
and therefore by the elements of $T_i$ in the quotients $H_i$.
Thus the relation is verified by the elements of $T$ in the limit $H$.
\end{proof}

In the case of abelian quotients, this gives the following result:

\begin{cor}\label{cor_betti}
Consider a sequence of marked groups $(G_k,S_k)$ which converge to a
marked group $(G,S)$.
Then the abelianization $(G_k^{ab},S_k^{ab})$ of $G_k$ converge to an
abelian quotient of $G$.
In particular, the first Betti number does no decrease at the limit.
\end{cor}

\begin{proof}
Clear since the rank of a converging sequence of free abelian groups
increases when taking a limit.
\end{proof}

\begin{rem*}
  The first Betti number may increase at the limit: 
one can find small cancellation presentations $\langle s_1,\dots,s_n| r_1,\dots,r_n\rangle$
with trivial abelianization and with arbitrarily large relators $r_1,\dots,r_n$.
Thus there are perfect groups arbitrarily close to free groups.
\end{rem*}

%%%%%%%%%%%%%%%%%%%%%%%%%%%%%%%%%%%%%%%%%%%%%%%%%%%%%%%%%%%%%%%%%%%%%%%%%%%%%%%%%%%%%%%%%%%%%%%%%%%%%%%%%

\section{Limit groups of Sela}
\label{limit_groups}

\subsection{Summary of simple properties of limit groups}\label{sec_summary}

\begin{dfn*}
  Denote by $\call_n$ the set of marked limit groups in $\calg_n$, \ie the
  closure of the set of markings of free groups in $\calg_n$.
\end{dfn*}

By definition, $\call _n$ is a compact subset of $\calg _n$.
Remember that the notation $[G]_{\calg_n}$ represents the set of markings
of the group $G$ in $\calg_n$ (see section \ref{sec_isomorphism}). Thus 
$$\call_n=\bigcup_{i\in\{0,\dots,n\}}\overline{[F_i]_{\calg_n}}.$$

\paragraph{Abelian vs non-abelian limit groups.}
There are actually three kinds of limit groups in $\call_n$:
the trivial group,
non-trivial free abelian groups (as limits of markings of $\bbZ$), 
and non-abelian limit groups.
Since being abelian is a closed property, the non-abelian limit groups are limits of free groups
of rank at least $2$ in $\calg _n$. 
But the closure of $[F_2]_{\calg_n}$ contains a marking of $F_k$ 
for all $l\in \{2,\dots,n\}$ (see example \ref{F_2convergetoF_k} in section \ref{examples_sequences}),
and since $\overline{[F_2]}_{\calg_n}$ is saturated,
the set of non-abelian limit groups is actually the closure of markings of $F_2$.
To sum up,
$$\call_n=\{1\}\dunion \overline{[\bbZ]}_{\calg_n} \dunion \overline{[F_2]}_{\calg_n} .$$\\

\paragraph{Properties.}
The following proposition summarizes some elementary properties
 of limit groups that we have already encountered, 
and which easily result from the topological setting.

\begin{SauveCompteurs}{elementary_properties}
\elementaryProperties
\end{SauveCompteurs}

Properties \ref{subgroups}, \ref{betti} and \ref{2_generated}
are clear in the context of fully residually free groups.
Property \ref{item_CSA} is also easy to see for fully residually free groups 
(theorem 2 of \cite{Baumslag_generalised}).

\begin{proof}
Properties \ref{item_CSA} and \ref{item_ordered} are proved in section~\ref{open_closed},
property~\ref{subgroups} in cor.~\ref{cor_limit_subgroups},
and property \ref{betti} in cor.~\ref{cor_betti}.
There remains to check property~\ref{2_generated}.

Using point \ref{subgroups}, this reduces to check that
any $2$-generated limit group is isomorphic to $F_2$, $\bbZ^2$, $\bbZ$ or $\{1\}$.
So consider a marked group $(G,\{a,b\})$ in $\calg_2$
which is a limit of free groups $(G_i,\{a_i,b_i\})$.
Assume that $a$ and $b$ satisfy a non-trivial relation.
Then so do $a_i$ and $b_i$ for $i$ large enough. 
Since $G_i$ is a free group,
$a_i$ and $b_i$ generate a (maybe trivial) cyclic group.
Since $\overline{[\bbZ,\{1\}]_{\calg_2}}=[\{1\}]_{\calg_2}\cup[\bbZ]_{\calg_2}\cup[\bbZ^2]_{\calg_2}$,
the point \ref{2_generated} is proved.
\end{proof}

%%%%%%%%%%%%%%%%%%%%%%%%%%%%%%%%%%%%%%%%%%%%%%%%%%%%%%%%%%%%%%%%%%%%%%%%%%%%%%%%%%%%%%%%%%%%%%%

\subsection{First examples of limit groups.}
\label{sec_construction}

Fully residually free groups have been studied for a long time 
(see for example \cite{Chi_introduction} and references). 
In this section, we review some classical constructions of fully
residually free groups which provide the first known examples of limit groups 
(see section \ref{first_approach}). We use the topological context to
generalize the Baumslag's extension of centralizers of free groups to
limit groups (see prop.~\ref{prop_extension_centralizer}, or \cite{MyRe_exponential2,BMR_discriminating_completions}).

\begin{dfn}\label{def_ext_cent}
Let $Z$ be the centralizer of an element in a group $G$ and $A$ be a finitely generated
free abelian group.

Then $G*_Z (Z\times A)$ is said to
be a \emph{free extension of centralizer} of $G$.
\end{dfn}

Such a group can be obtained by iterating extensions of the form $G*_Z (Z\times \bbZ)$,
which are called \emph{free rank one extension of centralizer}.
In the sequel, we might simply say \emph{extension of centralizer} instead of \emph{free}
extension of centralizer. More general extensions of centralizers are
studied in \cite{BMR_discriminating_completions}.

\paragraph{Free products.}
The set of limit groups is  closed under taking free products since a limit of free products
is the free product of the limit (Lemma \ref{free_products}):

\begin{lem}\label{lem_free_products}
The free product of two limit groups is a limit group.
\end{lem}

\paragraph{Extension of centralizers and double of free groups.}
The first non-free finitely generated fully residually free groups have been
constructed by Gilbert and Benjamin Baumslag
(\cite{Baumslag_residually} theorem 8, see also
\cite{Baumslag_generalised} theorem 1) by extension of
  centralizers.

\begin{prop}[Extension of centralizers of free groups,
  \cite{Baumslag_residually}]\label{prop_ec}
If $F$ is a free group and $C$ a maximal cyclic subgroup of $F$,
then for any free abelian group $A$,
the free extension of centralizer $F*_C (C\times A)$ is fully residually free.
\end{prop}

To prove this result, one shows that the morphisms $G\ra F$
whose restriction to $F$ is the identity map and which send each generator $e_i$ of $A$
to a power $c^{k_i}$ of a fixed non-trivial element $c$ of $C$ converge to the identity map
when the $k_i$'s tend to infinity.
This convergence is a consequence of the following lemma:

\begin{lem}[ {\cite[Proposition 1]{Baumslag_generalised}}
or {\cite[Lemma 7]{Baumslag_residually}}  ]
\label{lem_baum}
Let $a_1,\dots ,a_n$ and $c$ be elements in a free group $F_k$
such that $c$ does not commute with any $a_i$.
Then for any integers $k_0,\dots ,k_n$ large enough,
the element $c^{k_0}a_1c^{k_1}a_2\dots c^{k_{n-1}}a_nc^{k_n}$ is non-trivial in $F_k$.
\end{lem}

As a corollary of the proposition, we get the following result used to
prove that non-exceptional surface groups are limit groups. 

\begin{cor}[Double of free groups, \cite{Baumslag_residually}]
\label{double}
If $F$ is a free group and $u\in F$ is not a proper power, then the double
$F*_{u=\overline{u}}\overline{F}$ of $F$ over $u$ is fully residually free
\end{cor}

\begin{proof}
This double actually embeds in the amalgam $F*_C (C\times \langle t\rangle)$ as 
the subgroup generated by $F$ and $tFt\m$.
\end{proof}

\paragraph{Surface groups.}
The fundamental groups of the torus and of the sphere are trivially limit groups.
The fundamental group of the orientable surface of genus $2$ can be written as
the double $\langle a,b \rangle *_{[a,b]=[c,d]} \langle c,d\rangle$, and is thus a limit group.
Similarly, the fundamental group of the non-orientable surface of Euler Characteristic $-2$ is
$\langle a,b \rangle *_{a^2b^2=c^2d^2} \langle c,d\rangle$ and is a limit group.
Thus, since finitely generated subgroups of limit groups are limit groups,
all non-exceptional surface groups
(\ie distinct from the non-orientable surfaces of Euler characteristic
$1$,$0$ or $-1$) are limit groups (\cite{Baumslag_generalised}). 

The fundamental group of the projective plane (resp. of the Klein bottle)
is not a limit group since a limit group is torsion free
(resp.\ is commutative transitive);
the fundamental group $G=\langle a,b,c|a^2b^2c^2=1\rangle$ of the non-orientable surface of
Euler characteristic $-1$ is not 
a limit group since three elements in a free group satisfying
$a^2b^2c^2=1$ must commute (\cite{Lyndon_equation}, see also section \ref{sec_exampleMR} for a topological proof).

\paragraph{Extension of centralizers and double of limit groups.}
\label{sec_ecdl}
Baumslag's constructions for free groups can be generalized to limit
groups: the next propositions say that the class of limit groups is stable under
extension of centralizers and double over any centralizer.
This is proved in \cite{MyRe_exponential2} (see also \cite{BMR_discriminating_completions}).
This result is an elementary case of a more general construction
(MR-resolution) given in \cite{Sela_diophantine1}, see also section \ref{sec_gta}.

Recall that limit groups are commutative transitive (Corollary \ref{limit_csa}), 
so the centralizers of non trivial elements are precisely
the maximal abelian subgroups.

%\begin{SauveCompteurs}{cpteur_Extension}
\begin{prop}[\cite{MyRe_exponential2,BMR_discriminating_completions}.]\label{prop_extension_centralizer}
  Let $G$ be a limit group, $Z$ a maximal abelian subgroup of $G$, 
and $A$ be a finitely generated free abelian group.

Then the free extension of centralizer $G*_Z (Z\times A)$ is a limit group.
\end{prop}
%\end{SauveCompteurs}

\begin{proof}
  Consider a generating family $S=(s_1,\dots,s_n)$ of $G$, 
and a sequence of free groups $(G_i,S_i)$ converging to $(G,S)$.
View $Z$ as the centralizer of an element $x$, and for $i$ large
enough, let $Z_i$ be the centralizer of the corresponding element
$x_i$ in $G_i$.
Let $a_1,\dots,a_p$ be a basis of $A$ and consider $\Tilde S=(s_1,\dots,s_n,a_1,\dots,a_p)$ (resp.
$\Tilde S_i=(s_1^{(i)},\dots,s_n^{(i)},a_1,\dots,a_p)$)
a generating family of $\Tilde G=G*_Z (Z\times A)$ (resp. $\Tilde G_i=G_i*_{Z_i} (Z_i\times A)$).

By Baumslag's extension of centralizers, $\Tilde G_i$ is a limit group so we just need to check that
$(\Tilde G_i,\Tilde S_i)$ converge to $(\Tilde G,\Tilde S)$.

If we already know that abelian subgroups of limit groups are finitely generated (see \cite{KhMy_irreducible1,KhMy_irreducible2,Sela_diophantine1})
one can apply the convergence of amalgamated products (Proposition \ref{prop_amalgamated}):
it is clear that the inclusion maps $Z_i\ra Z_i\times A$ converge to the inclusion
$Z\ra Z\times A$ in the sense of definition \ref{dfn_gluings}.
If we don't assume the finite generation of centralizers, a direct argument based on 
the normal forms in an amalgamated product gives a proof of the result.
\end{proof}

As for free groups, the following corollary is immediate.

\begin{cor}\label{cor_double_limit}
The double of a limit group over the centralizer of one of its non
trivial element is a limit group.
\end{cor}

\begin{rem*}
The corollary also follows directly from the convergence of amalgamated products 
even when centralizers are not finitely generated (Proposition \ref{prop_amalgamated}).
\end{rem*}

%%%%%%%%%%%%%%%%%%%%%%%%%%%%%%%%%%%%%%%%%%%%%%%%%%%%%%%%%%%%%%%%%%%%%%%%%%%%%%%%%%%%%%%%%%%%%%%%%%%%%%%%%

\subsection{Morphisms to free groups and Makanin-Razborov diagrams.}
\label{sec_MR}

In this section we describe the construction of a Makanin-Razborov
diagram of a limit group. Such a diagram encodes all the morphisms from
a given finitely generated group to free groups.
This construction is due to Sela and uses two deep results in
\cite{Sela_diophantine1}.

\subsubsection{Finiteness results for limit groups.}\label{sec_finiteness}

%\coucou{liaison}
The main finiteness result for limit groups is their finite presentation.

\begin{thm}[\cite{KhMy_irreducible1,KhMy_irreducible2,Sela_diophantine1}]
Limit groups are finitely presented.  
Moreover abelian subgroups of limit groups are finitely generated.
\end{thm}

This result is proved using JSJ-decomposition by Z.~Sela
(the \emph{analysis lattice} of limit groups), and using embeddings into free
$\bbZ[X]$-groups by O.~Kharlampovich and A.~Myasnikov.
An alternative proof is given in \cite{Gui_limit} by studying free
actions on $\bbR ^n$-trees.

%%%%%%%%%%%%%%%%%%%%%%%%%%

As seen in section \ref{first_approach}, an elementary consequence of
the fact that limit groups are finitely presented
is the following corollary:

\begin{cor}
\label{limit_residually}
  A finitely generated group is a limit group if and only if it is fully residually free.
\end{cor}

A proof of this fact which does not use the finite presentation of limit groups is given in 
propositions \ref{limit_ultra} and \ref{ultra_residually}.
\\

Following Sela, for any finitely generated marked group $(E,S)$,
there is a natural partial order on the compact
set $\calg (E)$: we say that $(G_1,S_1) \leq (G_2,S_2)$ if and only if
the marked epimorphism $E \ra G_1$ factorizes through the marked
epimorphism $E \ra G_2$, \ie $G_1$ is a marked quotient of $G_2$.

\begin{lem}
\label{lem_fme}
Let $E$ be a finitely generated group.
  Any compact subset $K$ of $\calg (E)$ consisting of finitely
  presented groups has at most finitely many maximal elements,
and every element of $K$ is a marked quotient of one of them.
\end{lem}

The main ingredient in this result is Lemma \ref{fp_voisinage}
claiming that a finitely presented marked group has a neighbourhood consisting
of marked quotients of this group.

\begin{proof}
First, for any $(G,S)$ in $K$, there exists a maximal element $(G',S')\in K$
such that $(G,S)\leq(G',S')$. Indeed, apply Zorn Lemma to the set of marked groups $(G',S')\in K$
such that $(G,S)\leq(G',S')$ (this uses only the compacity of $K$). 

Now, for each maximal element $(G,S)\in K$, the set of marked quotients of $(G,S)$
is a neighbourhood of $(G,S)$ since $G$ is finitely presented. This gives a covering
of the compact $K$ from which one can extract a finite subcovering.
\end{proof}
%%%%%%%%%%%%%%%%

The next result is the first step of the construction of a
Makanin-Razborov diagram.
It claims the existence of finitely many maximal limit quotients of
any finitely generated group.

\begin{SauveCompteurs}{cpteur_Finitude}
\begin{prop}[Sela \cite{Sela_diophantine1}]\label{prop_maximal}
Let $E$ be a finitely generated group. Then there
exists a finite set of epimorphisms $E\ra G_1$,\dots ,$E\ra G_p$ from $G$
to limit groups $G_1,\dots ,G_p$ such that every morphism from $E$
to a free group factorizes through one of these epimorphisms.
\end{prop} 
\end{SauveCompteurs}

\begin{proof}  
Let $E$ be a finitely generated group.
Let $\call (E)$ be the closure in $\calg (E)$ of the set of
epimorphisms from $E$ to free groups. Thus $\call(E)$ is compact, and since every limit group is finitely presented,
Lemma \ref{lem_fme} applies.
\end{proof}

\begin{rem*}
We will give a proof of Proposition \ref{prop_maximal} which does not use the finite
presentation of limit groups in section \ref{sec_finitude}.  
\end{rem*}

In some sense, the next finiteness result means that if $G$ is a limit group, 
and $H$ is a limit group which is a strict quotient of $G$,
then $H$ is \emph{simpler} than $G$.

\begin{SauveCompteurs}{cpteur_sequence_term}
\begin{prop}\label{prop_decroissance}
  Consider a sequence of quotients of limit groups 
$$G_1\onto G_2\onto\dots\onto G_k\onto\dots$$
Then all but finitely many epimorphisms are isomorphisms.
\end{prop}
\end{SauveCompteurs}

A proof which does not use the finite presentation of
limit groups will be given in section \ref{sec_finitude}.

\begin{proof}
  Let $S_1$ be a marking of $G_1$ in $\calg_n$, and let $S_i$ the image of $S_1$ in $G_i$.
The sequence $(G_i,S_i)$ clearly converges (balls eventually stabilize) so consider $(G,S)$ 
the limit of this sequence ($G$ is the direct limit of $G_i$). 
As a limit of limit groups, $(G,S)$ is itself a limit group.
Thus $G$ is finitely presented, which implies that all but finitely many epimorphisms
are isomorphisms.
\end{proof}

\subsubsection{Shortening quotients.}
\label{sec_short}

We now introduce the second deep result of \cite{Sela_diophantine1}, namely the
fact that shortening quotients are strict quotients (\cite[Claim 5.3]{Sela_diophantine1}).
We won't give a proof of this result here.

Let $(G,S)$ be a freely indecomposable marked limit group, and let $\Gamma$ be a
 splitting of $G=\pi_1(\Gamma)$ over abelian groups.
A vertex $v$ of $\Gamma$ is said to be of \emph{surface-type} if
$G_v$ is isomorphic to the fundamental group of a compact surface $\Sigma$ with boundary,
such that the image in $G_v$ of each edge group incident on $v$ is conjugate
to the fundamental group of a boundary component of $\Sigma$ (and not to a proper subgroup),
and if there exists a two-sided simple closed curve in $\Sigma$ which is not nullhomotopic and not
boundary parallel. In other words, this last condition means that there is a non-trivial refinement of $\Gamma$
at a surface-type vertex corresponding to such a curve.
This excludes the case where $\Sigma$ is a sphere with at most
three holes or a projective plane with at most two holes.
Note that a surface with empty boundary is allowed only if no edge is incident on $v$.

A homeomorphism $h$ of $\Sigma$ whose restriction to $\partial \Sigma$ is the identity
naturally induces an outer automorphism of $G$ whose restriction to the fundamental group
of each component of $\Gamma\setminus\{v\}$ is a conjugation (see for instance \cite{Lev_automorphisms}).
We call any element of $\Aut(G)$ inducing such an outer automorphism of $G$ a 
\emph{surface-type modular automorphism} of $\Gamma$.
Similarly, if $G_v$ is abelian, any automorphism of $G_v$ which fix the
incident edge groups extends naturally to an outer automorphisms of $G$,
and we call \emph{abelian-type modular automorphism} of $\Gamma$ any element of $\Aut(\Gamma)$
inducing such an outer automorphism of $G$.
Let $\Mod(\Gamma)\subset \Aut(G)$ be the \emph{modular group} of $\Gamma$, \ie
the subgroup generated by inner automorphisms, the preimages in $\Aut(G)$ of Dehn twists along edges
of $\Gamma$, and by abelian-type and surface-type modular automorphisms.

Of course, $\Mod(\Gamma)$ depends on the splitting $\Gamma$ considered.
To define a modular group $\Mod(G)\subset \Aut(G)$, one could think of 
looking at the modular group of a JSJ decomposition of $G$. But
the JSJ decomposition of $G$ is not unique and it may be not invariant under
$\Aut(G)$. However, the \emph{tree of cylinders} of a JSJ decomposition
of $G$ (defined in appendix \ref{sec_CSA}) does not depend on the JSJ splitting
considered, and is thus invariant under automorphisms of $G$ (see \cite{GL_automorphismes}).
Therefore, we call  \emph{canonical splitting} of $G$ the splitting of $G$ corresponding to
the tree of cylinders of any JSJ decomposition of $G$. We denote this canonical splitting by $\Gamma_\can$,
and we let $\Mod(G)=\Mod(\Gamma_\can)$.
This modular group is maximal in the following sense: for any abelian splitting $\Gamma$ of $G$, one has
$\Mod(\Gamma)\subset\Mod(\Gamma_\can)$ (see \cite{GL_automorphismes}).

There is a natural action of $\Mod(G)$ on $\calg(G)$ by precomposition.
Given a morphism $h$ from $G$ onto a free group, denote by $[h]_{\Mod}$ its
orbit in $\calg(G)$ under $\Mod(G)$. 
For every such orbit $[h]_{\Mod}$, Sela introduces some preferred representants
called \emph{shortest morphisms}.  We give a slightly different definition.

For any morphism $h$ from $G$ onto a free group $F$, we define the
\emph{length} $l(h)$ of $h$:

$$ l(h) = \min _{B\text{ basis of }F}\, \max_{s\in S} \, |h(s)|_B$$

where $|.|_B$ is the word metric on $(F,B)$.

\begin{dfn}[Shortest morphisms, shortening quotients]
A morphism $h\in\calg(G)$ is called \emph{shortest} if 
$$l(h)=\inf_{h'\in [h]_{\Mod}} l(h').$$

The closure of the set of shortest morphisms in $\calg(G)$ is called the set of
\emph{shortening quotients}.
\end{dfn}

\begin{thm}[{\cite[Claim~5.3]{Sela_diophantine1}}]
\label{thm_short}
  Let $G$ be a freely indecomposable limit group.

Every shortening quotient is a strict quotient of $G$.
\end{thm}

The closure of shortest morphisms is called the set of \emph{shortening quotients}.

\begin{rem*}
  Our definition slightly differs from the definition by Sela.
In our definition, the length
is a well defined function on the subset of $\calg(G)$ consisting of marked free groups
(however this function depends on the choice of a generating set $S$ of $G$).
In other words, if $h:G\ra F$ is a morphism, and $\tau$ is
an automorphism of $F$, $h$ and $\tau\circ h$ represent the same element of $\calg(G)$.
Thus,  if $F$ has a preferred basis $B$,  $l(h)= \min _{\tau\in\Aut F}\, \max_{s\in S} \, |\tau\circ h(s)|_B$.
On the other hand,
Sela does not take the infimum on the set all automorphisms of $F$, but only on the
set of inner automorphisms. 
%This means that Sela's length
%is not constant on the set of representants of a given marking of $F$.
But the limit of a sequence of morphisms depends only on the marked group
they induce. Hence the theorem stated below follows from the one proved by Sela
using Sela's definition of length. 
\end{rem*}

\begin{cor}\label{cor_fmsq}
  Let $G$ be a freely indecomposable limit group.

Then there are finitely many
maximal shortening quotients, and any shortening quotient is
a quotient of one of them.
\end{cor}

\begin{proof}
This follows from lemma~\ref{lem_fme} since the set of shortening quotients of $G$
is a compact of finitely presented groups in $\calg(G)$.
\end{proof}

\subsubsection{Makanin-Razborov diagrams}\label{subsec_MR}

\begin{figure}[htbp]
\begin{center}
\includegraphics[width=14.5cm]{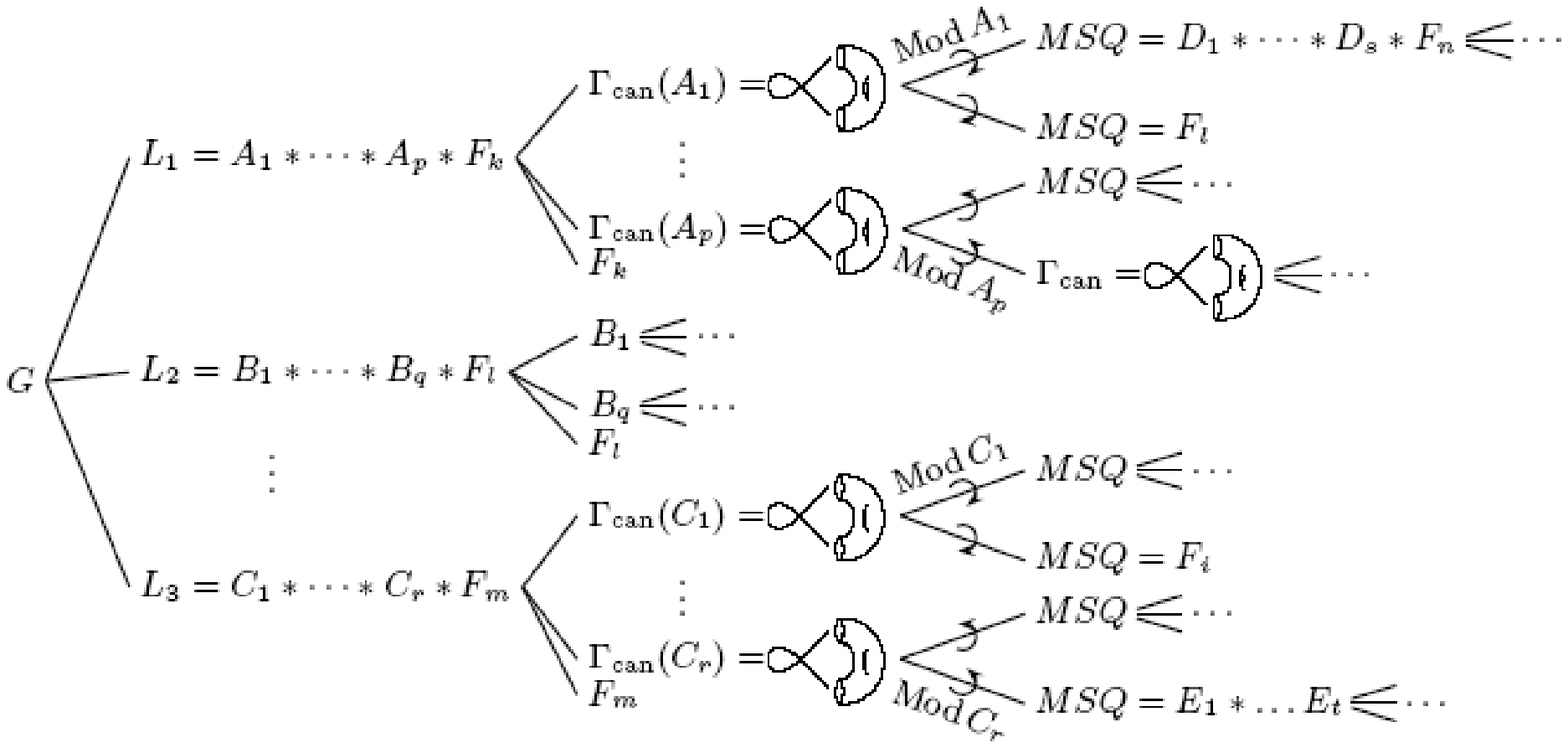}
\caption{Construction of a Makanin-Razborov Diagram}
\label{fig_MRD}
\end{center}
\end{figure}

The main application of Theorem \ref{thm_short} (see also Theorem 2
and 3 of \cite{KhMy_irreducible2}) is the construction of a Makanin-Razborov
diagram $\cald(G)$ of a finitely generated group $G$ (see figure
\ref{fig_MRD}).

This diagram is a labeled rooted tree where the root vertex is labeled by $G$, each other vertex is labeled by a limit group,
and each non-oriented edge is labeled by a morphism (which may go upwards on downwards).
Recall that in a rooted tree $\calt$, %the \emph{level} $l(v)$ of a vertex $v\in\calt$ is its distance to the root vertex.
a child of a vertex $v\in\calt$ is a vertex $u\in\calt$ adjacent to $v$, which is further from the root than $v$ is.

The children of the root vertex are labeled by the maximal limit quotients of $G$ in $\calg(G)$
(Prop.~\ref{prop_maximal}), and the edges originating from the root are labeled by
the natural morphisms from $G$ to its quotients.
We now construct inductively $\cald(G)$ by describing the children of any non-root vertex $v$ of the diagram. 
Let $G_v$ be the vertex group at $v$.

If $G_v$ is
  freely decomposable and is not a free group, write a Grushko decomposition $G_v=H_1*\dots*H_k*F_l$ 
  we define the children of $v$ to be $H_1,\dots,H_k$, $F_l$.
  We take as (upwards) edge morphisms the inclusions of $H_1,\dots,H_k,F_l$ into $G_v$.

If $G_v$ is freely indecomposable, we define the children of $v$ to be its
maximal shortening quotients.
The edges originating from $v$ are labeled by the natural morphisms from $G_v$ to 
the corresponding maximal shortening quotient.

If $G_v$ is a free group, 
  then $v$ is a leaf of $\cald(G)$. 

Since $\cald(G)$ is locally finite (Prop.~\ref{prop_maximal}) and has no infinite ray 
(Prop.~\ref{prop_decroissance}), $\cald(G)$ is finite.\\

The main feature of this diagram $\cald(G)$ is that any morphism from $G$ to a free group 
can be read in this diagram inductively in terms of morphisms of free groups to free groups
and of modular automorphisms in the following manner.
We call a \emph{Makanin-Razborov Diagram} such a diagram:

\begin{dfn*}[Makanin-Razborov Diagram]
  Given a finitely generated group $G$,
consider a finite rooted tree $\cald(G)$ whose 
root vertex is labeled by $G$, and whose other
vertices are labeled by
limit groups,
and such that each edge joining a vertex $u$ to one of its children $v$ is labeled 
either by a epimorphism  $G_u\onto G_v$ (\emph{downwards edge}) or by a monomorphism $G_v\hookrightarrow G_u$
(\emph{upwards edge}).

We say that $\cald(G)$ is a Makanin-Razborov diagram of $G$ 
if for any vertex $v$, any morphism $h:G_v\ra F$ 
to a free group $F$ can be understood in terms of morphisms from its children groups to $F$ 
in one of the four following ways: 
\begin{enumerate*}
\item $v$ is the root vertex, all edges originating at $v$ are downwards,
and $h$ factorizes through
one of the epimorphisms labeling theses edges;

\item 
$G_v$ is freely indecomposable, all edges from $v$  to its children are downwards,
and there exists a modular automorphism $\tau\in\Mod(G)$
such that $h\circ\tau$ factors through one of the epimorphisms labeling the edges between
$v$ and its children.

\item $G_v$ is freely decomposable but not free, all edges from $v$ to its children are upwards,
and $G_v$ has a non-trivial Grushko decomposition of the form $G_v=i_1(H_1)*\dots*i_k(H_k)*i_{k+1}(F_l)$
where $H_1,\dots,H_k,F_l$ are the groups labeling the children of $v$,
and $i_1:H_1\into G_v\dots$, $i_k:H_k\into G_v$, $i_{k+1}:F_l\into G_v$ are the edge monomophisms.
In that case, one has $\Hom(G_v,F)\simeq\Hom(H_1,F)\times\dots\times\Hom(H_k,F)\times\Hom(F_l,F)$
by the natural map $h\mapsto (h\circ i_1,\dots,h\circ i_{k+1})$,
therefore $h$ can be understood in terms of morphisms from its children groups $H_1,\dots,H_k,F_l$
to $F$. 

\item $G_v$ is a free group and $v$ has no child.
Note that a morphism $h:G_v\ra F$ is just a ``substitution''.
\end{enumerate*}
\end{dfn*}

\subsection{Examples of Makanin-Razborov Diagrams}\label{sec_exampleMR}

In this section, we give some examples of Makanin-Razborov diagrams.
Except in the first few cases, we won't actually
describe the result of Sela's construction (in particular, we will not make
explicit the set of shortening quotients), but we will rather describe another 
Makanin-Razborov diagram.

\paragraph{Trivial examples.}
if $G$ is a free group, then $G$ is its only maximal limit quotient,
and its Makanin-Razborov Diagram is reduced to $G\ra G$.
If $G$ has finite abelianization,
then its only limit quotient is the trivial group $L=\{1\}$,
and the Makanin-Razborov Diagram of $G$ is reduced to $G\ra \{1\}$.

\paragraph{Abelian groups.}
We now consider the case where $G$ is abelian.
In this case, $G$ has a unique maximal limit quotient $L$ obtained
by killing torsion elements. If $G$ is virtually cyclic, then $L$
is a free group and the Makanin-Razborov diagram of $G$ is the segment $G-L$.
Otherwise, any two kernels of epimorphisms from $L$ to $\bbZ$ differ by an automorphism of $L$.
Since the modular group of $L$ is its full automorphism group,
this means that the diagram $G\ra L\ra \bbZ$ is a Makanin-Razborov Diagram of $G$.

For $G=\bbZ^p$ endowed with its standard marking $(e_1,\dots,e_p)$,
we can easily work out the output of Sela's construction by computing maximal shortening quotients.
If $h:G\ra \bbZ$ is an epimorphism, its length is
 $l(h)=\max \{|h(e_i)|\, |\, i=1,\dots,p\}$ (see section \ref{sec_short}).
But there is an automorphism $\tau$ of $G$ such that $h\circ \tau (e_1)=1$
and $h\circ\tau(e_i)=0$ for $i=2,\dots,p$. It means that all shortest morphisms of $(G,(e_1,\dots,e_p\})$
have length 1. Thus, a shortest morphism consists in sending each generator $e_i$ to
$0$, $1$, or $-1$ in $\bbZ$.
In particular, there are finitely many shortest morphisms $h_1,\dots,h_n$.
Therefore, every shortening quotient is a shortest
morphism. Moreover, each of them is maximal.
Since there are several maximal shortening quotients, the output of Sela's construction of the 
Makanin Razborov diagram has several terminal vertices labeled by $\bbZ$, where the edge morphisms correspond to
 $h_1,\dots,h_n$.
The fact that all the morphisms $h_1,\dots,h_n$ are all in the same orbit under $\Aut L$ means that we can keep
only one of them to get a Makanin-Razborov diagram.

\paragraph{Surface groups.}
We now describe a Makanin-Razborov diagram of a surface group $G=\pi_1(\Sigma)$
(but without describing the maximal shortening quotients of
$G$). This problem has been studied by many authors 
(see 
\cite{Piollet_solutions},
\cite{Stallings_problems},
\cite{ComEd_solutions},
\cite{GriKu_classification}).
Let's first introduce a definition.
\begin{dfn}
A \emph{pinching} of a surface $\Sigma$ is
a family $\calc$ of finitely many disjoint simple closed curves such that
\begin{itemize*}
\item \label{2side} each curve in $\calc$ is two-sided;
\item \label{nonsep} $\Sigma\setminus\calc$ is connected;
\end{itemize*}
\end{dfn}

Corresponding to a pinching $\calc$ of $\Sigma$, 
there is natural free quotient of $G=\pi_1(\Sigma)$:
the quotient $F_\calc$ of $G$ by the normal subgroup $N_\calc$ generated by the fundamental group
of the connected component of $\Sigma\setminus\calc$ is free of rank $\#\calc$.
As a matter of fact, let $G=\pi_1(\Gamma_\calc)$ be the graph of
groups decomposition corresponding to $\calc$.
The graph $\Gamma_\calc$ has one vertex corresponding to the
connected component of $\Sigma\setminus \calc$. Edges
of $\Gamma_\calc$ correspond to the connected components of $\calc$.
Vertex and edge groups are
the fundamental groups of the corresponding subsets of $\Sigma$.
The underlying graph $\calg_\calc$ of $\Gamma_\calc$ is thus a rose having one edge for each curve
of $\calc$, and $F_\calc$ is the fundamental group of this rose.

A pinching is \emph{maximal} if it cannot be enlarged into a pinching.
Clearly, if $\calc\subset\calc'$ are pinchings, then $h_\calc$ factors through $h_{\calc'}$.
This is why we will only need to consider maximal pinchings of $\Sigma$.

\begin{prop}[\cite{Piollet_solutions},
  \cite{ComEd_solutions}, \cite{GriKu_classification},\cite{Stallings_problems}]
\label{prop_diag_surface}
  Let $G$ be the fundamental group of a closed compact surface $\Sigma$.
Let $h:G\ra F$ be a morphism to a free group.
Then there exists a maximal pinching $\calc$ of $\Sigma$, such that $h$ factors through $h_\calc$.

Moreover, there are only finitely many maximal pinchings up to homeomorphism of $\Sigma$.
\end{prop}

\begin{rem*}
  Actually, if $\Sigma$ is orientable, or if $\Sigma$ has odd Euler characteristic, then there
is exactly one maximal pinching up to homeomorphism of $\Sigma$.
For an orientable surface of genus $g$, view this surface as the boundary of a handlebody $H$,
then the corresponding morphism to $F_g$ is induced by the inclusion $\partial H\subset H$.

Also note that this proposition implies Lyndon's result that the
fundamental group of the non orientable closed surface $\Sigma$ of Euler characteristic $-1$ is not a limit group
since any maximal pinching $\calc$ in $\Sigma$ consists of only one curve so $F_\calc$ is cyclic.
\end{rem*}

This proposition means that there exist pinchings
$\calc_1,\dots,\calc_n$ such that, for every morphism $h$ from $G$ to a free group,
there exists a modular automorphism $\tau$ of $G$ such that
$h\circ\tau$ factors
through one of the morphisms $h_{\calc_i}$.
Thus, if the Euler characteristic of $\Sigma$ is at most $-2$ (so that $G$ is a limit group), the diagram 
$$\xymatrix{&&F_{\calc_1}\\
G\ar[r]&G\ar[ur]^{h_{\calc_1}}\ar[dr]_{h_{\calc_n}}&\vdots\\
&&F_{\calc_n}}$$
is a Makanin-Razborov diagram for $G$.
For surfaces of characteristic at least $-1$, the only maximal limit quotient of $G$ is the
torsion free part of its abelianization.

\begin{proof}[Proof of the proposition]
We assume that $G$ is endowed with its standard generating set
so that $G=\langle a_i,b_i|\prod [a_i,b_i]=1\rangle$ in the orientable case,
 or $G=\langle a_i|\prod a_i^2=1\rangle$ in the non-orientable case.

Consider a morphism $h$ from $G$ to a free group $F$. 
We want to represent $h$ by a topological map.
Let $\Sigma$ be the Cayley $2$-complex corresponding to the presentation of $G$ above. 
Note that $\Sigma$ is a surface, and its a cellulation has only one $0$-cell $*$, and one $2$-cell.
Subdivide the $2$-cell into triangles to obtain a one-vertex ``triangulation'' of $\Sigma$.
Identify $F$ with the fundamental group of a rose $\Gamma$, and denote by $*$ the only vertex of $\Gamma$.
For each $1$-cell $e$ of $\Sigma$, we still denote by $e$ the corresponding element of $\pi_1(\Sigma,*)$.
We define $f:\Sigma\ra\Gamma$ as follows. Send $*$ to $*$,
for each $1$-cell $e$ of $\Sigma$, we let $f$ send $e$ to the reduced
path in $\Gamma$ representing $h(e)$.
For each $2$-cell $\sigma$, we can define $f$ so that the preimage of any point $x\in \Gamma\setminus\{*\}$
is a disjoint union of finitely many disjoint arcs, such that the endpoints of each arc lies
in the interior of two distinct sides of $\sigma$ (this is a track \`a la Dunwoody). 
This can be achieved by lifting $f$ to the universal covering, and by extending $\tilde f$ on $2$-cells
according to the model shown on figure \ref{fig_lamination}.

\begin{figure}
\begin{center}
\includegraphics[width=7cm]{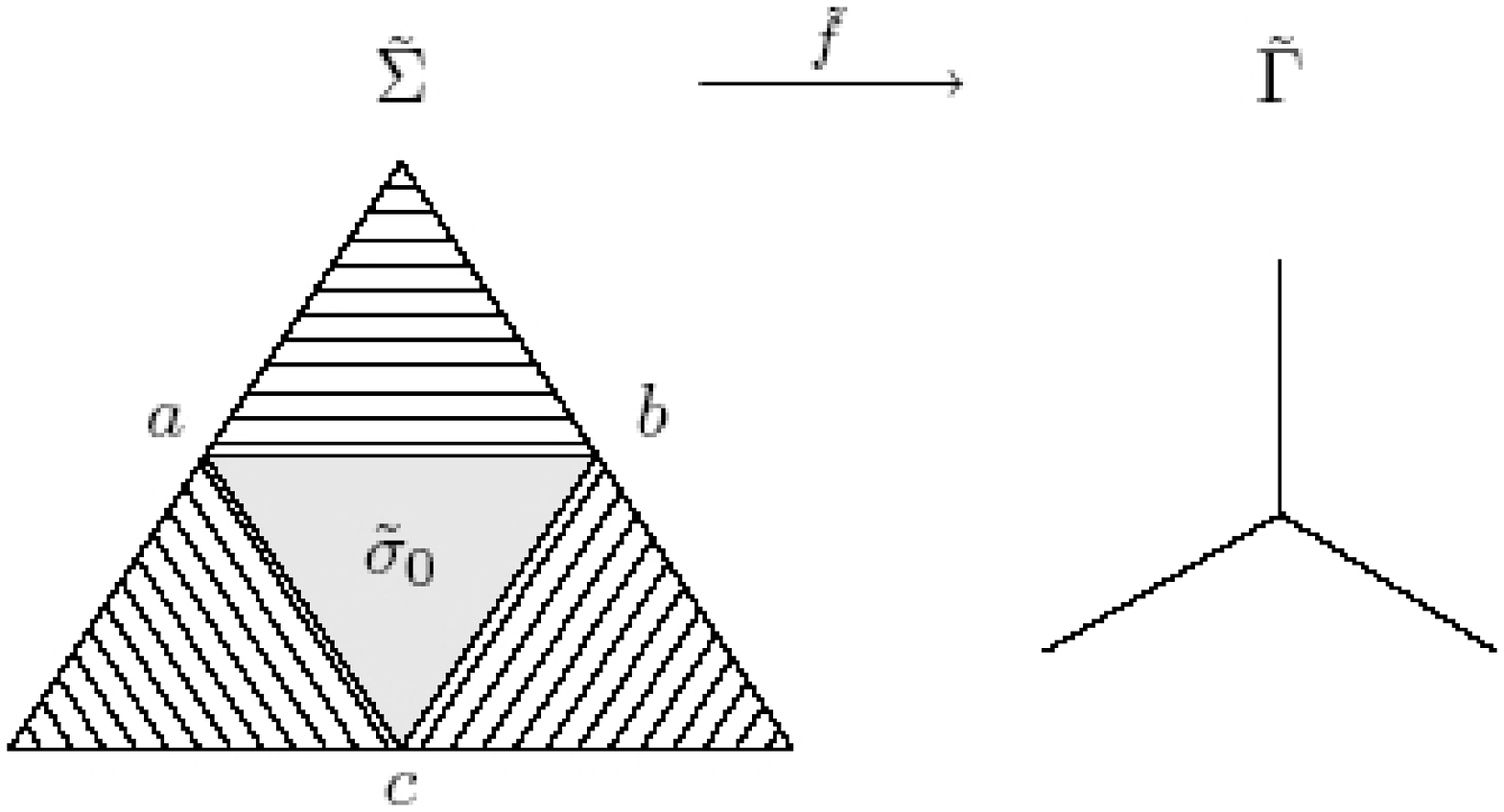}
\caption{Fibers of $\Tilde f$}\label{fig_lamination}
\end{center}
\end{figure}

Now, let $D\subset \Gamma$ be the set of midpoints of edges of $\Gamma$, and let $\calc_0=f\m(D)$.
By construction, $\calc$ is the disjoint union of finitely many
two-sided simple closed curves of $\Sigma$ (however, $\Sigma\setminus\calc_0$ may be disconnected).
Let $G=\pi_1(\Gamma_{\calc_0})$ be the corresponding graph of groups decomposition of $G$.
Clearly, each vertex group of $\Gamma_{\calc_0}$ lies in the kernel of $h$.
Denote by $\calg_{\calc_0}$ the graph underlying $\Gamma_{\calc_0}$, let $F_{\calc_0}$ be the free group 
$F_{\calc_0}=\pi_1(\calg_{\calc_0})$, and let $h_0:G\onto F_0$ be the natural map consisting in 
killing vertex groups. Thus $h$ factors through $h_0$.
Now let $T$ be a maximal subtree of $\Gamma_0$, and let $\calc_1\subset\calc_0$ be the set of curves
corresponding to edges outside $T$. Clearly, $\calc_1$ is a pinching of $\Sigma$.
Denote by $h_1:G\ra F_{\calc_1}=\pi_1(\calg_{\calc_1})$ the corresponding morphism.
Clearly, $h_1$ and $h_0$ have the same kernel so $h$ factors through $h_1$.
Now let $\calc$ be a maximal pinching containing $\calc_1$.
Then $h$ factors through the corresponding morphism $h_\calc$.

To conclude the proof of the proposition, there remains to check that the
set of maximal pinching is finite modulo homeomorphisms of $\Sigma$.

Indeed, to recover $\calc$ up to homeomorphism, it suffices to
consider the surface $\Sigma_\calc$ obtained by cutting $\Sigma$ along $\calc$,
and to know the gluing homeomorphisms between the boundary components
of $\Sigma_\calc$. 
Note that the parity of the Euler characteristic of $\Sigma_\calc$ is
the same as the one of $\Sigma$.
Since $\calc$ is a maximal pinching, every two-sided curve of $\Sigma$
disconnects $\Sigma$. 
Therefore, $\Sigma_\calc$ is either a sphere or a projective plane 
with an even number of holes (the fact that we get a sphere or a
projective plane depends only on the parity of the
Euler characteristic of $\Sigma$, not on $\calc$). 
Since any permutation of the boundary components can be realized by an
homeomorphism of $\Sigma_\calc$, up to homeomorphism there is only one
way to gather boundary components of $\Sigma_\calc$ into pairs. 
There remains to choose a gluing homeomorphism between the boundary
components in each pair.
If $\Sigma_\calc$ is a punctured sphere, there are two choices for
each pair: either the gluing homeomorphism preserves the orientation,
or not.
In particular, if $\Sigma$ is orientable, there is exactly one choice for the gluing homeomorphisms.
If $\Sigma$ is not orientable, there are exactly $c$ choices where $c=\#\calc=1-\chi(\Sigma)/2$
(choose the number of orientation-reversing homeomorphisms, it has to be between one and $\#\calc$).
If $\Sigma$ is a punctured projective plane, there is a homeomorphism
of $\Sigma_\calc$ fixing all the boundary components of $\Sigma_\calc$
except one, and which restricts to an orientation reversing
homeomorphism on the last one.
Thus in this case, the two obvious choices differ by an homeomorphism
of $\Sigma_\calc$). 
Finally, there are finitely many possible maximal families $\calc$ up
to homeomorphism of $\Sigma$ (and even exactly one in the orientable
case, or when $\Sigma$ has odd Euler characteristic). 

\end{proof}

%%%%%%%%%%%%%%%%%%%%%%%%%%%%%%%%%%%%%%%%%%%%%%%%%%%%%%%%%%%%%%%%%%%%%%%%%%%%%%%%%%%%%%%%%%%%%%%%%%%%%%%%%%%
%
%%%%%%%%%%%%%%%%%%%%%%%%%%%%%%%%%%%%%%%%%%%%%%%%%%%%%%%%%%%%%%%%%%%%%%%%%%%%%%%%%%%%%%%%%%%%%%%%%%%%%%%%%%%
%
%%%%%%%%%%%%%%%%%%%%%%%%%%%%%%%%%%%%%%%%%%%%%%%%%%%%%%%%%%%%%%%%%%%%%%%%%%%%%%%%%%%%%%%%%%%%%%%%%%%%%%%%%%%
%
%%%%%%%%%%%%%%%%%%%%%%%%%%%%%%%%%%%%%%%%%%%%%%%%%%%%%%%%%%%%%%%%%%%%%%%%%%%%%%%%%%%%%%%%%%%%%%%%%%%%%%
%
%%%%%%%%%%%%%%%%%%%%%%%%%%%%%%%%%%%%%%%%%%%%%%%%%%%%%%%%%%%%%%%%%%%%%%%%%%%%%%%%%%%%%%%%%%%%%%%%%%%%%%

\section{Constructing limit groups, fully residually free towers.}
\label{sec_gta}%got them all

Following Sela, the goal of this section is to describe
 how to construct inductively any limit
group as a graph of simpler limit groups.
We show two main ways of doing this: the first one (due to Kharlampovich and Myasnikov)
claims that any limit groups occurs as a subgroup of a group obtained from a free group
by a finite sequence of extension of centralizers (\cite[Th.4]{KhMy_irreducible2}). 
The proof we give is different from the one by Kharlampovich and Myasnikov since it
relies on Sela's techniques.
The second way of building any limit group (without passing to a subgroup) is by iterating
a construction which we call \emph{generalized double} (see definition below).
This characterization is derived from Sela's characterization of limit groups as 
\emph{strict MR-resolutions} (\cite[Th.5.12]{Sela_diophantine1}).
The arguments of this section follow the proof of Theorem 5.12 of \cite{Sela_diophantine1} up
to some technical adjustments (see the remarks following
Proposition \ref{prop_general_to_simple} and Proposition
\ref{prop_cut_surf}).

\newcommand{\subICE}{\text{sub-}\mathcal{ICE}}
\begin{dfn}
  A group is an \emph{iterated extension of centralizers of a free group}
  if it is obtained from a finitely generated free
  group by a finite sequence of free extensions of
  centralizers.

We denote by $\subICE$ the class of finitely generated subgroups of
iterated extensions of centralizers of a free group.
\end{dfn}

\begin{rem*}
  In general, one cannot obtain an iterated extension of centralizers of a free group by performing
on a free group all the extensions of centralizers simultaneously.
\end{rem*}

Clearly, the class $\subICE$ contains only limit groups.
Furthermore, the class $\subICE$ is closed under taking finitely generated subgroups, 
under free product, under free extension of centralizers, and in particular under double over
a maximal abelian subgroup.

Our first goal in this section will be the following theorem:
\begin{SauveCompteurs}{caracUn}
\caracUn
\end{SauveCompteurs}

\begin{cor}[{\cite[Cor. 6]{KhMy_irreducible2}}]
  Any limit group has a free action on a $\bbZ^n$-tree (where $\bbZ^n$ has the lexicographic ordering).

Any limit group has a free properly discontinuous action (maybe not cocompact) on a $CAT(0)$ space.
\end{cor}

\begin{proof}[Proof of the corollary]
If follows from \cite[Th. 4.16]{Bass_non-archimedean}
that if a group $G$ has a free action on a $\bbZ^n$-tree, 
then a free rank one extension of centralizers
of $G$ has a free action on a $\bbZ^{n+1}$-tree.
Similarly, if a group $G$ has a free properly discontinuous action on a 
$CAT(0)$ space, then so does a free extension of centralizers of $G$
 \cite{BridsonHaefliger_metric}.
The corollary is then clear since both properties claimed in the corollary pass to subgroups.
\end{proof}

\begin{SauveCompteurs}{dfnGeneDouble}
 \dfnGeneDouble 
\end{SauveCompteurs}

We will also say that $G$ is a generalized double over $\phi$.

\begin{rem*}
The double considered in corollary \ref{cor_double_limit}
is a particular case of generalized double:
if $G=A*_{C=\ol C} \ol A$ where $C$ a maximal abelian subgroup of $A$, one can take
$L=A$, and $\phi$ is the natural morphism sending $A$ and $\ol A$ on $A$.

Free rank one extension of centralizers is also a particular case of a generalized double:
$A*_C (C\times \bbZ)$ is isomorphic to the HNN extension $G=A*_C$
(where the two embeddings of $C$ are the inclusion), and one can take
$L=A$, and $\phi:G\ra A$ the morphism killing the stable letter of the HNN extension.
 
We will prove in next section that a generalized double over a limit group is a limit group.
More general constructions of limit groups are given in Definition \ref{dfn_sgolg} and Proposition \ref{prop_sgolg_subice},
and in Definition \ref{dfn_ggolg} and Proposition \ref{prop_general_to_simple}.
\end{rem*}

\begin{dfn}[Iterated generalized double]
  A group is an \emph{iterated generalized double} if it belongs to the smallest class of groups $\IGD$ 
containing finitely generated free groups, and stable
under free products and generalized double over a group in $\IGD$.
\end{dfn}

The second goal of this section is the following Theorem deriving from Sela's work:
\newcommand{\caracDeux}{%
\begin{thm}[Second characterization of limit groups] (Compare Sela's MR-resolution).
  A group is a limit group if and only if it is an iterated generalized double.
\end{thm}
}
\begin{SauveCompteurs}{caracDeux}
  \caracDeux
\end{SauveCompteurs}

The argument is structured as follows.
First, in section \ref{sec_subdouble}, we prove that a generalized double 
over a limit group $L$ is a subgroup of an extension of centralizers of $L$. In particular,
a generalized double is a limit group. Then, in sections \ref{sec_sgolg} and \ref{sec_twisting}, 
we extend this result to a more general situation: \emph{simple graphs of limit groups}, and we show that
those simple graphs of limit groups can be obtained using iteratively the generalized double construction
or by iteratively taking subgroups of extensions of centralizers.
Thanks to the finiteness results mentioned in section \ref{sec_finiteness},
to conclude, it will suffice to prove the following \emph{key result:}
 any non-trivial, freely indecomposable limit group can be 
written as a simple graph of limit groups over a strict quotient.
This is proved in section \ref{sec_key_result}
 using the fact that shortening quotient are strict quotients.

\subsection{Generalized double as a subgroup of a double\label{sec_subdouble}}

\begin{prop}\label{prop_subdouble}
  Let $G=A*_C B$ (resp.\ $G=A*_C$) be a generalized double over a limit group $L$.

Then $G$ is a limit group. More precisely, 
$G$ is a subgroup of a double of $L$ over a maximal abelian subgroup of $L$
(resp.\ $G$ is a subgroup of a free rank one extension of centralizers of $L$).

In particular, if $L$ is a subgroup of an iterated extension of centralizers of a free group,
then so is $G$.
\end{prop}

 We first prove the following simple lemma.
Remember that a $G$-tree is \emph{$k$-acylindrical} if the set of fix points
of any element of $G\setminus\{1\}$ has diameter at most $k$.
Accordingly, a graph of groups $\Gamma$ is $k$-acylindrical if the action of $\pi_1(\Gamma)$
on the Bass-Serre tree of $\Gamma$ is $k$-acylindrical. We will also say that $\Gamma$
is \emph{acylindrical} if it is $k$-acylindrical for some $k$.

\begin{lem}\label{lem_rewrite}
If $G=A*_C B$ is a generalized double over a limit group $L$,
then this splitting is necessarily $1$-acylindrical.

If $G=A*_C$ is a generalized double over a limit group $L$,
then either this splitting is $1$-acylindrical or it can
be rewritten so that the two embedding of $C$ into $A$ coincide. 
In the latter case, this splitting is not $k$-acylindrical for any $k$,
and $G$ is isomorphic to a free rank one extension of centralizer of $A$.
\end{lem}

\begin{proof}
Since vertex groups are CSA (they embed into the limit group $L$),
each edge group is malnormal in the neighboring vertex groups since it
is maximal abelian. Acylindricity follows in the case
of an amalgamated product.

Consider now the case of an HNN extension. Denote by $C_1$ and $C_2$ the
images of $C$ in $A$, and by $t$ the stable letter of the HNN extension.
If this splitting is not $1$-acylindrical, 
then there exists $a\in A$ such that $C_1\cap aC_2 a\m$ is non-trivial.
Since $A$ embeds into $L$, $A$ is commutative transitive so $C_1=aC_2a\m$.
Let $c_1\in C_1$ and $c_2=tc_1t\m\in C_2$ and let $c'_1=(at)c_1(at)\m\in C_1$.
Since $\phi(c'_1)$ and $\phi(c_1)$ commute and $L$ is CSA, $\phi(at)$ commutes with
$\phi(c_1)$ and $\phi(c_1)=\phi(c'_1)$.
Since $\phi$ one-to-one in restriction to $A$, one gets
that $at$ commutes with $C_1$.
Therefore, changing $t$ to $at$, the HNN extension can be rewritten as
 $G=\langle A,t\ |\ tct{}\m=c,\ c\in C\rangle$.
The lemma follows.
\end{proof}

\begin{figure}[htbp]
\begin{center}
%\fontsize{7pt}{7pt}\selectfont%
%\input{sous-double.pst}%   
\includegraphics[width=10cm]{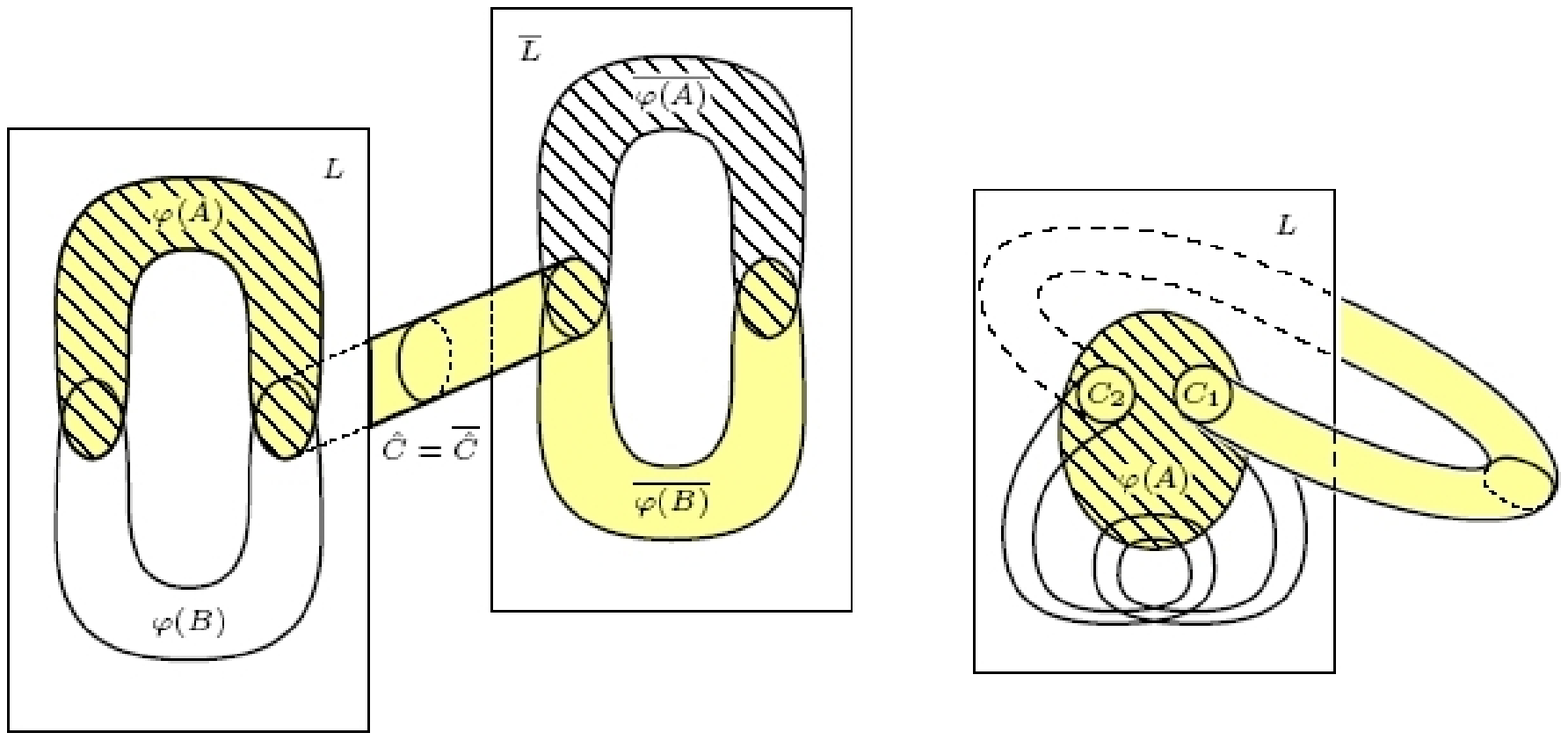}
\caption{A generalized double inside a double or an extension of centralizers}%
\end{center}
\end{figure}

\begin{proof}[Proof of the proposition]
Suppose first that $G = A*_C B$ with $\phi : G \ra L$ one to one in
restriction to $A$ and $B$, and $L$ a limit group.
We identify $A$, $B$ and $C$ with their natural images in $G$.
Let $\Hat C$ be the maximal abelian subgroup of $L$ containing $\phi (C)$.
One has $\Hat C\cap \phi(A)=\Hat C\cap \phi(B)=\phi(C)$.
Consider the double $D=L *_{\Hat C = \overline{\Hat C}} \overline{L}$.
Then the map $\psi:G\ra D$ whose restriction to $A$ and $B$ is $\phi$ 
and $\ol\phi$ respectively (with obvious notations)
is one-to-one. In other words, $G\simeq \phi(A)*_{\phi(C)=\ol \phi(C)} \ol{\phi(B)}$.
In particular, if $L$ lies in $\subICE$, so does its double $D$, so $G\in\subICE$.

Suppose now that $G=A*_C$.
If the HNN extension is not acylindrical, then $G$ is a free rank one extension of
centralizers  $G=A*_C (C\times \bbZ)$. Let $\Hat C$ be the maximal abelian subgroup of $L$
containing $\phi(C)$, and let $D$ be the free rank one extension of
centralizers $D=L*_{\Hat C} (\Hat C\times \bbZ)$.
Then the map $\psi:G\ra D$ whose restriction to $A$ is $\phi$ and sending
$\bbZ$ to $\bbZ$ is one-to-one. Thus, if $L$ lies in $\subICE$, then so does $G$.

Suppose finally that the HNN extension is acylindrical. Let $C_1$ and $C_2$ be the two
images of $C$ in $A$, and still identify $A$ with its natural image in $G$.
Let $\Hat C_1$ and $\Hat C_2$ be the maximal abelian subgroups of $L$
containing $\phi (C_1)$ and $\phi (C_2)$.
Since $C_1$ and $C_2$ are conjugate in $G$, and since $L$ is
commutative transitive, $\Hat C_1$ and
$\Hat C_2$ are also conjugate in $L$ by an element $t$.
Consider the group $D=L *_{\Hat C_1}$ where one embedding of $\Hat C_1$ is the inclusion,
and the second embedding is the conjugation by $t$. 
Since this HNN extension is not acylindrical, $D$ is a free
rank one extension of $L$.
Finally, the map $\psi:G\ra D$ whose restriction to $A$ is $\phi$ and sending
the stable letter of $G$ to the stable letter of $D$ is one-to-one.
In other words, $G\simeq \phi(A)*_{C_1} \subset D$.
\end{proof}

The following result will be used in the next section.
It controls how centralizers grow in a generalized double.

\begin{lem}\label{lem_maxsg}
  Let $G=A*_C B$ (resp.\ $G=A*_C$) be a splitting satisfying the
  hypothesis of the generalized double. Let $C'$ be a maximal
  abelian subgroup of a vertex group.
If $G$ is an amalgamated product or an acylindrical HNN extension,
then $C'$ is maximal abelian in $G$.
If $G$ is an HNN extension which is not acylindrical,
$C'$ is maximal abelian in $G$ if and only if $C'$ is not conjugate to
the edge group $C$.
\end{lem}

The proof is straightforward and left to the reader.

\subsection{Simple graphs of limit groups\label{sec_sgolg}}

In this section, we extend the notion of generalized double to some more general
graphs of limit groups to give more general constructions of limit groups.

\begin{dfn}[Simple graph of limit groups]\label{dfn_sgolg}
  A group $G$ is a \emph{simple graph of limit groups over a limit group $L$} if
$G$ is the fundamental group of a graph of groups $\Gamma$ such that:
\begin{itemize*}
\item each vertex group is finitely generated;
\item each edge group is a non-trivial abelian group whose images under both edge morphisms
are maximal abelian subgroups of the corresponding vertex groups;
\item $G$ is commutative transitive;
\item there is an epimorphism $\phi:G\ra L$ such that $\phi$ is one-to-one in restriction to each vertex group.
\end{itemize*}
\end{dfn}

We will also say that $G$ is a simple graph of limit groups over $\phi$.

\begin{rem*}
Corollary \ref{cor_acyl} in Appendix \ref{sec_CSA} shows that the requirement that $G$ is
commutative transitive would be implied by the stronger hypothesis
that $\Gamma$ is acylindrical.
\end{rem*}

\begin{SauveCompteurs}{prop_sgolg_subice} 
\begin{prop}\label{prop_sgolg_subice}
A simple graph of limit groups over a limit group $L$ is a limit group.

Moreover, if $L$ is a subgroup of an iterated extension of centralizers of a free group,
then so is $G$.

In other words, the class of groups $\subICE$ is stable under simple graph of limit groups.
\end{prop}
\end{SauveCompteurs}

\begin{proof}
  We proceed by induction on the number of edges of the graph of groups $\Gamma$.
We denote by $\phi:G\ra L$ the corresponding morphism.
If there is no edge, the proposition is trivial.
Assume first that $\Gamma$ contains an edge $e$ with distinct endpoints (\ie if $\Gamma$ has at
least two vertices).
Let $H=A*_C B$ be the amalgam carried by $e$.

Assume first that $\Gamma\setminus e$ has two connected components, and denote by 
$\Gamma_A$ and $\Gamma_B$ the components containing the vertex group $A$ and $B$ respectively.
Consider the double $D=L*_{\Hat C=\ol{\Hat C}}\ol L$
where $\Hat C$ is the maximal abelian subgroup of $L$ containing $\phi(C)$.
By Proposition \ref{prop_subdouble}, the map $\psi:H\ra D$ whose restriction to $A$ is $\phi$,
and whose restriction to $B$ is $\ol \phi$ (with obvious notations) is one-to-one.
The map $\psi$ has a natural extension to $G$ which coincide with $\phi$ (resp.\ $\ol\phi$) on
the fundamental group of $\Gamma_A$ (resp.\ $\Gamma_B$). 
One can then apply induction hypothesis to the graph of groups $\Gamma_0$ 
obtained by collapsing $e$, together with the morphism $\psi:G\ra D$:
$\psi$ is one to one in restriction to the new vertex group $H$, and each edge group of $\Gamma_0$
is maximal abelian in its neighbouring vertex groups because of Lemma \ref{lem_maxsg}.

If $\Gamma\setminus e$ is connected, 
write $G$ as the HNN extension $\pi_1(\Gamma\setminus e)*_C$ obtained from $\Gamma$ by collapsing $\Gamma\setminus e$.
Denote by $s$ the stable letter of this HNN extension.
Consider the HNN extension $D=L*_{\Hat C}$
where $\Hat C$ is the maximal abelian subgroup of $L$ containing $\phi(C)$ and
where both edge embeddings are the inclusion. Denote by $t$ the stable letter of this HNN extension.
Let $\psi:G\ra D$ whose restriction to 
$\pi_1(\Gamma\setminus e)$ is $\phi$, and sending $s$ on $t$.
One easily checks as in the proof of the previous proposition 
that $\psi$ is one-to-one in restriction to $H$.
As above, using Lemma \ref{lem_maxsg}, one can apply the induction hypothesis 
to the graph of groups $\Gamma_0$ 
obtained by collapsing $e$.

Assume now that $\Gamma$ has only one vertex and assume that there is an edge $e$ in $\Gamma$
such that the HNN extension $A*_C$ carried by $e$ is acylindrical.
Similarly, write $G$ as the HNN extension $G=\pi_1(\Gamma\setminus e)*_C$,
define the morphism $\psi:G\ra D=L*_{\Hat C}$ whose restriction to $\pi_1(\Gamma\setminus e)$
is $\phi$ and sending the stable letter of $\pi_1(\Gamma\setminus e)*_C$ to the stable letter of $D$.
The previous proposition shows that $\psi$ is one-to-one in restriction to $H=A*_C$,
and one can use induction hypothesis thanks to Lemma \ref{lem_maxsg}.

Finally, if $\Gamma$ has only one vertex and all the edges of $\Gamma$ carry a non-acylindrical
HNN extension, then $G$ is an iterated extension of centralizers of the vertex group of $\Gamma$,
which is a subgroup of $L$.
\end{proof}

\subsection{Twisting generalized doubles\label{sec_twisting}}

In this section, we give an alternative proof (due to Sela) that 
a generalized double is a limit group. Actually, we prove the more precise
result (which will be needed in the sequel)
that given a generalized double $G$ over $\phi:G\ra L$, there is a sequence of Dehn
twists $\tau_i$ such that the markings of $L$ defined by $\phi\circ\tau_i$
converge to $G$.

\begin{prop}[{\cite[Th.5.12]{Sela_diophantine1}}]\label{prop_gd}
  Consider a generalized double $G=A*_C B$ (resp.\ $G=A*_C$) over
  $\phi:G\ra L$.
  
  Then there exists a sequence of Dehn twists $(\tau_i)$ on $G$ such
  that $\phi\circ \tau_i$ converge to $\id _G$ in $\calg(G)$.
\end{prop}

\begin{proof}
  Using the fact that $L$ is fully residually free, consider a
  sequence of morphisms $\psi _i$ from $L$ to free groups $F^{(i)}$
  converging to $\id _L$, so that $\phi_i = \psi _i \circ \phi$
  converge to $\phi$ in $\calg(G)$.  We prove that for any finite set
  $g_1,\dots ,g_k$ of non-trivial elements of $G$, there exist $i$
  such that for $n$ large enough, the images of $g_1,\dots ,g_k$ under
  $\phi_i \circ \tau^n$ are all non-trivial in the free groups
  $F^{(i)}$.  Since each $\phi _i$ factorizes through $\phi$, this
  will imply that $\phi \circ \tau ^n (g_j)$ is non trivial for any
  $j$, which will prove the convergence of $\phi \circ \tau ^n$ to
  $\id _G$.  To save notation, we will treat only the case of one
  element $g\in G\setminus\{1\}$, the case $k>1$ being identical.
  
  We first consider the case of an amalgamated product.  Take $c$ an
  element of $C$ with non-trivial image in $L$ under $\phi$ and denote
  by $\tau$ the Dehn twist along $c$.  Write $g$ as a reduced form
  $g=a_1b_1... a_{p}b_{p}$, with $a_j\in A\setminus C$ and $b_j\in
  B\setminus C$ (except the maybe trivial elements $a_1$ and $b_p$).
  In particular, $C$ being maximal abelian in both $A$ and $B$, for
  all $j$, $c$ does not commute with $a_j$ nor $b_j$.  For $i$ large
  enough, $\phi _i (c)$ does not commute with $\phi_i(a_j)$, nor
  $\phi_i(b_j)$ for all $j$.  Baumslag's Lemma \ref{lem_baum} shows
  that for $n$ large enough, the image $\phi_i(a_1) \phi_i (c)^n
  \phi_i(b_1)...  \phi_i(a_{p}) \phi_i (c)^n \phi_i(b_{p})$ of $g$
  under $\phi_i \circ \tau^n$ is non-trivial.
  
  Let us now consider the slightly more subtle case of an HNN
  extension.  Denote by $C_1$ and $C_2$ the images of $C$ in $A$ under
  the two embeddings $j_1$ and $j_2$.  We write $G=\langle A,t\ |\ 
  tj_1(c)t\m=j_2(c),\ c\in C\rangle$.
  
  We first assume that $C_1\cap aC_2 a\m=\{1\}$ for all $a\in A$ (this
  means that the HNN extension is $1$-acylindrical).  We prove that
  for all $c_1\in C_1\setminus\{1\}$, $\phi(c_1)$ does not commute
  with any element of the form $\phi(at)$ for any $a\in A$.  Indeed,
  if $\phi([at,c_1])=1$, then $\phi(a\m c_1 a)=\phi(t c_1
  t\m)=\phi(c_2)$ where $c_2=j_2(j_1\m(c_1))$. Since $\phi$ is
  one-to-one in restriction to $A$, one gets $a\m c_1 a=c_2$, a
  contradiction.
  
  Let $c_1$ be a non-trivial element of $C_1$, and let $\tau$ the Dehn
  twist along $c_1$: $\tau$ restricts to the identity on $A$ and
  $\tau(t)=tc_1$.  Consider $g\in G\setminus \{1\}$, and let us prove
  that there exist $i$ such that for $n$ large enough, $\phi_i \circ
  \tau^n (g)\neq 1$ in the free group $F^{(i)}$.  The element $g$ can
  be written as a reduced form $g=a_0 t^{\eps_1}a_1 t^{\eps_2}a_2
  \dots t^{\eps_p} a_p$ with $\eps_j=\pm 1$, and where $a_j\notin C_1$
  if $\eps_j=-\eps_{j+1}=1$ (resp.\ $a_j\notin C_2$
if $\eps_j=-\eps_{j+1}=-1$).%, and 
%$a_0 \notin C_1$ if $\eps _1 =-1$ (resp.\ $a_p\notin C_2$
%if $\eps _p =-1$)

Choose $i$ large enough so that $\phi _i (c)$ does not commute with
the image of any $\phi_i(a_j t)$.  We have
$$\tau^n(g)=a_0 (tc_1^n)^{\eps_1}a_1 (tc_1^n)^{\eps_2}a_2 \dots
(tc_1^n)^{\eps_p} a_p$$
so the words $w_j$ appearing between two
powers of $c_1$ are of one of the following forms: $a_j$, $a_jt$, $t\m
a_j$ and $t\m a_j t$.  The reduced form guarantees that $\phi _i
(w_j)$ do not commute with $\phi _i (c_1)$.  Baumslag's Lemma
\ref{lem_baum} then concludes, the case where $\phi _i (a_0)$ or $\phi
_i (a_p)$ commutes with $\phi _i (c_1)$ being easy to handle.

In the case when the HNN extension is non-acylindrical, then it can be
rewritten as $G=\langle A,t\ |\ tct{}\m=c,\ c\in C\rangle$ (see Lemma
\ref{lem_rewrite}).  One checks as above that if an element of the
form $t^ka$ with $a\in A$ commutes with an element $c_1\in C$, then
$a\in C$.  Now choose a non-trivial element $c_1\in C$, and $\tau$ the
Dehn twist along $c_1$ sending $t$ to $tc_1$.  The argument above can
be adapted to this case by writing each element of $G$ as a reduced
word of the form $a_0t^{k_1}a_1t^{k_2}\dots t^{k_p}a_p$ where
$a_0,\dots a_p\in A$, $a_1,\dots, a_{p-1}\notin C$, and $k_2,\dots
k_p\neq 0$.  The argument above concludes the proof.
\end{proof}

\begin{prop}[{\cite[Th.5.12]{Sela_diophantine1}}]\label{prop_sgolg_twist}
  Consider $G=\pi_1(\Gamma)$ a simple graph of limit groups over
  $\phi:G\ra L$.
  
  Then there exists a sequence of multiple Dehn twists $\tau_i$ on
  $\Gamma$ such that $\phi\circ \tau_i$ converges to the identity in
  $\calg(G)$.
  
  Moreover, if $\phi$ is not one-to-one, then $G$ can be written as a
  generalized double over an epimorphism $\phi':G\ra L'$ which is not
  one-to-one.
\end{prop}

\begin{rem*}
The \emph{moreover} part of the proposition will be used
to prove that every limit group is an iterated generalized double.
\end{rem*}

\begin{proof}
  We argue by induction on the number of edges of $\Gamma$.  If there
  is only one edge, then we are in the situation of a generalized double and
  the proposition results from Proposition \ref{prop_gd}.
  
  Consider an edge $e$ of $\Gamma$, and let $H=A*_C B$ (resp.
  $H=A*_C$) be the subgroup of $G$ corresponding to the amalgam or HNN
  extension carried by $e$.  By Proposition
  \ref{prop_gd}, there exists Dehn twists $\tau_i$ along $e$ such that
  $\phi\circ\tau_i{}_{\ts| H}$ converges to $\id_H$ in $\calg(H)$.
  
  In the compact $\calg(G)$, extract a subsequence of
  $\phi\circ\tau_i$ converging to an epimorphism $\psi:G\ra L_0$. The
  group $L_0$ is a limit group as a limit of markings of the limit
  group $L$.  Let $\ol\Gamma$ be the graph of groups obtained from
  $\Gamma$ by collapsing $e$.  The map $\psi$ is one-to-one in
  restriction to $H$ since $\phi\circ\tau_i{}_{\ts| H}$ converges to
  $\id_H$, and $\psi$ is one-to-one in restriction to any other vertex
  group $G_v$ since $\phi\circ\tau_i$ is one-to-one on $G_v$ for all
  $i$.
  
  For the first part of the proposition, there remains to check that
  the edge groups of $\ol\Gamma$ are maximal abelian in their
  neighbouring vertex group to conclude using induction hypothesis.
  This is true if the endpoints of $e$ are distinct, or if the HNN
  extension carried by $e$ is $1$-acylindrical, since the maximal
  abelian subgroups of the vertex groups are maximal abelian in $H$
  (Lemma \ref{lem_maxsg}).
  
  Therefore, we can assume that no edge of $\Gamma$ holds an amalgam
  or a $1$-acylindrical HNN extension.  This means that $\Gamma$ is a
  multiple HNN extension of the form $G=\langle A,t_1,\dots,t_n|\, t_i
  c t_i\m= c,\ c\in G_{e_i}\rangle$ where $e_1,\dots,e_n$ are the
  edges of $\Gamma$.  In this case, we take $e=e_1$, $H=\langle
  A,t_1|\, t_1 c t_1\m= c,\ c\in G_{e_1}\rangle$, and $\ol\Gamma$ the
  graph of groups obtained by collapsing $e_1$ as above.  The fact
  that $G$ is commutative transitive implies that for $i\neq 1$,
  $G_{e_i}$ is maximal abelian in $H$.  Indeed, if $g\in H$ commutes
  with $G_{e_i}$, then $g$ commutes with $t_i$, so $g$ must act by
  translation on the axis of $t_i$ in the Bass-Serre tree of
  $\ol\Gamma$.  Since $g\in H$, $g$ is elliptic, so $g$ fixes the axis
  of $t_i$.  In particular, $g\in G_{e_i}$, so $G_{e_i}$ is maximal
  abelian in $H$.  Thus the induction hypothesis concludes the proof
  of the first part of the proposition.

  To check the \emph{moreover} part, we just need to take care of the
  case where $\psi$ is one-to-one.  In this case, consider a connected
  component $\Gamma'$ of $\Gamma\setminus e$.  We claim that $\phi$ is
  one-to-one in restriction to the fundamental group $G'$ of
  $\Gamma'$.  Indeed, $\tau_i$ restricts to a conjugation on $G'$,
  hence for all $g\in G'$, if $\phi(g)=1$, then $\phi(\tau_i(g))=1$
  for all $i$, and since $\phi\circ\tau_i$ converges to $\psi$ which
  is one-to-one, one gets $g=1$.  Therefore, by collapsing the
  connected components of $\Gamma\setminus e$, we obtain a 1-edge
  graph of groups such that $\phi$ is one-to-one in restriction to its
  vertex groups.  If the edge group is maximal abelian in both
  neighbouring vertex groups, then $G$ is a generalized double over
  $\phi$, and we are done.  Otherwise, the following claim concludes
  since a free rank one extension of centralizers is a generalized
  double over a strict quotient.
\end{proof}

\begin{claim}
  Let $G$ be a group which decomposes as a graph of groups
  $\Gamma$ with finitely generated vertex groups and non-trivial abelian edge groups
where each edge group is maximal abelian in its neighbouring vertex groups.
Assume that $G$ is a limit group
(in other words, $G=\pi_1(\Gamma)$ is a simple graph of limit groups
over the identity).  Assume that there exists a maximal abelian
subgroup $C$ of a vertex group of $\Gamma$ such that $C$ is not
maximal abelian in $G$.

Then $G$ can be written as a free rank one extension of centralizers.
\end{claim}

\begin{proof}
  We proceed by induction on the number of edges of $\Gamma$.  If
  $\Gamma$ has no edge, then the claim holds as the hypothesis is
  impossible.
  
  Assume now that $\Gamma$ contains an edge $e$ such that the $1$-edge
  subgraph of groups $\Gamma_e$ of $\Gamma$ containing $e$ is
  acylindrical.  By Lemma \ref{lem_maxsg}, every maximal abelian
  subgroup of a vertex group of $\Gamma_e$ is maximal abelian in
  $\pi_1(\Gamma_e)$.  Therefore, the graph of groups $\ol\Gamma$ obtained from
  $\Gamma$ by collapsing $e$ satisfies the hypotheses of the claim and
  induction hypothesis conclude.
  
  Otherwise, for every edge $e$ of $\Gamma$, the HNN extension
  $\Gamma_e$ is a free rank one extension of centralizers, and the
  result follows.
\end{proof}

\subsection{Statement of the key result and characterizations of limit groups}\label{sec_characterization}

The following key result will be proved in next section.

\begin{thm}[Key result (see {\cite[Th.5.12]{Sela_diophantine1}})]\label{key_result}
  Any non-trivial freely indecomposable limit group $G$ is a simple
  graph of limit groups over a strict quotient, \ie over a morphism
  $\phi:G\ra L$ which is not one-to-one.
\end{thm}

The key result allows to deduce the characterizations of limit groups:

\begin{UtiliseCompteurs}{caracUn}
  \caracUn
\end{UtiliseCompteurs}

\begin{UtiliseCompteurs}{caracDeux}
  \caracDeux
\end{UtiliseCompteurs}

\begin{proof}[Proof of the two characterization theorems]
  We have already seen that the classes $\subICE$ and $\IGD$ consist of limit
  groups.
  
  Let $G$ be a non-trivial limit group. We are going to construct
  inductively a labeled rooted tree $\calt$, where each vertex is
  labeled by a non-trivial limit group, and where the root is labeled
  by $G$.  If a vertex $v$ of $\calt$ holds a group $H$ which is
  freely decomposable, we define its children to be its freely
  indecomposable free factors. In particular, if $v$ is labeled by a
  free group, then $v$ is a leaf of $\calt$ 
  (remember that $\bbZ$ is freely decomposable so free
  groups have no freely indecomposable free factors).  If a vertex $v$
  of $\calt$ holds a non-trivial freely indecomposable limit group
  $H$, the key result provides a strict quotient $L$ of $H$ such that
  $H$ is a simple graph of limit groups over $L$. In this case, we
  attach a single child to $v$ labeled by $L$.
  
  This tree $\calt$ is locally finite, and has no infinite ray by
  the finiteness property in Proposition \ref{prop_decroissance}.  
  Thus $\calt$ is finite.
  
  Since labels of leaves of $\calt$ are free groups, they belong to
   $\subICE$, and since $\subICE$ is
  stable under free products and simple graphs of limit groups (Prop.~\ref{prop_sgolg_subice}), we
  deduce that $G$ belongs to $\subICE$.

  To prove that $G$ belongs to $\IGD$, we consider a tree $\calt'$,
  which similar to $\calt$ except in the case of freely indecomposable groups: 
  if a vertex $v$ of $\calt'$ holds a non-trivial freely indecomposable limit group
  $H$, the key result provides a strict quotient $L$ of $H$ such that
  $H$ is a simple graph of limit groups over $L$, and Proposition \ref{prop_sgolg_twist}
  gives another strict quotient $L'$ of $H$ such that $H$ is a generalized double
  over $L'$. We then attach a single child to $v$ labeled by $L'$.
  The same finiteness argument concludes that $\calt'$ is finite and that
  $G\in\IGD$ since $\IGD$ is stable under generalized double and under free product.
\end{proof}

\subsection{Proof of the key result}\label{sec_key_result}

Our aim in this section is to prove the key result (Th.~\ref{key_result}) 
\ie that any non-trivial freely indecomposable limit group $G$ is a simple graph of groups over
a strict quotient (definition \ref{dfn_sgolg}).

Since $G$ is a limit group, consider a sequence of epimorphisms $\phi_i$ from $G$ to free groups
converging to the identity in $\calg(G)$. For each index $i$, consider 
 $\tau_i$ in $\Mod(G)$ such that $\sigma_i=\phi_i\circ\tau_i$
a shortest morphism in $[\phi_i]_{\Mod}$ (see section \ref{sec_short}).
Up to taking a subsequence, $\sigma_i$ converges to a shortening quotient $\sigma:G\ra L$.
Since shortening quotients are strict quotients (Theorem \ref{thm_short}), $\sigma$ is not one-to-one. 

Next proposition will gather some properties of this morphism $\sigma$.
We first need a definition.

\begin{dfn}[Elliptic abelian neighbourhood]
Consider a graph of groups $\Gamma$ over abelian groups whose fundamental group $G$ is commutative transitive.

Consider a non-trivial elliptic subgroup $H\subset G$. The \emph{elliptic abelian neighbourhood}
of $H$ is the subgroup $\Hat H\subset G$ generated by all the elliptic elements of $G$
which commute with a non-trivial element of $H$.
\end{dfn}

\begin{rem*}
If $H$ is abelian (in particular, if $H$ is an edge group of $\Gamma$),
the elliptic abelian neighbourhood $\Hat H$ of $H$ 
is precisely the set of elliptic elements of $G$ commuting with $H$ (since this set is a group).  
  In particular, if $\Gamma$ is acylindrical, then the elliptic abelian neighbourhood
of an abelian group is its centralizer.
\end{rem*}

\begin{claim}
For a vertex $v$ of $\Gamma$, $\Hat G_v$ is the subgroup
of $G$ generated by $G_v$ and all the groups $\Hat G_e$ for $e$ incident on $v$.   
\end{claim}

\begin{proof}
  Let $g\in G_v$, and $h\in G$ and elliptic element commuting with $g$.
We just need to prove that if $g$ does not fix an edge then $h\in G_v$.
But since $h$ commutes with $g$, $h$ preserves $\Fix g$, and $\Fix g=\{v\}$
by hypothesis.
\end{proof}

\begin{prop}\label{prop_apply_mod}
  Let $G$ be a non-trivial freely indecomposable limit group.
Let $\Gamma_\can$ be the canonical splitting of $G$ (see section \ref{sec_short}).
Then, either $G$ can be written as a non-trivial free extension of centralizers, or
there exists an epimorphism $\sigma$ from $G$ to a limit group $L$ 
which is not one-to-one, and such that
\begin{itemize*}
\item $\sigma$ is one-to-one in restriction to each edge group;
  \item for each vertex $v\in\Gamma$ of surface type, $\sigma(G_v)$ is non-abelian;
  \item for each non surface type vertex $v$, $\sigma$ is one-to-one
in restriction to the elliptic abelian neighbourhood $\Hat G_v$ of $G_v$.
\end{itemize*}
\end{prop}

To prove the proposition, we will use the following simple lemma.

\begin{lem}\label{lem_conj}
  Consider a sequence of morphisms $\phi_i\in \calg(G)$ converging to $\id$ and
let $\tau_i$ be a sequence of endomorphisms of $G$ such that $\phi_i\circ\tau_i$
converge to $\psi$.

Assume that there is a subgroup $H\subset G$ such that for all index $i$, $\tau_{i|H}$
 coincides with the conjugation by an element of $G$.

Then $\psi$ is one-to-one in restriction to $H$.
\end{lem}

\begin{proof}
  Let $h\in H$, and assume that $\psi(h)=1$. Then for $i$ large enough,
$\phi_i\circ\tau_i(h)=1$, therefore $\phi_i(h)=1$. At the limit, one gets $\id_G(h)=1$ and $h=1$.
\end{proof}

\begin{proof}[Proof of the proposition]
  If $\Gamma_\can$ contains an abelian vertex group $G_v$ such that the
group generated by incident edge groups is contained in a proper free factor of $G_v$,
then $G$ can clearly be written as a non-trivial free extension of stabilizers.
Thus, from now on, we can assume that for each abelian vertex group $G_v$,
the subgroup generated by incident edge groups has finite index in $G_v$.

Therefore, each element $\tau$ of the modular group of $\Gamma$ coincides with a conjugation
in restriction to each non-surface type vertex group and to each edge group of $\Gamma$.

Consider a morphism $\sigma$ as defined above: $\sigma$ is a limit of shortest morphisms
$\phi_i \circ\tau_i$ where $\phi_i$ is a sequence of morphisms to free groups converging to
the identity.

Let $G_v$ be a non-surface type vertex group of $\Gamma_\can$.
In view of Lemma \ref{lem_conj}, to prove that $\sigma$ is one-to-one in restriction to $\Hat G_v$,
we just need to show that any modular automorphism $\tau$ restricts to a conjugation on $\Hat G_v$.

We first prove that for each edge group $G_e$, $\tau$ coincides with a conjugation
on $\Hat G_e$. Indeed, since $\Hat G_e$ is elliptic, let $G_v$ be a vertex group
containing a conjugate of $\Hat G_e$. If $G_v$ is not of surface type, then 
this is clear since the restriction of $\tau$ to $G_v$ is a conjugation. If $G_v$ is of surface type,
then $\Hat G_e$ is conjugate to an edge group of $\Gamma$, and the restriction of $\tau$ to
$\Hat G_e$ is a conjugation.

We now prove that for each non-surface type vertex group $G_v$, $\tau$ coincides with a conjugation 
in restriction to $\Hat G_v$.
Remember that $\Hat G_v$ is generated by $G_v$ and by the groups $\Hat G_e$
for $e$ incident on $v$.
Moreover, $\tau$ coincides with a conjugation $i_g$ on $G_v$,
and with some (maybe different) conjugation $i_h$ on $\Hat G_e$.
But $i_{g}\m\circ i_h$ fixes $G_e$, so $g\m h$ commutes with $G_e$.
Since $G$ is commutative transitive, $i_h$ coincides with $i_g$ on  $\Hat G_e$ so $\tau$ coincides with $i_g$ on 
$\Hat G_v$. 

Finally, if $G_v$ is a surface vertex group, then $\sigma (G_v)$ is non abelian as a limit of the
non abelian groups $\phi_i \circ \tau_i (G_v)$ in $\calg (G_v)$.
\end{proof}

In a simple graph of limit groups, the edge groups are asked to be maximal
abelian in both of their adjacent vertex groups, and the morphism $\phi$
is asked to be one-to-one in restriction to \emph{all} vertex groups. 
However, those properties need not be satisfied by the canonical splitting $\Gamma_\can$ and
by the morphism $\sigma$.
The goal of next proposition is to show how those assumptions can be dropped. 
The key result follows immediately.
%\coucou{j'ai cru un instant que les groupes d'aretes de $\Gamma_can$
%etaient toujours max abeliens dans les groupes de sommets. En fait c'est presque le cas
%mais pas tout a fait: si un groupe d'arete $G_e$ de $Gamma_\can$ n'est pas max abelien
%dans $G_v$, alors $G_v$ est abelien (et correspond a un cylindre). Mais, meme
%si on suppose que $G$ n'est pas une extension libre de centralisateurs,
%il se peut que $G_e$ soit un sous-groupe strict (d'indice fini par exemple) de $G_v$}

\begin{dfn}[General graph of limit groups]\label{dfn_ggolg}
  A group $G$ is a \emph{general graph of limit groups over a limit group $L$} if
$G$ is the fundamental group of a graph of groups $\Gamma$ whose vertex groups are finitely generated and such that
\begin{itemize*}
\item $G$ is commutative-transitive;
\item each edge group is a non-trivial abelian group;
\item there is an epimorphism $\phi:G\ra L$ such that 
\begin{itemize*}
\item $\phi$ is one to one in restriction to each edge group;
  \item for each vertex $v\in\Gamma$ of surface type, $\phi(G_v)$ is non-abelian;
  \item for each non surface type vertex $v$, $\phi$ is one-to-one
in restriction to the elliptic abelian neighbourhood $\Hat G_v$ of $G_v$.
\end{itemize*}
\end{itemize*}
\end{dfn}

This proposition gives a general statement that precises
the statements in Definition 5.11 of \cite{Sela_diophantine1}.

\begin{prop}\label{prop_general_to_simple}
Let $G=\pi_1(\Gamma)$ be a general graph of limit groups over $\phi:G\ra L$.

Then $G$ can be written as a \emph{simple} graph of limit groups over the same morphism $\phi$.

In particular, $G$ is a limit group.
\end{prop}

\begin{rem*}
  \begin{enumerate*}
    \item
  If for some abelian vertex group $G_v$, we allow $\phi$ to be one-to-one only in restriction 
to the direct summand of the incident edge groups, then 
it is not true that there exists a sequence of Dehn twists $\tau_i$ 
such that $\phi\circ\tau_i$ converges to $\id_G$ in $\calg(G)$. As a matter of fact, such a Dehn twist
restricts to a conjugation on each abelian vertex group. Sela's proof misses this point.

\item The statements in
definition 5.11 and theorem 5.12 of \cite{Sela_diophantine1} seem to
be slightly incorrect.  A simple counterexample is
%the following:
%take a free abelian group $A$, and $\phi:A\ra A$ an automorphism
%which is not the identity. Then the mapping torus of $\phi$ is not a limit group as it is not CSA.
%Another kind of example is 
the following double~: $G=(_C * S) *_{C=\ol{C}} \overline{(S *_C)}$ 
where $S$ is the fundamental group of a punctured
torus, and $C$ is conjugate to the fundamental group of its boundary component
(see figure below).

\centerline{\fontsize{7pt}{7pt}\selectfont \input contrexemple.pst}

  The fundamental group of this graph of groups is a double of a surface
group with extended centralizer, and it is not commutative transitive
(it contains a subgroup isomorphic to $\bbZ^2 *_\bbZ \bbZ ^2 \simeq F_2 \times \bbZ$). 

To avoid stating technical conditions on the centralizers of edges, 
we include the hypothesis that $G$ is commutative transitive in the
result above. In view of the characterization of CSA graph of groups
given in Corollary \ref{cor_acyl} in appendix \ref{sec_CSA}, one could replace this hypothesis
with the stronger assumption that $\Gamma$ is acylindrical.
  \end{enumerate*}
\end{rem*}

\begin{proof}
There are two main steps to prove this proposition. First, we cut surfaces occuring in $\Gamma$
so that $\phi$ is one-to-one in restriction to the elliptic abelian neighbourhood
of \emph{all} vertex groups of the new graph of groups (Proposition \ref{prop_cut_surf}).
In second step, we \emph{pull centralizers} so that edge groups become maximal abelian
in neighbouring vertex groups (Proposition \ref{prop_pull}). The proposition follows.
\end{proof}

\subsubsection{Step 1: cutting surfaces}

\begin{prop}\label{prop_cut_surf}
  Let $G=\pi_1(\Gamma)$ be a general graph of limit groups over
  $\phi:G\ra L$.
  
  Then one can refine $\Gamma$ into a graph of groups $\Gamma_1$ such
  that
\begin{itemize*}
\item $\Gamma_1$ is a general graph of limit groups over $L$;
\item $\phi$ is one-to-one in restriction to the elliptic abelian
  neighbourhood of \emph{each} vertex group.
\end{itemize*}
\end{prop}

\begin{rem*}
  The proof follows Sela when all the fundamental groups of boundary components
of surface type vertices are maximal abelian in  $G$. 
The general case needs the additional easy Lemma \ref{lem_small_surfaces}.
\end{rem*}

The proof is based on the following elementary lemma of Sela.

  \begin{lem}[{\cite[Lemma 5.13]{Sela_diophantine1}}]\label{lem_ll}
    Let $S$ be the fundamental group of a surface $\Sigma$ (maybe with
    boundary) with Euler characteristic at most $-1$.  Let $\phi : S
    \ra L$ be a morphism to a limit group $L$ with non abelian image,
    and which is one-to-one in restriction to the fundamental groups
    of its boundary components.
    
    Then there exists a family of disjoint simple closed curves
    $c_1,\dots, c_p$ of $\Sigma$, such that $\phi(c_i)$ is non-trivial
    for all $i$, all the connected components of $\Sigma\setminus
    (c_1\cup\dots\cup c_p)$ is either a pair of pants or a punctured
    M\"obius band, and $\phi$ is one to one in restriction to the
    fundamental group of each of these components.
  \end{lem}

  \begin{rem*}
    Note that the fundamental group of a surface of Euler
    characteristic -1 with non-empty boundary is a free group of rank
    2.  Since its image under $\phi$ is a non-abelian limit group,
    $\phi$ is one-to-one in restriction to this fundamental group as
    soon as its image is non abelian (see point \ref{2_generated} in
    Prop.~\ref{elementary_properties}). The idea to prove the lemma is
    to find an essential simple closed curve whose image in $L$ is
    non-trivial, and such that the images in $L$ of the connected
    components of the complement are non-abelian.  Then one iterates
    the procedures on connected components of the complement.
  \end{rem*}

\begin{proof}[Proof of the proposition.]
  Using Lemma \ref{lem_ll}, we refine the graph of groups $\Gamma$
  into a graph of groups $\Gamma_1$ by splitting the surface type
  vertices occuring in $\Gamma$ 
along the simple closed curves given by the lemma.
 
  Call \emph{new vertices} of $\Gamma_1$ all the vertices coming from
  the subdivision of surface type vertices of $\Gamma$, and \emph{old
    vertices} the other ones.  
We want to prove that for each (old or new) vertex group $G_v$,
$\phi$ is one-to-one in restriction to its elliptic abelian neighbourhood $\Hat G_v$.

  The elliptic abelian neighbourhood of an old vertex group in
  $\Gamma_1$ is not larger than in $\Gamma$ since elliptic elements
  in $\Gamma_1$ are elliptic in $\Gamma$.  Thus $\phi$ is still
  one-to-one in restriction to $\Hat G_v$ for each old vertex $v$ of
  $\Gamma_1$.
  
According to Lemma \ref{lem_ll}, we
  also know that $\phi$ is one-to-one in restriction to each
  new vertex group $G_v$ of $\Gamma_1$. 
  This implies that for each edge $e$ of $\Gamma_1$, $\phi$ is
  one-to-one in restriction to $\Hat G_e$: $\Hat G_e$ is elliptic, so it is conjugate
  into some vertex group of $\Gamma_1$.
  
  There remains to prove that for each new vertex group $G_v$ of
  $\Gamma_1$, $\phi$ is one-to-one in restriction to its elliptic
  abelian neighbourhood $\Hat G_v$.
  
  We remark that at least one of the edges $e$ incident on $v$
  corresponds to a subdivision curve (otherwise $v$ would be an old
  vertex).  Therefore, for this edge $e$, one has $\Hat G_e=G_e$.
  Since $G_v$ is the fundamental group of a pair of pants or of a
  punctured M\"obius band, next lemma concludes.
\end{proof}

\begin{lem}[Embedding of abelian neighbourhood of small surfaces]\label{lem_small_surfaces}
  Let $\Sigma$ be a pair of pants or a punctured M\"obius band.  Let
  $G$ be a group containing $S=\pi_1(\Sigma)$.
  Denote by $b_1,b_2,b_3$ (or by $b_1,b_2$ in the case of a punctured M\"obius band) 
  some generators of the fundamental groups of the boundary components of $\Sigma$.
  For each index $i$, consider an abelian group $B_i\subset G$ containing $b_i$ such
that for at least one $i$, one has $B_i=\langle b_i\rangle$. Let $\Hat S$
be the subgroup of $G$ generated by $S$ and the abelian groups $B_i$.
  
  Let $\phi:\Hat S\ra L$ be a morphism to a limit group $L$ which is 
  one-to-one in restriction to each group $B_i$.
 and such that $\phi (\Hat S)$ is non-abelian.
  
  Then $\phi$ is one-to-one in restriction to $\Hat S$.
\end{lem}

\begin{proof}
  We first consider the case of a pair of pants, so that $S$ has a
  presentation of the form $\langle b_1,b_2,b_3 \,|\, b_1b_2b_3=1 \rangle$
where each $b_i$ is a generator of the fundamental group of a boundary components of $\Sigma$. 
Assume for instance that $B_3=\langle b_3\rangle$.
Then $\Hat S$ is generated by $B_1$ and $B_2$ and the following claim concludes.
  
  In the case of a a punctured M\"obius band, one has a presentation of the form
 $S=\langle a,b_1,b_2\,|\, a^2b_1b_2=1 \rangle$. If $B_2=\langle b_2\rangle$, 
then $\Hat S$ is generated by $B_1$ and $\langle c\rangle$, and the following claim also
  concludes.
\end{proof}

\begin{claim}\label{claim_one-to-one}
  Let $A,B$ be two abelian groups, and $L$ be a limit group. If a
  morphism $\phi:A*B\ra L$ has non-abelian image and is one-to-one in
  restriction to $A$ and $B$, then $\phi$ is one-to-one on $A*B$.
\end{claim}

\begin{proof}[Proof of the claim]
  Consider $a_1b_1\dots a_nb_n$ a reduced word in $A*B$.  Let
  $\rho:L\ra F$ be a morphism into a free group such that
  $\rho\circ\phi(a_i)$ and $\rho\circ\phi(b_i)$ are non-trivial and do
  not commute.  Since $\rho\circ\phi (A)$ and $\rho\circ\phi (B)$ are
  abelian in $F$, there exist $\alpha,\beta\in F$ such that the
  elements $\rho\circ\phi(a_i)$ are powers of $\alpha$ and
  $\rho\circ\phi(b_i)$ are powers of $\beta$. Thus $\alpha$ and
  $\beta$ do not commute, and they freely generate a free group.  The
  image under $\rho\circ\phi$ of the word $a_1b_1\dots a_nb_n$ is a
  reduced word in $\langle \alpha,\beta \rangle$ and is thus
  non-trivial.
\end{proof}

\subsubsection{Step 2: pulling centralizers}

\begin{prop}[Pulling centralizers]\label{prop_pull}
  Consider $G=\pi_1(\Gamma)$ be a splitting of a commutative transitive group with
  abelian edge groups.
  
  Then there exists a splitting $G=\pi_1(\Gamma')$ with the same
  underlying graph as $\Gamma$ such that
\begin{itemize*}
\item each edge group of $\Gamma'$ is maximal abelian in the
  neighbouring vertex groups
\item each edge group $G'_e$ of $\Gamma'$ is the elliptic abelian
  neighbourhood $\Hat G_e$ in $\Gamma$ of the corresponding edge group
  $G_e$
\item each vertex group $G'_v$ of $\Gamma'$ is the elliptic abelian
  neighbourhood $\Hat G_v$ in $\Gamma$ of the corresponding edge group
  $G_v$
\end{itemize*}
\end{prop}

The graph of groups $\Gamma'$ will be obtained from $\Gamma$ by a
finite sequence of the following operation:

\begin{dfn}[Pulling centralizers across an edge]
  Consider a graph of groups $\Gamma$ with abelian edge group, and an
  oriented edge $e\in E(\Gamma)$ with $u=o(e)$ and $v=t(e)$.  Consider
  the graph of groups $\Hat \Gamma$ with same underlying graph and
  fundamental group as $\Gamma$ and obtained from $\Gamma$ by the
  following operation: replace $G_e$ by $\Hat G_e$, and replace $G_v$
  by $\Hat G_v$.  The edge morphisms are the natural ones.
  
  The graph of groups $\Hat \Gamma$ is said to be obtained from
  $\Gamma$ by \emph{pulling centralizers across $e$}.
\end{dfn}

\begin{figure}[htbp]
\begin{center}
\includegraphics[width=10cm]{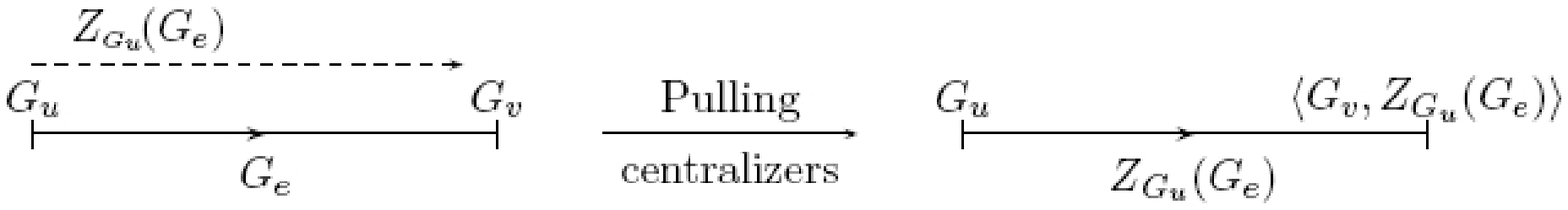}
\caption{Pulling centralizers across an edge}
\label{fig_pulling2}
\end{center}
\end{figure}

\begin{proof}[Proof of the Proposition]
\newcommand{\OK}{full}% si t'as mieux comme terme...
We say that $e$ is \emph{\OK} if $G_e=\Hat G_e$.  If all edges of $\Gamma$
are \OK, then one can take $\Hat\Gamma=\Gamma$ so we argue by
induction on the number of edges which are not \OK.

Let $T$ be the Bass-Serre tree of $\Gamma$.  We already know that for
any edge $e$, $\Hat G_e$ fixes a point in $T$.  Thus, if an edge $e$
is not \OK, let $v_0$ be the vertex fixed by $\Hat G_{e}$ closest to
$e$.  Then $G_{e}$ fixes the arc joining $e$ to $v_0$, and the edge
$e_0$ of this arc incident to $v_0$ satisfies $G_{e_0} \subsetneqq
\Hat G_{e_0}=\Hat G_{e}$, and $\Hat G_{e_0}\subset G_{v_0}$.

Denote by $\Gamma'$ the graph of groups obtained by pulling the
centralizer $\Hat G_{e_0}$ of $G_{e_0}$ across $e_0$ in $\Gamma$.
Clearly, each edge group $G'_e$ (resp. vertex group $G'_v$) in
$\Gamma'$ is contained in the elliptic abelian neighbourhood $\Hat
G_e$ (resp. $\Hat G_v$) of the corresponding group in $\Gamma$.

To conclude we just need to check that $\Gamma'$ has strictly fewer
non-\OK\ edges than $\Gamma$.  Note that pulling centralizers increases
the set of elliptic elements, so one may imagine that $e_0$ might not
become \OK\ or that some other edge $e$ which used to be \OK\ becomes
not \OK\ after the operation. 
We prove that this does not occur by
proving 
that for all edge $e$ of $\Gamma$, the set of
elliptic elements in $Z(G_e)$ does not increase when pulling
centralizers.
Thus, the following claim  will conclude the proof.
\end{proof}

\begin{claim*}
 Consider two non-trivial commuting elements $h,g\in G$ such that
$h$ is hyperbolic and $g$ is elliptic in $T$. 
Then $h$ is still
  hyperbolic in the Bass-Serre tree $T'$ of $\Gamma'$.
\end{claim*}

\begin{proof}
The operation of pulling centralizers might be seen at the level of
Bass-Serre tree as follows: $T'$ is the quotient of $T$ under the
smallest equivariant equivalence relation $\sim$ such that $e_0\sim
g_0.e_0$ for all $g_0\in \Hat G_{e_0}$. More precisely, two edges
$\eps_1,\eps_2$ are folded together if and only if there exists $k\in
G$ such that $k.\eps_1=e_0$ and $k.\eps_2=g_0.e$ for some $g_0\in \Hat
G_{e_0}$.  

Now consider two non-trivial commuting elements $h,g\in G$ such that
$h$ is hyperbolic and $g$ is elliptic in $T$.  If $h$ is elliptic in
$T'$, then there are two distinct edges $\eps_1,\eps_2$ in the axis of $h$
which are folded together.  Up to conjugating $h$, we can assume that
$\eps_1=e_0$ and $\eps_2=g_0.e_0=g_0.\eps_1$ for some $g_0\in \Hat G_{e_0}$.  On
the other hand, since $[g,h]=1$, $g$ fixes pointwise the axis of $h$
in $T$, so $g\in G_{e_0}$. By commutative transitivity, the element
$g_0\in \Hat G_{e_0}$ also commutes with $h$. Since $g_0$ is elliptic,
$g_0$ fixes the axis of $h$, and thus fixes $\eps_1$ and $\eps_2$.
This contradicts the fact that $g_0.\eps_1=\eps_2$.
\end{proof}

This terminates the proof of Proposition \ref{prop_general_to_simple} showing
that a general graph of limit groups over $L$ can be turned into a simple graph of limit groups.
Thus, the key result and the two characterizations of limit groups follow.

\subsection{Fully residually free towers}\label{sec_frft}

As a corollary of the fact that general graph of limit groups over limit groups
are limit groups, one way to construct limit
groups will be by \emph{gluing retracting surfaces} (Proposition
\ref{prop_grs}). Together with the extension of centralizers, 
this construction is the building block
for Sela's \emph{fully residually free towers} (Definition \ref{dfn_frft}).
In topological terms,
consider a space $X$ obtained by gluing a surface $\Sigma$ onto 
a space $L$ whose fundamental group is a limit group
by attaching the boundary components of $\Sigma$ to non-trivial loops of $L$.
If there is a retraction of $X$ onto $L$ which sends
$\Sigma$ to a subspace of $L$ with non-abelian fundamental group,
then $\pi_1(X)$ is a limit group.

\begin{figure}[htbp]
\begin{center}
\includegraphics[width=6cm]{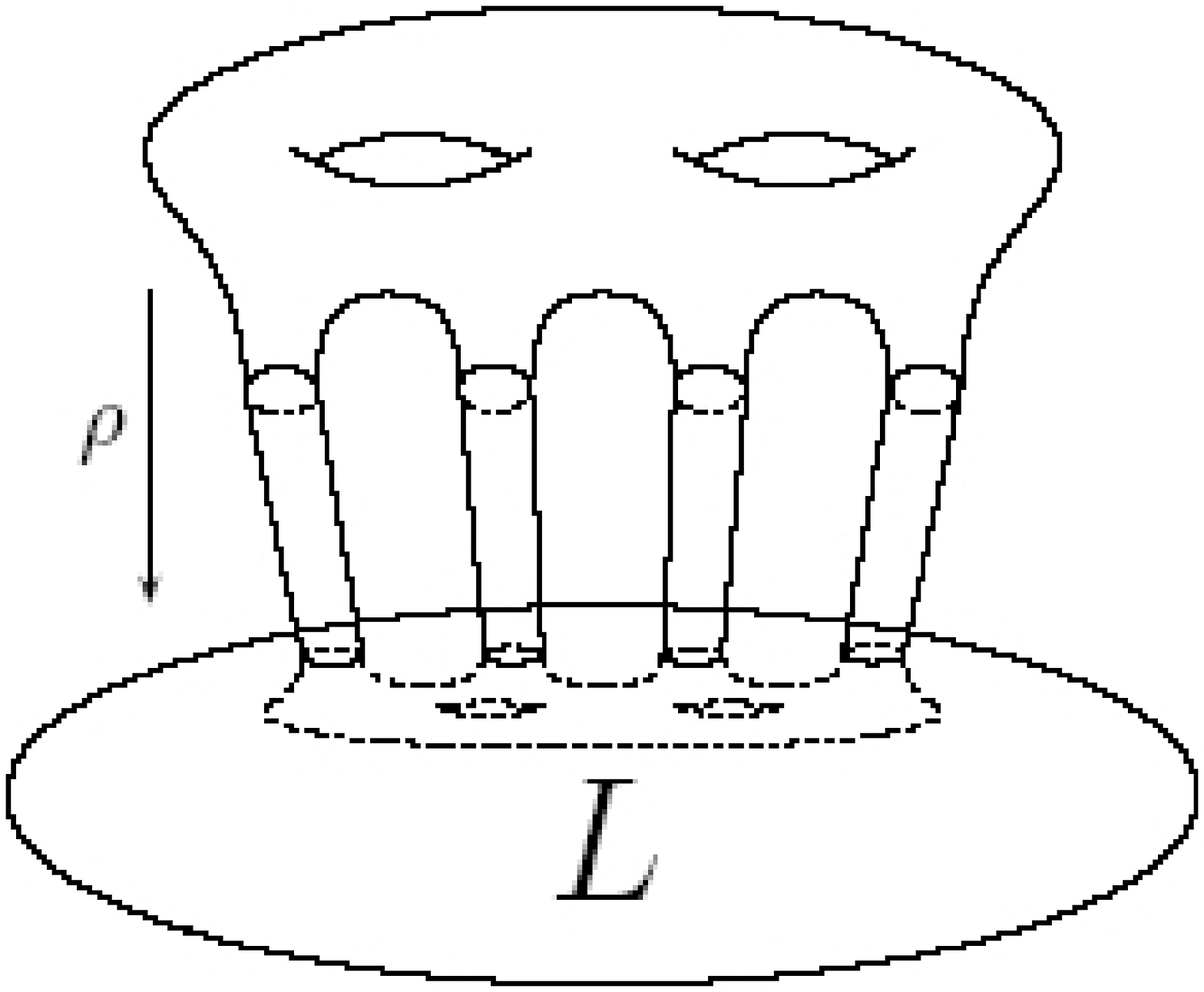}
\caption{Gluing retracting surfaces}
\label{fig_retracting}
\end{center}
\end{figure}

\begin{prop}[Gluing retracting surfaces]\label{prop_grs}
Let $L$ be a limit group, and $\Sigma$ be a surface with boundary with
Euler characteristic at most $-2$ or a punctured torus or a punctured Klein bottle.
%\coucou{il la faut la bouteille de Kleine? Pourquoi pas d'autres surfaces ?
%les seules autres surfaces a bord sont: - le disque, l'anneau, le pantalon,
%le ruban de mobius, et le ruban de mobius troue.
%Le disque l'anneau et le mobius n'ont pas de quotient non-abelien: ils ne peuvent
%satisfaire les hypotheses. Le pantalon et le ruban de mobius troue peuvent satisfaire les hypotheses,
%mais dans ce cas, ils se plongent dans L, ... et alors ? pourquoi les interdire ? 
%ce ne sont certes pas des sommets de Type Surface au sens de la def donnee au dessus, mais
%il faudrait comprendre la preuve de la suite pour savoir pourquoi on les exclut.}
Consider a morphism $\rho:\pi_1(\Sigma)\ra L$ with non-abelian image
which is one-to-one in restriction to the fundamental groups $C_1,\dots, C_n$ of the
boundary components of $\Sigma$.

Consider the graph of groups $\Gamma$ with two vertex groups $L$ and
$\pi_1(\Sigma )$, and
$n$ edge groups $C_1,\dots, C_n$, the two edge morphisms being
the identity and the restriction of $\rho$.

Then $\pi_1(\Gamma)$ is a limit group.
\end{prop}

\begin{rem*}
In fully residually free towers, the groups $\rho(C_i)$ will be asked to
be maximal cyclic in $L$ (see definition \ref{dfn_frft}).
\end{rem*}

\begin{proof}
 Consider a retraction $\phi:\pi_1(\Gamma)\ra L$ such that $\phi$ restricts to the identity on $L$, and 
coincides with $\rho$ on $\pi_1(\Sigma)$. 
Clearly, the centralizer of each edge group is contained in the vertex group $L$.
In particular, the elliptic abelian neighbourhood $\Hat L$ of $L$ coincides with $L$.
Thus, to show that Proposition \ref{prop_general_to_simple} applies, 
we only need to check that $\pi_1(\Gamma)$ is commutative transitive.
This is for instance a consequence of corollary \ref{cor_acyl} given in appendix \ref{sec_CSA}
since $\Gamma$ is clearly $2$-acylindrical since the fundamental group of a boundary component of $\Sigma$
is malnormal in $\pi_1(\Sigma)$.
\end{proof}

We are now ready to give Sela's definition of fully residually free towers.

%We recall here the examples given by Sela.
\begin{dfn}[Fully residually free towers {\cite[Def. 6.1]{Sela_diophantine1}}]\label{dfn_frft}
A finitely generated group is a \emph{fully residually free tower}
if it belongs to the smallest class of groups  $\calt$ containing all the finitely generated
free groups and surface groups and stable under 
free products, 
free extension of centralizers, 
and gluing of retracting surfaces on maximal cyclic subgroups 
(see Proposition \ref{prop_ec} and Corollary \ref{prop_grs}).
\end{dfn}

\begin{SauveCompteurs}{cpteur_frftowers}
\begin{thm}[\cite{Sela_diophantine1}]
Fully residually free towers are fully residually free.
\end{thm}
\end{SauveCompteurs}

If we never extend centralizers, this construction gives only hyperbolic
limit groups.
Sela announces in his sequence of six papers
the following answer to Tarski's problem. We state it as a conjecture
since the referring process is not yet completed.

\begin{conj*}[{\cite[Th.7]{Sela_diophantine6}}]
  A finitely generated group is elementary equivalent to a non
abelian free group if and only if it is a non elementary hyperbolic
fully residually free tower.
\end{conj*}

\section{Some logic}\label{sec_logic}

The goal of this section is to give an intuitive feeling of basic
logical notions, in order to state the following result: a finitely
generated group is a limit group if and only if it has the same
\emph{universal theory} as free groups.  This is by no means a
substitute to a real introduction to model theory.  For more precise
information, see for instance
\cite{Chang_Keisler_model,Hodges_shorter}.  As a motivation, we first
present a short introduction to the \emph{elementary theory} of a
group, and the Tarski problem.

\subsection{The Tarski Problem}\label{sec_tarski}

Z.~Sela and Kharlampovich-Myasnikov have recently announced a positive
solution to this problem (\cite{Sela_diophantine1},
\cite{KhMy_irreducible1}, \cite{KhMy_irreducible2}).  It is in the
solution to this problem that Z.~Sela introduced limit groups.

\begin{pb*}[Tarski, 1945]
  Do finitely generated non-abelian free groups have the same
  elementary theory?
\end{pb*}

The language of groups uses the following symbols:
\begin{enumerate*}
\item The binary function group multiplication ``$.$'', the unary
  function inverse ``$\m$'', the constant ``$1$'', and the equality
  relation ``$=$'',
\item variables $x_1,\dots,x_n,\dots$ (which will have to be
  interpreted as individual group elements),
\item logical connectives ``$\et$'' (meaning \emph{and}), ``$\ou$''
  (meaning \emph{or}), $\non$ (meaning \emph{not}), the quantifiers
  $\forall$ (\emph{for all}) and $\exists$ (\emph{there exists}), and
  parentheses: ``('' ``)''.
\end{enumerate*}

In the language of groups, the \emph{terms} (used to be interpreted as
elements of the group) are words in the variables, their inverses, and
the identity element.  For instance, $((x_1.1).(x_2\m))$ will be
interpreted as an element of the group.  Of course, since we will be
working in a group, we may rather rewrite this term as $x_1x_2\m$,
dropping parentheses for convenience.  To make a \emph{formula}
(interpreted as true or false) from terms, one can first compare two
terms using ``=''.  For instance, $x_1x_2\m=x_3$ is a formula (called
an \emph{atomic} formula).  Then, one can use logical connectors and
quantifiers to make a new formula from other formulae.  For instance,
$\forall x_1 ((x_1x_3=1) \et (\exists x_2\ x_1x_2=x_3))$ is a formula.
Note that this formula has a free variable $x_3$: when interpreted,
the fact that it is true or false will depend on $x_3$.  A formula
with no free variable is called a \emph{sentence}, like for instance
$\forall x_3 \forall x_1 ((x_1x_3=1) \ou (\exists x_2\ x_1x_2=x_3))$.
Given a group $G$ and a sentence $\sigma$ (with no free variable), we
will say that $G$ \emph{satisfies} $\sigma$ if $\sigma$ is true if we
interpret $\sigma$ in $G$ (in the usual sense).  We denote by $G\satis
\sigma$ the fact that $G$ satisfies $\sigma$.  For instance,
$G\satis\forall x_1,x_2 [x_1,x_2]=1$ if and only if $G$ is abelian.

%Given a formula $\phi(x_1,\dots,x_p)$ with free variables
%$x_1,\dots,x_p$,
%one need to specify some values $a_1,\dots,a_p\in G$
%in a group $G$ to interpret it: we will denote by
%$G\satis_{(x_1=a_1,\dots,x_p=a_p)}\phi(x_1,\dots,x_p)$
%the fact that $\phi(a_1,\dots,a_p)$ holds in $G$.
%For instance, $G\satis_{(x_1=a_1)} \forall x_2 [x_1,x_2]=1$
%means that $a_1$ lies in the center of $G$.

Note that the following statement $\forall x_1 \exists k\in\bbN\ x_1^k
= 1$ is not allowed, because it quantifies over an integer, and not a
group element; similarly, one cannot quantify over subsets, subgroups
or morphisms.  In fact the quantifier does not mention to which group
the variables belong, it is the interpretation which specifies the
group.  The \emph{elementary theory} of $G$, denoted by $\Elem(G)$ is
the set of sentences which are satisfied by $G$.

The Tarski problem is a special aspect of the more general problem to
know which properties of a group can be read from its elementary
theory.  Two groups are said \emph{elementarily equivalent} if they
have the same elementary theory.  For instance, the fact of being
abelian, can be expressed in one sentence, and can therefore be read
from the elementary theory of a group.  There is no sentence saying
that a group is torsion free.  The sentence $\forall x (x\neq 1 \implq
x.x\neq 1)$ says that a group has no 2-torsion. Similarly, for any
integer $k$, there is a sentence saying that a group has no
$k$-torsion. Thus, the property of being torsion free can be read in
an infinite set of sentences, so a group elementarily equivalent to a
torsion-free group is torsion-free.

\begin{example*}
  $\bbZ$ and $\bbZ ^2$ don't have the same elementary theory.  Indeed
  one can encode in a sentence the fact there are at most 2 elements
  modulo the doubles: the sentence
  $$\forall x_1,x_2,x_3\ \exists x_4\ 
  (x_1=x_2+2x_4)\ou(x_1=x_3+2x_4)\ou(x_2=x_3+2x_4)$$
  holds in $\bbZ$
  and not in $\bbZ^2$ since $\bbZ/2\bbZ\leq 2$ has only two elements,
  and $\bbZ^2/2\bbZ^2$ has four elements.
\end{example*}

\subsection{Universal theory}
\label{sec_universal}

A \emph{universal formula} is a formula which can be written $\forall
x_1\dots\forall x_p\ \phi(x_1,\dots ,x_p)$ for some quantifier free
formula $\phi(x_1,\dots,x_p)$.
%Note that $\forall x_1 \non (\forall x_2 [x_1,x_2]=1)$
%is \emph{not} a universal formula because of the symbol
%$\non$ occuring before the second $\forall$.
If it has no free variables, a universal formula is called a universal
\emph{sentence}.  The \emph{universal theory} of a group is the set of
universal sentences satisfied by $G$.  Similarly, one can define
\emph{existential} formulae and sentences, and the existential theory
of a group.  Note that two groups which have the same universal theory
also have the same existential theory since the negation of a
universal sentence is equivalent to an existential statement.

Note that if $H<G$, a quantifier free formula which is satisfied for
every tuple of elements of $G$ is also satisfied for every tuple of
elements of $H$, hence $\Univ(G)\subset \Univ(H)$.  Similarly, if $H<
G$, $\Exist(H)\subset \Exist(G)$.  As a corollary, for $n\geq 3$,
since $F_n$ contains $F_2$ and $F_2$ contains $F_n$, one gets
$\Univ(F_2)=\Univ(F_n)$.

\begin{thm}[\cite{Gaglione_Spellmann_even_more,Chi_introduction,Remeslennikov_exist_siberian}]
\label{limit_universal}
A finitely generated group $G$ has the same universal theory as a
non-abelian free group (resp. an infinite cyclic group) if and only if
$G$ is a non-abelian limit group (resp. an abelian limit group, \ie a
free abelian group).
\end{thm}

This theorem was stated by Remeslennikov in the context of
non-standard free group (see section \ref{sec_non-standard}).  For
instance, this theorem means that a non-exceptional surface group has
the same universal theory as $F_2$.  It also means that $\bbZ$ and
$\bbZ ^2$ have the same universal theory.

\begin{proof}
  Let's first prove that a non-abelian limit group $\Gamma$ has the
  same universal theory as $F_2$.  Since a non-abelian limit group
  contains $F_2$ (prop.~\ref{elementary_properties}), one has
  $\Univ(\Gamma)\subset \Univ(F_2)$.  There remains to check that
  $\Univ(\Gamma)\supset\Univ(F_2)$.  This follows from the following
  proposition:

\begin{SauveCompteurs}{cpteursUniv}
\begin{prop}
\label{universal=closed}
Let $\sigma$ be a universal sentence. Then the property
$G\satis\sigma$ is closed in $\calg_n$.  Equivalently, if a sequence
of marked groups $(G_i,S_i)$ converge to a marked group $(G,S)$, then
$\Univ (G)\supset \limsup\Univ (G_i).$
\end{prop}
\end{SauveCompteurs}

%Equivalently, for any existential sentence $\sigma$, the property $G\satis \sigma$ is open in $\calg_n$.

\begin{proof}[Proof of the proposition]
  We will prove that for any existential sentence $\sigma$, the
  property $G\satis \sigma$ is open.
  
  Consider the sentence $\exists x_1,\dots, x_p\ \phi(x_1,\dots,x_p)$
  where $\phi(x_1,\dots,x_p)$ is quantifier free. Using distributivity
  of $\et$ with respect to $\ou$, one easily checks that
  $\phi(x_1,\dots,x_p)$ is equivalent to a formula
  $\Sigma_1(x_1,\dots,x_p)\ou\dots\ou\Sigma_q(x_1,\dots,x_p)$, where
  each $\Sigma_i$ is a system of equations or inequations in the
  following sense: a set of equations or inequations of the form
  $w(x_1,\dots,x_p)=1$ or $w(x_1,\dots,x_p)\neq 1$ separated by the
  symbol ``$\et$'' (\emph{and}) where $w(x_1,\dots ,x_p)$ is a word on
  $x_1\pmu ,\dots,x_p\pmu$.
  
  Consider a marked group $(G,S)\in \calg_n$ with $S=(s_1,\dots,s_n)$,
  and consider an existential sentence $\sigma$ of the form
  $$\exists x_1,\dots, x_p\ 
  \Sigma_1(x_1,\dots,x_p)\ou\dots\ou\Sigma_q(x_1,\dots,x_p)$$
  such
  that $G\satis\sigma$. So consider $a_1,\dots,a_p\in G$ and
  $i\in\{1,\dots,q\}$ such that $\Sigma_i(a_1,\dots,a_p)$ holds.
  Consider $R$ large enough so that the ball of radius $R$ in $(G,S)$
  contains $\{a_1,\dots,a_p\}$ and so that for each word $w$ occuring
  in $\Sigma_i$, the corresponding word on $\{a_1,\dots,a_p\}$ can be
  read in this ball (for instance, one can take $R$ to be the maximal
  length of the words times the maximal length of the $a_i$'s in
  $(G,S)$).  Now assume that a marked group $(H,S')$ has the same ball
  of radius $R$ as $(G,S)$. Clearly, this implies that the
  corresponding elements $a'_1,\dots, a'_p$ in the ball of $(H,S')$
  satisfy $\Sigma_i$, so that $H\satis \sigma$.
\end{proof}

We now prove that if $G$ has the same existential theory as $F_2$,
then it is a non-abelian limit group.  It is clearly non-abelian since
the property of being non-abelian expresses as an existential
sentence.  Let $S=(s_1,\dots,s_n)$ be a finite generating family of
$G$, let $R>0$, and let $B$ be the ball of radius $R$ of $(G,S)$. We
aim to find a generating set of a free group having the same ball.
For this purpose, we are going to encode the ball in a system of
equations and inequations.  Let $w_1,\dots,w_p$ be an enumeration of
all the words on $x_1\pmu,\dots,x_n\pmu$ of length at most $R$.  For
each of these words $w_i$, consider $g_i=w_i(s_1,\dots,s_n)\in B$.  We
consider the following system $\Sigma(x_1,\dots,x_p)$ of equations and
inequations: for each $i,j\in\{1,\dots p\}$, we add to $\Sigma$ the
equation $w_i=w_j$ or the inequation $w_i\neq w_j$ according to the
fact that $g_i=g_j$ or $g_i\neq g_j$. Of course,
$\Sigma(s_1,\dots,s_p)$ holds. Thus
$$G\satis \exists x_1,\dots,x_p\ \Sigma(x_1,\dots,x_p).$$
Since $F_2$
has the same existential theory, let $s'_1,\dots,s'_p\in F_2$ such
that $\Sigma(s'_1,\dots,s'_p)$ holds in $G$. Let $F=\langle
s'_1,\dots,s'_p\rangle < F_2$.  Thus $F$ is a free group, and
$(F,(s'_1,\dots,s'_p))$ has the same ball of radius $R$ as $(G,S)$.

Finally, the abelian part of the theorem states that all the finitely
generated free abelian groups have the same universal theory.  The
same proofs as above work in the abelian context using the facts that
$\bbZ^p\supset \bbZ$ and that
$[\bbZ^p]_{\calg_n}\subset\ol{[\bbZ]_{\calg_n}}$.
\end{proof}

Note that the second part of the proof actually shows the following
more general statement:
\begin{prop}\label{prop_univ2}
  If $\Univ(G)\supset\Univ(H)$, then for all generating family $S$ of
  $G$, $(G,S)$ is a limit of marked subgroups of $H$.
\end{prop}

This kind of result is usually stated using ultra-products. The next
section details the relation between convergence in $\calg_n$ and
ultra-products.

\section{A little non-standard analysis}
\label{sec_non-standard}

\subsection{definitions}
\begin{dfn}[ultrafilter]
  An \emph{ultrafilter} on $\bbN$ is a finitely additive measure of
  total mass $1$ (a \emph{mean}) defined on all subsets of $\bbN$, and
  with values in $\{0,1\}$. In other words, it is a map
  $\omega:\calp(\bbN)\ra\{0,1\}$ such that for all subsets $A,B$ such
  that $A\cap B=\es$, $\omega(A\cup B)=\omega(A)+\omega(B)$,
  $\omega(\bbN)=1$.
  
  An ultrafilter is \emph{non-principal} if it is not a Dirac mass,
  \ie if finite sets have mass $0$.
\end{dfn}

We will say that a property $P(k)$ depending on $k\in\bbN$ is true
\emph{$\omega$-almost everywhere} if $\omega(\{k\in\bbN|P(k)\})=1$.
Note that a property which is not true almost everywhere is false
almost everywhere.
%Also note that the ultrafilter is non-principal if and only if finite sets have measure $0$ 
%(and cofinite sets have measure $1$).
Given an ultrafilter $\omega$ (which will usually be supposed to be
non-principal), and a family of groups $(G_k)_{k\in\bbN}$, there is a
natural equivalence relation $\eqv_\omega$ on $\prod_{k\in \bbN} G_k$
defined by equality $\omega$-almost everywhere.  When there is no risk
of confusion, we may drop the reference to the ultrafilter $\omega$.

\begin{dfn}[ultraproduct, ultrapower]
  The ultraproduct with respect to $\omega$ of a sequence of groups
  $G_k$ is the group
  $$\left(\prod_{k\in\bbN} G_k\right)\Big/\eqv_\omega.$$
  When starting
  with a constant sequence $G_k=G$, the ultraproduct is called an
  \emph{ultrapower}, and it is often denoted by $\ultra{G}$ (though
  depending on the ultrafilter $\omega$).
\end{dfn}

The main interest of ultraproducts and ultrapowers is {\L}os Theorem,
which claims that ultrapowers of a group $G$ have the same elementary
theory as $G$ (see for instance \cite{Bell_Slomson_models}).

\begin{thm}[{\L}os] 
  Let $G$ be a group, and $\ultra{G}$ an ultrapower of $G$.  Then $G$ and
  $\ultra{G}$ have the same elementary theory.
  
  More generally, for every formula $\phi(x_1,\dots,x_n)$,
  $\ultra{G}\satis_{(x_1=(a_{1,k})_{k\in\bbN},\dots,x_n=(a_{n,k})_{k\in\bbN})}$
  if and only if for almost every $k\in \bbN$,
  $G\satis_{(x_1=a_{1,k},\dots,x_n=a_{n,k})}$.
\end{thm}

\subsection{Ultraproducts and the topology on the set of marked
  groups}
\label{ultraproduct_topology}

The link between ultraproducts and convergence of groups in $\calg_n$
is contained in the following lemma.

\begin{SauveCompteurs}{cpteursUltra}
\begin{prop}\label{prop_ultra_topology}
  \begin{enumerate}
  \item\label{item_one} Consider a sequence of marked groups
    $(G_k,S_k)\in\calg_n$ which accumulates on $(G,S)\in\calg_n$.
    Then there exists some non-principal ultrafilter $\omega$ such
    that $G$ embeds in the ultraproduct $\prod_{k\in\bbN}
    G_k/\eqv_\omega$.
  \item\label{item_any} If the sequence above is convergent, then
    $(G,S)$ embeds in \emph{any} ultraproduct $\prod_{k\in\bbN}
    G_k/\eqv_\omega$ (assuming only that $\omega$ non-principal).
  \item\label{item_soussuite} Consider a finitely generated subgroup
    $(G,S)$ of an ultraproduct $\prod_{k\in\bbN} G_k/\eqv_\omega$
    (where $\omega$ is non-principal).  Then there exists a
    subsequence $(i_k)_{k\in\bbN}$, and marked subgroups
    $(H_{i_k},S_{i_k})<G_{i_k}$, such that $(H_{i_k},S_{i_k})$
    converge to $(G,S)$.
  \end{enumerate}
\end{prop}
\end{SauveCompteurs}

\begin{proof}[Proof of the lemma.]
  \ref{item_one}.  Denote $S_k=(s_1^{(k)},\dots,s_n^{(k)}$. Assume
  that there is a subsequence $(G_{i_k},S_{i_k})$ converging to
  $(G,S)$.  Consider an ultrafilter $\omega$ such that the subsequence
  $\{i_k|k\in\bbN\}$ has full $\omega$-measure (the existence of
  $\omega$, which uses the axiom of choice, is proved in
  \cite[p.39]{Bourbaki_topologie}).  Let $U=\prod_{k\in\bbN}
  G_k/\eqv_\omega$ be the corresponding ultraproduct, and let
  $\ol{S}=(\ol{s_1},\dots,\ol{s_n})$ the family of elements of $U$
  defined by $\ol{s_p}=(s_p^{(k)})_{k\in\bbN}\in U$ for
  $p\in\{1,\dots,n\}$.  Let $\ol{G}$ be the subgroup of $U$ generated
  by $\ol{S}$.  We prove that $(G,S)$ is isomorphic to
  $(\ol{G},\ol{S})$ as a marked group.  For any word $w$ on the
  generators $s_1,\dots,s_n$ and their inverses, if $w$ is trivial in
  $(G,S)$ then it is trivial in $(G_{i_k},S_{i_k})$ for all but
  finitely many $k$'s, and hence for $\omega$-almost every $k\in\bbN$,
  thus $w$ is trivial in $(\ol{G},\ol{S})$.  If $w$ is non-trivial in
  $(G,S)$ then it is non-trivial in $(G_{i_k},S_{i_k})$ for all but
  finitely many $k$'s, and hence for $\omega$-almost every $k\in\bbN$,
  thus $w$ is non-trivial in $(\ol{G},\ol{S})$.  This proves that
  $(G,S)\simeq(\ol{G},\ol{S})$.
  
  \ref{item_any}. If the sequence above is convergent, then taking
  $i_k=k$ in the argument above, allows to choose any non-principal
  ultrafilter.
  
  \ref{item_soussuite}. Let $S=(s_1,\dots,s_n)$, where each $s_p$ is
  an element of $U$.  So write each generator $s_p$ as a sequence
  $(s_p^{(k)})_{k\in\bbN}$. Let $S_k=(s_1^{(k)},\dots,s_n^{(k)})$, and
  let $H_k$ be the subgroup of $G_k$ generated by $S_k$.  We prove
  that $(H_k,S_k)$ accumulates on $(G,S)$.  Consider a ball of radius
  $R$ in $(G,S)$ and consider a word $w$ on $s_1,\dots,s_n$ of length
  at most $R$. The word $w$ is trivial in $(G,S)$ if and only if it is
  trivial in $(H_k,S_k)$ for $\omega$-almost every $k$. Thus the set
  of indices $k$ such that the ball of radius $R$ of $(G_k,S_k)$
  coincides with the ball of radius $R$ of $(G,S)$ has full measure as
  an intersection of finitely many full measure subsets. This set of
  indices is therefore infinite since $\omega$ is non-principal, thus
  $(G_k,S_k)$ accumulates on $(G,S)$.
\end{proof}

The following corollary is immediate:

\begin{cor}
\label{limit_ultra}
A group is a limit group if and only if it is a finitely generated
subgroup of an ultraproduct of free groups, and any ultraproduct of
free groups contains all the limit groups.
\end{cor}

\subsection{Application to residual freeness}\label{sec_residual}

We review here the following result of Remeslennikov which proves that
limit groups are residually free (\cite{Remeslennikov_exist_siberian},
see also \cite[lem.5.5.7]{Chi_book}).

\begin{prop}[Remeslennikov]
\label{ultra_residually}
A finitely generated subgroup of an ultrapower $\ultra{F_2}$ is fully
residually free.
\end{prop}

Using corollary \ref{limit_ultra}, this result gives an elementary
proof (without use finite presentation of limit groups) that limit
groups are fully residually free.

\begin{proof}
  Fix an ultrafilter $\omega$, the corresponding ultrapower $\ultra{F_2}$ of
  $F_2$, and $G<\ultra{F_2}$ a finitely generated subgroup.  It is well
  known that for any odd prime $p$, the kernel of $\phi:SL_2(\bbZ)\ra
  SL_2(\bbZ/p\bbZ)$ is a non-abelian free group.  Thus $F_2$ embeds in
  $\ker \phi\subset SL_2(\bbZ)$.  Therefore, $\ultra{F_2}$ embeds in the
  kernel of the natural morphism $\ultra{\phi}:SL_2(\ultra{\bbZ})\ra
  SL_2(\ultra{(\bbZ/p\bbZ)})$ where $\ultra{\bbZ}$ is the ring obtained by taking
  the $\omega$-ultrapower of the ring $\bbZ$.  In particular $\ultra{F_2}$
  embeds in $SL_2(\ultra{\bbZ})$.  $G$ being finitely generated, $G$ embeds
  in $SL_2(R)$ for a ring $R$ which is finitely generated subring of
  $\ultra{\bbZ}$.

\begin{lem}[Remeslennikov]
  Consider a finitely generated subring $R$ of $\ultra{\bbZ}$.  Then $R$, as
  a ring, is fully residually $\bbZ$, \ie for $a_1,\dots,a_k\in
  R\setminus\{0\}$, there exists a ring morphism $\rho:R\ra\bbZ$ such
  that $\rho(a_i)\neq 0$ for all $i=1,\dots,k$.
\end{lem}

Let's conclude the proof of the proposition using the lemma.  Consider
finitely many elements $g_1,\dots,g_k\in G\setminus\{1\}\subset
SL_2(R)$, and let $a_1,\dots,a_{k'}$ be the set of non-zero
coefficients of the matrices $g_j-\Id$ ($j=1,\dots,k$) (there is at
least one non-zero coefficient for each $g_j$ since $g_j\neq \Id$).
Consider a morphism $\rho:R\ra\bbZ$ given by the lemma.  The induced
morphism $\psi:SL_2(R)\ra SL_2(\bbZ)$ maps the elements $g_j$ to
non-trivial elements, and since $G\subset\ker \ultra{\phi}$,
$\psi(G)\subset\ker\phi$ which is free.  Thus $\psi(G)$ is free and
$G$ is fully residually free.
\end{proof}

\begin{proof}[Proof of the lemma]
  Let $t_1,\dots,t_n\in \ultra{\bbZ}$ such that $R=\bbZ[t_1,\dots,t_n]$.
  Consider the corresponding exact sequence $J\hookrightarrow
  \bbZ[T_1,\dots,T_n]\onto R$, where $\bbZ[T_1,\dots,T_n]$ is the ring
  of polynomials with $n$ commuting indeterminates.  Since
  $\bbZ[T_1,\dots,T_n]$ is Noetherian, the ideal $J$ is generated by
  finitely many polynomials $f_1,\dots,f_q$.  Let $a_1,\dots,a_k\in
  R\setminus\{0\}$, and let $g_1,\dots,g_k$ some preimages in
  $\bbZ[T_1,\dots,T_n]\setminus J$.  Note that $(t_1,\dots,t_n)$ is a
  solution of the system of equations and inequations
  $$\left\{\begin{array}{rl}
      f_i(x_1,\dots,x_n)=0 &(i=1,\dots,q) \\
      g_j(x_1,\dots,x_n)\neq 0 &(j=1,\dots,k)
\end{array}\right.$$
Now one can invoke {\L}os theorem, or just remember that each $t_i$ is
a sequence of integers modulo the ultrafilter $\omega$ to check that
almost all the components of $t_i$ provide a solution
$(x_1,\dots,x_n)$ to this system in $\bbZ$.  The morphism
$\bbZ[T_1,\dots,T_n]\ra \bbZ$ sending $T_i$ to $x_i$ induces the
desired morphism $\rho:R\ra\bbZ$.
\end{proof}

\subsection{Maximal limit quotients}\label{sec_finitude}

The result of this section is lemma 1.1 of \cite{Razborov_systems} and lemma 3 of
\cite{KhMy_irreducible2}. It also appears in \cite{BMR_algebraicI}.
 We give a short proof inspired by \cite{Chatzidakis_limit}.
 It is a clever way to get the existence
of maximal limit quotients without using the finite presentation of
limit groups.  It just uses the fact that limit groups are residually
free.

\begin{rem*}
  One could avoid using the fact that limit groups are residually free
  by replacing the mentions of $\emph{residually free}$ groups by
  groups that are \emph{residually a limit group}, and by replacing
  the field $\bbC$ by an ultrapower $\ultra{\bbC}$ to ensure that limit
  groups embed in $SL_2(\ultra{\bbC})$ in lemma \ref{prop_decroissance2}.
\end{rem*}

\begin{prop}\label{prop_decroissance2}
  Consider a sequence of quotients of finitely generated groups
  $$G_1\onto G_2\onto\dots\onto G_k\onto\dots$$
  If every group $G_i$
  is residually free, then all but finitely many epimorphisms are
  isomorphisms.
\end{prop}

\begin{proof}
  Take $S_1=(s_1,\dots,s_n)$ a finite generating family of $G_1$, and
  let $S_k=(s_1^{(k)},\dots,s_n^{(k)})$ its image in $G_k$ under the
  quotient map.  Let $V_k\subset SL_2(\bbC)^n$ be the variety of
  representations of $(G_k,S_k)$ in $SL_2(\bbC)$, \ie
\begin{multline*}V_k=\Big\{ 
  (M_1,\dots,M_n)\in SL_2(\bbC)^n\ |\\
  \forall \text{ relation $r$ of $(G_k,S_k)$, } r(M_1,\dots,M_n)=\Id
  \in SL_2(\bbC) \Big\}\end{multline*}

Note that $V_k$ is an affine algebraic variety in $(\bbC)^{4n}$, and
that $V_1\supset V_2\supset \dots \supset V_k\dots$. By noetherianity,
for all but finitely many indices $k$, one has $V_k=V_{k+1}$.  There
remains to check that if $G_{k+1}$ is a strict quotient of $G_k$, then
$V_{k+1}$ is strictly contained in $V_k$.  So consider a word $r$ on
$S\pmu$ which is trivial in $G_{k+1}$ but not in $G_k$.

Since $G_k$ is residually free, there exists a morphism $\phi:G_k\to
F_2$ such that $\phi(r)\neq 1$.  Since $F_2$ embeds in $SL_2(\bbC)$,
there exists a representation $\rho:G_k\ra SL_2(\bbC)$ such that
$\rho(r)\neq 1$.  This representation provides a point in
$V_k\setminus V_{k+1}$.
\end{proof}

Remember that any group $G$ has a largest residually free quotient
$RF(G)$: $RF(G)$ is the quotient of $G$ by the intersection of the
kernels of all morphisms from $G$ to free groups.  The following
corollary says that a residually free group is presented by finitely
many relations plus all the relations necessary to make it residually
free.

\begin{cor}[\cite{Chatzidakis_limit}]\label{lem_finitely_presented_modulo_RF}
  If $(G,S)$ is residually free, then there exists finitely $S$-words
  $r_1,\dots,r_p$ such that $G=RF(H)$ where $H$ is the group presented
  by $\langle S;r_1,\dots,r_p\rangle$.
\end{cor}

\begin{proof}
  Enumerate the relations $r_i$ of $(G,S)$, and take $G_k=RF(\langle
  S; r_1,\dots,r_k\rangle)$.  The previous lemma says that for $k$
  large enough, $(G,S)=(G_k,S)$.
\end{proof}

\begin{cor}
  Given a residually free marked group $(G,S)$, there is a
  neighbourhood $V_{(G,S)}$ of $(G,S)$ such that every residually free
  group in $V_{(G,S)}$ is a quotient of $(G,S)$.
\end{cor}

\begin{proof}
  Take $r_1,\dots r_p$ relations of $(G,S)$ as in the corollary above.
  One can take $V_{(G,S)}$ to be the set of marked groups $(G',S')$
  such that the relations $r_1,\dots,r_p$ hold in $(G',S')$.
\end{proof}

We can now give an elementary proof (without using finite presentation
of limit groups) of Proposition \ref{prop_maximal}:

\begin{cor}
  Let $G$ be a finitely generated group.  Then there exist finitely
  many quotients $\Gamma_1,\dots,\Gamma_k$ of $G$, such that each
  $\Gamma_i$ is a limit group and such that any morphism from $G$ to a
  free group factors through one $\Gamma_i$.
\end{cor}

\begin{proof}
  Consider a marking $(G,S)$ of $G$, and consider the set
  $K\subset\calg_n$ of marked quotients of $(G,S)$ which are limit
  groups.  This is clearly a compact subset of $\calg_n$, and each of
  its points is residually free.  Now cover $K$ by finitely many
  $V_{(\Gamma_i,S_i)}$ as in the previous corollary to get the result.
\end{proof}

\begin{rem*}
  The proof given in \cite{Chatzidakis_limit} is more logical in nature, and it is
  quite appealing.  It relies on the following ideas.  Let
  $r_1,\dots,r_p$ be some relations of $(G,S)$.  Say that an $S$-word
  $r'$ is \emph{deducible from $r_1,\dots,r_p$ modulo $\Univ(F_2)$} if
  the statement
  $$\forall s_1,\dots, s_n,\ \left\{\begin{matrix}r_1=1\\ \vdots \\ 
      r_p=1\end{matrix}\right. \rightarrow r'=1$$
  holds in $F_2$.  Say
  that a marked group $(G,S)$ is \emph{closed under deduction mod
    $\Univ(F_2)$} if for all relations $r_1,\dots,r_p$ of $(G,S)$, any
  word $r'$ which is deducible from $r_1,\dots,r_p$ modulo
  $\Univ(F_2)$ is a relation of $(G,S)$.  One can easily show that a
  group is closed under deduction mod $\Univ(F_2)$ if and only if it
  is residually a limit group (and hence if and only if it is
  residually free).  Then the statement of lemma
  \ref{lem_finitely_presented_modulo_RF} says that a group $(G,S)$
  which is closed under deduction mod $\Univ(F_2)$ is finitely
  presented mod $\Univ(F_2)$.
\end{rem*}

\appendix
\section{Reading property CSA from a graph of groups.}\label{sec_CSA}

%Remember that property CSA (Definition \ref{dfn_CSA})
%means that maximal abelian subgroups are malnormal.
%A group is CSA, then it is commutative transitive.
This appendix explains how to read CSA property (Definition
\ref{dfn_CSA}) on a graph of groups $\Gamma$ with abelian edge groups.
 
Let us start with two basic cases.  In an amalgamated product $G=A*_C
B$ with abelian edge group $C$, if $C$ is not maximal in $A$ nor in
$B$, then $G=A*_C B$ is not commutative transitive.  Thus for example
$_\bbZ * \bbZ *_\bbZ \bbZ * _\bbZ \simeq \bbZ ^2 *_\bbZ \bbZ ^2$ (with
obvious embeddings) is not a limit group, since it is not commutative
transitive.  Similarly, the amalgam $\bbZ^2 *_Z \bbZ *_Z \bbZ ^2$ is
not commutative transitive and is therefore not a limit group.

The second basic case concerns the HNN extension $\langle A,t|
tat\m=\phi(a)\rangle$ for an injective endomorphism $\phi:A\ra A$.
This group may be commutative transitive for non trivial $\phi$ (for
example the Baumslag-Solitar groups are commutative transitive), but
such an extension is CSA if and only if $\phi$ is the identity.  We
are going to prove that those two basic phenomena are the
only obstructions for getting a CSA group (see Cor.~\ref{cor_CSA}).\\

Let $\Gamma$ be a graph of group ant let $T$ be the Bass-Serre tree of
$\Gamma$.  Consider the equivalence relation on the set of
non-oriented edges of $T$ generated by $e\sim e'$ if $e$ and $e'$ have
a common endpoint and their stabilizers commute.  We call a
\emph{cylinder} the subtree of $T$ defined by an equivalence class.

Before stating a CSA criterion for a graph of groups, we first give
the following necessary condition.

\begin{lem}
\label{lem_neccsa}
Let $\Gamma$ be a graph of groups with non-trivial abelian edge
groups, ant let $T$ be its Bass-Serre tree.

If $\pi_1(\Gamma)$ is CSA, then the vertex groups of $\Gamma$ are CSA
and the global stabilizer of every cylinder in $T$ is abelian.
\end{lem}

\begin{proof}
  Since CSA property is stable under taking subgroups, if
  $G=\pi_1(\Gamma)$ is CSA then vertex groups are CSA.  Consider a
  cylinder $C$ in $T$.  By commutative transitivity of $G$, the
  stabilizers of all the edges of $C$ commute. Since edge groups are
  non-trivial, there exist a unique maximal abelian subgroup $A$ of
  $G$ containing these edge stabilizers.  Let $g\in G$ such that
  $g.C=C$. Then $gAg\m=A$, so $g\in A$ because $G$ is CSA.
\end{proof}

Before proving the converse of this lemma, we first introduce the
\emph{tree of cylinders} of $\Gamma$, and give some properties of this
tree. A similar construction is used in \cite{GL_automorphismes} to
get splittings invariant under automorphisms.

A trivial but fundamental property of cylinders is that two distinct
cylinders intersect in at most one vertex (they have no edge in
common).  Moreover, the set of cylinders is $G$-invariant.  Therefore,
there is a natural bipartite $G$-tree $T_\calc$ called \emph{the tree
  of cylinders} defined as follows: $V(T_\calc)=V_1(T_\calc)\dunion
V_0(T_\calc)$ where $V_1(T_\calc)$ is the set of cylinders of $T$, and
$V_0(T_\calc)$ is the of vertices of $T$ belonging to at least two
distinct cylinders, and there is an edge $\eps=(x,C)$ between $x\in
V_0(T_\calc)$ and $C\in V_1(T_\calc)$ if $x\in C$. The fact that this
graph is a tree is straightforward (see the notion of
transverse coverings in \cite{Gui_limit}).

\begin{dfn}\label{dfn_acyl}
  An action of a group $G$ on a tree is \emph{$k$-acylindrical} if for
  all $g\in G\setminus\{1\}$, the set of fix points of $g$ has
  diameter at most $k$. Similarly, a graph of groups is
  $k$-acylindrical if the action of $\pi_1(\Gamma)$ on its Bass-Serre
  tree is $k$-acylindrical. We say that $\Gamma$ is
  \emph{acylindrical} if it is $k$-acylindrical for some $k$.
\end{dfn}

We gather a few simple facts about the tree of cylinders.

\begin{lem}
\label{lem_treecyl}
Let $\Gamma$ be a graph of groups with CSA vertex groups and abelian
edge groups. Suppose that the global stabilizers of cylinders of
$\Gamma$ are abelian.

Then:

\begin{enumerate*}
\item the stabilizer of each vertex $x\in V_0(T_\calc)$ is CSA
\item the stabilizer of each vertex $C\in V_1(T_\calc)$ is
  abelian
\item \label{cl_max} The stabilizer of any cylinder intersect a vertex
  group in maximal abelian subgroup of this vertex group, that is for
  all edge $\eps\in E(T_\calc)$ incident on $x\in V_0(T_\calc)$,
  $G_\eps$ is maximal abelian in $G_x$
\item \label{cl_acyl} if $\eps,\eps'\in E(T_\calc)$ are such that
  $G_\eps\cap G_{\eps'}\neq\{1\}$ then $\eps$ and $\eps'$ have a
  common endpoint in $V_1(T_\calc)$ (in particular, $T_\calc$ is
  $2$-acylindrical)
\item \label{cl_ab} for every abelian subgroup $A\subset G$, either
  $A$ fixes a point in $T_\calc$, or $A$ is a cyclic group acting
  freely on $T_\calc$.
  \end{enumerate*}
\end{lem}

\begin{proof}
  The first two claims result from the definitions.  For claim
  \ref{cl_max}, consider an edge $\eps=(x,C)\in E(T_\calc)$.  Assume
  that $g\in G_x$ commutes with $G_\eps$.  Consider an edge $e\in C$
  incident on $x$.  Since $G_e\subset G_\eps$, $g$ commutes with
  $G_e$, so $G_e=G_{g.e}$ and $g.e\in C$. This implies that $g.C=C$,
  and that $g$ fixes the edge $\eps=(x,C)$.
  
  For claim \ref{cl_acyl}, assume that $\eps=(x,C),\eps'=(x',C')\in
  E(T_\calc)$ are such that $G_\eps\cap G_{\eps'}\neq \{ 1\} $.  We
  want to prove that $C=C'$.  If $x=x'$, claim \ref{cl_max} states
  that $G_{\eps}$ and $G_{\eps'}$ are maximal abelian in $G_x$. Since
  $G_\eps\cap G_{\eps'}\neq \{ 1\}$, $G_\eps = G_{\eps'}$. Let $e$ and
  $e'$ be any edges of $C$ and $C'$ adjacent to $x$. We have $G_e
  \subset G_\eps$ and $G_{e'} \subset G_{\eps'}$, thus $G_e$ and
  $G_{e'}$ commute. This proves $C=C'$.  If $x\neq x'$, then any non
  trivial element $h\in G_\eps\cap G_{\eps'}$ fixes $[x,x']$. Thus
  $[x,x']$ is contained in a cylinder $C''$.  Let $\eps '' = (x,C'')$.
  Then $h\in G_\eps \cap G_{\eps''}$, and the previous case shows that
  $C=C''$. Similarly, $C'=C''$.

  For claim \ref{cl_ab}, consider an abelian group $A$.  Suppose $A$
  contains a hyperbolic element $h$.  For any element $g\in
  A\setminus\{1\}$, let $\Fix (g)$ be its set of fix points in
  $T_\calc$.  Since $\Fix (g)$ is $h$-invariant and bounded by
  $2$-acylindricity, one has $\Fix g=\es$. This means that the action
  of $A$ is free, so $A$ is a cyclic group.  Suppose $A$ contains an
  elliptic element $h\neq 1$.  Let $F$ be the set of fix points of $h$
  in $T_\calc$.  Acylindricity shows that $F$ is bounded. Since it is
  $A$-invariant, $A$ fixes a point in $F$.
\end{proof}

\begin{prop}[a CSA criterion]
  Consider a graph of groups $\Gamma$ with torsion-free vertex groups
  and non-trivial abelian edge groups.
  
  Then $\pi_1(\Gamma)$ is CSA if and only if vertex groups of $\Gamma$
  are CSA and the global stabilizer of every cylinder in the
  Bass-Serre tree of $\Gamma$ is abelian.
\end{prop}

\begin{rem*}
  One could replace the assumption that vertex groups are torsion free
  by the assumption that $\pi_1(\Gamma)$ contains no infinite dihedral
  subgroup acting faithfully on $T$. This assumption is more natural
  since the infinite dihedral group is not CSA.
\end{rem*}

\begin{proof}
  Lemma \ref{lem_neccsa} shows one part of the equivalence.  We now
  assume that vertex stabilizers are CSA and that the stabilizer of
  each cylinder is abelian, and we have to prove that $G$ is CSA.  Let
  $A$ be a maximal abelian subgroup of $G$ and assume that $gAg\m\cap
  A\neq \{1\}$.  Suppose first that $A$ acts freely on $T_\calc$.
  Denote by $l$ its axis and let $G_l$ be the global stabilizer of
  $l$.  Clearly, $g\in G_l$ so we are reduced to prove that $G_l$ is
  abelian, since it will follow that $A=G_l\ni g$.  Because of the
  acylindricity of $T_\calc$, no element of $G_l$ fixes $l$ so $G_l$
  acts faithfully by isometries on $l$, and $G_l$ is either cyclic or
  dihedral. But since vertex groups are torsion free (this is the only
  place where we use this assumption), $G_l$ cannot be dihedral.
  
  If $A$ fixes a point in $T_\calc$, let $F$ be its set of fix points
  in $T_\calc$.  If $g$ fixes a point in $F$, then we are done since
  $g$ and $A$ both belong to a vertex stabilizer of $T_\calc$ which is
  known to be CSA.  Let $h\in A\cap gAg\m\setminus\{1\}$. Since $h$
  fixes pointwise $F\cup g.F$, and thus its convex hull, $F\cup g.F$
  is contained in the $1$-neighbourhood of a vertex $v_1\in
  V_1(T_\calc)$ (claim \ref{cl_acyl} of the fact).  Since $v_1$ is the
  only vertex of this neighbourhood lying in $V_1(T_\calc)$ ($T_\calc$
  is bipartite), if $v_1\in F$ , then $v_1$ is fixed by both $g$ and
  $A$, and this case was already settled.  Thus one can assume that
  $v_1\notin F$ so $F$ consists in a single vertex $v_0\in
  V_0(T_\calc)$.  One can also assume that $g.v_0\neq v_0$. Note that
  in this case, $v_1$ is the midpoint of $[v_0,g.v_0]$.  Let $\eps$ be
  the edge joining $v_1$ to $v_0$. One has $h\in G_\eps\subset
  G_{v_0}$, and $h\in A\subset G_{v_0}$.  Since $G_\eps$ is maximal
  abelian in $G_{v_0}$ (claim \ref{cl_ab} in the fact), one has
  $G_\eps=A$. Thus $A$ fixes $\eps$, so $A$ fixes $v_1$, so $v_1\in F$
  which has been excluded.
\end{proof}

We now translate our criterion into a more down-to-earth
property (compare \cite[Definition 5.11]{Sela_diophantine1})

\begin{dfn}[Cylinders]
  Let $\Gamma$ be a graph of groups with abelian edge groups and CSA
  vertex groups.  Denote by $Cyl(\Gamma)$ the following
  (non-connected) graph of groups.
  
  Edges of $Cyl(\Gamma)$ are the edges of $\Gamma$, and they hold the
  same edge groups. We define the vertices of $Cyl(\Gamma)$ by
  describing when two oriented edges have the same terminal vertex:
  $e$ and $e'$ have the same terminal vertex in $Cyl(\Gamma)$ if they
  have the same terminal vertex $v$ in $\Gamma$ and if there exists
  $g\in G_v$ such that $i_e(G_e)$ and $g.i_{e'}(G_{e'}).g\m$ commute.
  The corresponding vertex group is the maximal abelian group
  containing $i_e(G_e)$ (which is well defined up to conjugacy) and
  the edge morphism is $i_e:G_e\ra i_e(G_e)$, which is well defined
  since any conjugation preserving the maximal abelian group $A\subset
  G_v$ containing $G_e$ fixes $A$ because $G_v$ is CSA.
\end{dfn}

The connected components of $Cyl(\Gamma)$ correspond to the orbits of
cylinders in $T_\calc$ in the following sense:

\begin{lem}
  Let $C$ be a cylinder of $T$, $G_C$ its global stabilizer. Consider
  the graph of groups $\Lambda=C/G_C$.  Then $\Lambda$ corresponds to
  a connected component of $Cyl(\Gamma)$.  The fundamental group of a
  connected component of $Cyl(\Gamma)$ is conjugate to the stabilizer
  of the corresponding cylinder in $T_\calc$. In particular they are
  maximal abelian subgroups of $G$.
\end{lem}

\begin{proof}
  The proof is straightforward and left as an exercise.
\end{proof}

\begin{cor}\label{cor_CSA}
  Consider a graph of groups $\Gamma$ with torsion-free vertex groups
  and non-trivial abelian edge groups.
  
  Then $\pi_1(\Gamma)$ is CSA if and only vertex groups of $\Gamma$
  are CSA and each connected component $\Lambda$ of $Cyl(\Gamma)$ is
  of one of the following form:
\begin{itemize}
\item either $\Lambda$ is a trivial splitting: there is a vertex
  $v_0\in\Lambda$ such that the injection
  $G_{v_0}\subset\pi_1(\Lambda)$ is actually an isomorphism.  This
  translates into the fact that $\Lambda$ is a tree of groups and that
  for all vertex $v\neq v_0$ and for all edge $e$ with $t(e)=v$ and
  separating $v$ from $v_0$, the edge morphism $i_e:G_e\ra G_{t(e)}$
  is an isomorphism.
\item or $\Lambda$ has the homotopy type of a circle
  $\calc\subset\Lambda$ (as a simple graph), the edge morphisms of
  edges of $\calc$ are isomorphisms, the composition of all the edge
  morphisms around $\calc$ is the identity, and the injection
  $\pi_1(\calc)\subset\pi_1(\Lambda)$ is actually an isomorphism.
  This last property translates into the fact that for all edge
  $v\notin\calc$ and for all edge with $t(e)=v$ and separating $v$
  from $\calc$, the edge morphism $i_e:G_e\ra G_{t(e)}$ is an
  isomorphism.
\end{itemize}
\end{cor}

\begin{proof}
  It is clear that each of the cases implies that $\pi_1(\Lambda)$ is
  abelian since vertex groups of $\Lambda$ are abelian.
  
  Conversely, if $\pi_1(\Lambda)$ is abelian, then the fundamental
  group of the graph underlying $\Lambda$ has to be abelian, so
  $\Lambda$ has the homotopy type of a point or of a circle.
  
  In the first case, $\pi_1(\Lambda)$ is generated by finitely many
  vertex stabilizers.  But it is an easy exercise to check that two
  commuting vertex stabilizers must fix a common vertex, so whole
  group $\pi_1(\Lambda)$ must fix a vertex.
  
  If $\Lambda$ has the homotopy type of a circle $\calc\subset
  \Lambda$, then the action of $\pi_1(\Lambda)$ on its Bass-Serre tree
  is non-trivial, but since $\pi_1(\Lambda)$ is abelian, its minimal
  invariant subtree is a line $l$, and $\pi_1(\Lambda)$ acts by
  translations on $l$. Moreover, $\calc=l/\pi_1(\Lambda)$ is such that
  the injection $\pi_1(\calc)\subset\pi_1(\Lambda)$ is an isomorphism,
  and the edge morphisms of edges of $\calc$ are isomorphisms, and
  that the composition of all the edge morphisms around $\calc$ is the
  identity (as it is induced by a conjugation in the abelian group
  $\pi_1(\Lambda)$).
\end{proof}

\begin{cor}\label{cor_acyl}
  Consider a graph of groups $\Gamma$ with CSA vertex groups and
  abelian edge groups.  If every edge group is maximal abelian in the
  two neighbouring vertex groups, and if the Bass-Serre tree of
  $\Gamma$ is acylindrical, then $\pi_1(\Gamma)$ is CSA.
\end{cor}

\begin{proof}
%  Every edge morphism in the graph of cylinders is an isomorphism.
  Acylindricity of $\Gamma$ implies that each component $\Lambda$ of
  $Cyl(\Gamma)$ is a tree of groups.  Moreover, every edge morphism in
  $\Lambda$ is an isomorphism, therefore the splitting corresponding
  to $\Lambda$ is trivial.  Thus the first condition of the previous
  corollary holds.
\end{proof}

\bibliographystyle{alpha} \bibliography{published,unpublished}

\noindent Christophe Champetier\\
Institut Fourier\\
B.P. 74\\
38402 Saint-Martin d'H{\`e}res\\
FRANCE\\
\emph{e-mail:}\texttt{Christophe.Champetier@ujf-grenoble.fr}\\

\noindent Vincent Guirardel\\
Laboratoire E. Picard, UMR 5580\\
B\^at 1R2, Universit\'e Paul Sabatier\\
118 rte de Narbonne,\\
31062 Toulouse cedex 4\\
FRANCE\\
\emph{e-mail:}\texttt{guirardel@picard.ups-tlse.fr}\\
\end{document}